\documentclass[11pt]{article}

\usepackage{authblk}        
\usepackage{hyperref}       
\usepackage{geometry}       
\title{\bfseries An asymptotic model of Poisson--Nernst--Planck--Stokes systems \\ in narrow channels}

\author[1]{Christine Keller \thanks{Email: \texttt{christine.keller@wias-berlin.de}, ORCID: 0009-0002-6579-5900}}
\author[2]{Andreas M\"unch \thanks{Email: \texttt{andreas.muench@maths.ox.ac.uk},            ORCID: 0000-0002-8325-3809}}
\author[1]{Barbara Wagner \thanks{Email: \texttt{barbara.wagner@wias-berlin.de},     ORCID: 0000-0001-8306-3645}}

\affil[1]{Weierstrass Institute \\ Anton-Wilhelm-Amo-Str. 39 \\ 10117 Berlin \\ Germany}
\affil[2]{Mathematical Institute\\ University of Oxford\\Woodstock Road \\  Oxford, OX2 6GG, UK}

\date{\today}

\usepackage{setspace}

\usepackage{graphicx} 

\usepackage{tikz}
\usepackage{xfp} 

\usepackage{float}

\usepackage{xcolor}
\usepackage{mathrsfs}
\usepackage{amsmath}
\usepackage{amssymb}
\usepackage{graphics}
\usepackage{subcaption}
\usepackage{adjustbox} 

\usepackage{dsfont}
\usepackage{comment}
\usepackage{multirow}
\usepackage{enumitem}
\usepackage{cite}
\usepackage[siunitx]{circuitikz}
\usepackage{amsfonts}
\usepackage{amsthm}
\graphicspath{{img/}}
\usepackage{svg}
\usepackage{array}
\newcolumntype{M}[1]{>{\raggedright\arraybackslash}m{#1}}
\usepackage{longtable}
\usepackage{mathtools}
\usepackage[textsize=tiny]{todonotes}

\usepackage{algorithm}
\usepackage{algpseudocode}

\usepackage{listings}
\definecolor{codegreen}{rgb}{0,0.6,0}
\definecolor{codegray}{rgb}{0.5,0.5,0.5}
\definecolor{codepurple}{rgb}{0.58,0,0.82}
\definecolor{backcolour}{rgb}{0.95,0.95,0.92}
\lstdefinestyle{mystyle}{
    backgroundcolor=\color{backcolour},   
    commentstyle=\color{codegreen},
    keywordstyle=\color{magenta},
    numberstyle=\tiny\color{codegray},
    stringstyle=\color{codepurple},
    basicstyle=\ttfamily\footnotesize,
    breakatwhitespace=false,         
    breaklines=true,                 
    captionpos=b,                    
    keepspaces=true,                 
    numbers=left,                    
    numbersep=5pt,                  
    showspaces=false,                
    showstringspaces=false,
    showtabs=false,                  
    tabsize=2}
\usepackage{hyperref}
\usepackage{cleveref}
\usepackage{myHeader}
\usepackage{tikz}
\usetikzlibrary{shapes.geometric, arrows}
\tikzstyle{arrow} = [thick,->,>=stealth]

\tikzstyle{startstop} = [rectangle, rounded corners, 
minimum width=3cm, 
minimum height=1cm,
text centered, 
text width=4cm,
draw=black]

\usepackage{pgfplots}
\numberwithin{equation}{section}

\begin{document}
\maketitle
\begin{abstract}
Ion transport through narrow channels is described by the coupled Poisson--Nernst--Planck--Stokes equations (PNPS) on a continuum scale. 
However, direct numerical simulations in two or three dimensions of  boundary value problems for small aspect ratio geometries, a crucial characteristic of nanopores, can quickly become computationally intensive and thus limit the insights into the underlying mechanisms that control electrokinetic phenomena. 
Taking advantage of the small aspect ratio, we derive a systematic asymptotic reduction of the PNPS system.
In contrast to existing one-dimensional reductions, which assume a Debye length much smaller than the channel radius, our analysis identifies a 
distinguished asymptotic regime in which the Debye length is allowed to be comparable to the channel width. Our approach has a significantly larger range of validity and contains existing  approximations such as the Helmholtz-Smoluchowski approximation as limiting cases. 
The derived asymptotic model extends also to a generalized PNPS system, where finite-size constraints and solvation effects are taken into account and thus applies to  other well-known models such as the Bikerman-Freise model.
Using our asymptotic model we  demonstrate that the ion current can undergo a number of different flow transitions and in particular predict that positively charged ions can be pushed against their electrostatic gradient. Furthermore, we show how finite-size effects can influence the ion current and enhance ion selectivity. Finally, we revisit case studies of protein-based channels from the literature to illustrate the predictive potential of our asymptotic model.
\end{abstract}


\newpage
\section{Introduction}
Nanopores are nanoscale channels that play an important role in both biology and technical applications. Biological nanopores, also referred to as ion channels, are pore-forming proteins in cell membranes and are involved in many physiological functions. They regulate the exchange of ions across the otherwise impermeable lipid bilayer and are able to discriminate between different molecules based on size, charge and structure. Synthetic nanopores on the other hand are made of solid-state materials such as silicon, nitride, graphene, glass or polymers and highly customizable. They are used in a variety of applications such as protein and DNA sequencing \cite{branton2008potential,Saharia2021,Wang2021}, water filtration and desalination \cite{CohenTanugi2012,Koehler2018, Surwade2015} and many more. 

The transport in ion channels is controlled by a highly complex interplay of forces, where the relative magnitude of each contribution is set by the specific geometry of the channels, electrolyte composition, concentration gradients, the pressure exerted on the channel, and surface conditions such as charge patterns on the channel surface. 
An external electric field therefore exerts a direct force on the molecule (electrophoresis) while simultaneously pulling the surrounding fluid (electroosmosis).  By adjusting salt concentration, pore diameter, or the chemistry of the pore surface one can deliberately shift the dominant driving force from electrophoresis to electroosmotic flow (EOF), or vice‑versa, and even to pressure‑driven hydrodynamic flow \cite{CohenTanugi2012,Koehler2018,Saharia2021}.
These electro‑hydrodynamic phenomena occur naturally in biological ion channels as well as in synthetic nanopores  \cite{Gillespie2008, Li2024, Willems2020}.
Geometric asymmetries, for example, produce ion‑current rectification, i.e., asymmetries in the measured current-voltage curves \cite{siwy2004conical, harrell2004dna, dal2019confinement}.   
Surface charge patterns in channels can also qualitatively alter electroosmotic flows due to the coupling between hydrodynamic flow and electrostatic interactions \cite{stroock2000patterning, bocquet2010nanofluidics, jubin2018dramatic, Green2022}.
A quantitative understanding of how to tune the interplay of electrophoresis, electroosmosis, and hydrodynamic flow is therefore essential for optimizing molecular translocation in sequencing, sensing, and the design of new bio‑inspired micro‑ and nanofluidic devices.

Current-voltage (IV) relations provide important information about ion transport in nanopores. They can provide information about conductivity or the effects of EOF under varying conditions such as different ion concentrations, surface charges and geometries.
The interpretation of IV curves is therefore crucial for many applications and can help in the development and optimization of devices. 

Mathematical modeling can be used to quantify how ion concentration gradients, applied voltage and surface charge influence the flow of ions through the pore. To model ion electro-diffusion an extended form of the Poisson-Nernst-Planck (PNP) equations, including finite-size and solvation effects, are commonly used. This set of equations combines the Poisson equation which describes the electrostatic potential with the Nernst-Planck equations, describing the movement of ions due to diffusion and drift due to electrostatic forces. Coupling PNP with Stokes equations allows to capture the movement of the solvent (e.g. water) due to the electric field \cite{Curk2024,Dreyer2013, jubin2018dramatic,Green2021,Green2022}. 

Solutions of the corresponding boundary value problems of such generalized Poisson-Nernst-Planck-Stokes (PNPS) type system are typically pursued numerically in two or three dimensions and often in a axially symmetric setting. However, as shown in Li and Muthukumar \cite{Li2024} the shape of a protein channel for example can be quite irregular and at the same time the specific geometry has a significant impact on the signaling in sequencing. Thus, even on the continuum level, parameter studies can be computationally time consuming.

Asymptotic analysis is a valuable tool for identifying the predominant physical mechanisms governing the system in specific regimes, such as high or low ion concentration, strong electric fields, confined geometries, etc..
The simplification of the full model can facilitate more profound analysis of the system's behavior as a function of specific parameters (e.g. ion concentration, electric potential, nanopore size, etc.). Consequently, significant insights into the limiting behavior of ion transport and electrostatics in nanopores under various conditions can be obtained. Thereby facilitating the derivation of scaling laws, the prediction of current-voltage characteristics, and the optimization of nanopore designs.

A variety of approaches exist for the reduction of model complexity. 
Ion channels typically have a small aspect ratio and this fact has been used to reduce the PNPS system by introducing area-averaged 1D models for axially symmetric cases, see for example \cite{rankin2016effect, fair1971reverse}. Similar ideas had also been pursued for the simpler PNP system by Nonner and Eisenberg \cite{Nonner1998} (see also \cite{Gillespie1999}). 
The underlying assumptions are that the variations of the concentrations and the electrostatic potential in the transverse and radial directions are negligible. Consequently, the averaged concentrations and potentials are calculated over the cross-sectional area of the channel. However, already for the PNP system strong radial variations of ion concentrations due to surface charges, for example, can not be captured with that approach as could be shown in Matejczyk et al. \cite{matejczyk2018}. In fact, for the PNP model they showed that the assumption that the influence of the surface charge on the ion concentration and the voltage can be averaged over the cross section of the pore is only valid if the Debye length of the electrolyte is much larger than the pore width. In their analysis they present an asymptotic solution to the model for ion channels with radius comparable to the Debye length of the electrolyte. This asymptotic solution allows the characterization of the behavior of the ion channel in terms of the solution to a 1D model. By comparison to numerical solutions in their case studies they demonstrate that the 1D area averaged PNP equations do not provide a good approximation to the full two-dimensional equations except for extremely dilute electrolytes or very narrow pores so that the Debye length of the electrolyte is much larger than the radius of the channel. 
Matejczyk \cite{Matejczyk2019} also derived an asymptotic reduction for the generalized form of the PNP equations, using a slightly different scaling and a slightly different approach to the asymptotics as in \cite{matejczyk2018}.
Inspired by this studies and the success of the long-wave approximation in thin liquid films \cite{craster2009dynamics} or in shallow water waves \cite{whitham2011linear}, we derive an asymptotic model for the axially symmetric boundary value problem for a generalized PNPS system in the distinguished limit where the Debye length is comparable to the pore radius. 
It is in particular the limit where the charge selectivity is most pronounced \cite{schoch2008transport}.

In the following sections we first present the boundary value problem for the generalized PNPS system in section 2, followed by the scaled version for axisymmetric geometries in the asymptotic limit of large-aspect ratio in section 3. 
We then systematically derive the asymptotic model for the generalized PNPS system for large-aspect ratio pore geometries in section 4. 

Apart from showing convergence of the full two-dimensional axially symmetric problem to the quasi-1D asymptotic model in section 5, we apply our model for a number parameter studies that predict universal transport behavior across a range of potential landscapes that point to underlying physical mechanisms controlling the flow. In addition, we show that the ion current can undergo different flow transitions and predict that positively charged ions can be pushed against their electrostatic gradient. Furthermore, by comparing different variants of the model, we show that finite-size effects can influence the ion current and enhance ion selectivity, particularly in narrow domains. 
Finally, we revisit case studies of protein-based channels from the literature to demonstrate the predictive potential of our asymptotic model. 

\section{The generalized Poisson-Nernst-Planck-Stokes system}
\begin{figure}[H]
    \centering\
   \begin{tikzpicture}[scale=0.5]

    \draw[->,line width=0.25mm] (-0.5,0) -- (27,0) node[right] {$z$}; 
    \draw[->,line width=0.25mm] (0,-3.5) -- (0,3.5) node[above] {$r$}; 
    
    \draw[line width =0.5mm,domain=0:25,samples=100] plot (\x, {\fpeval{2 + (0.8 - 2)*exp(-(\x - 6)*(\x - 6) / (2*1.5*1.5)) + (0.2 - 0.8)/(1 + exp(-(\x - 9.5)/0.4)) + (1.7 - 1.0)/(1 + exp(-(\x - 15.5)/0.5))}});
    
    \draw[->,line width=0.25mm] (11,3.5) -- (9.5,2) 
    node[right] at (10,4) {pore wall $S^\twall$};

    \draw[->,line width=0.25mm] (11,3.5) -- (12.5,-1.2);

    \draw[line width =0.5mm,domain=0:25,samples=100] plot (\x, {\fpeval{-(2 + (0.8 - 2)*exp(-(\x - 6)*(\x - 6) / (2*1.5*1.5)) + (0.2 - 0.8)/(1 + exp(-(\x - 9.5)/0.4)) + (1.7 - 1.0)/(1 + exp(-(\x - 15.5)/0.5)))}});

    \draw[line width =0.5mm] (0,2) -- (0,-2) 
    node[left] at (0,1){outlet $S^\tout$}; 
    
    \draw[line width =0.5mm,dashed] (0,0) -- (25,0);

    \draw[->,line width=0.25mm] (20,3.3) -- (20,0.2)
    node[right] at (17,4) {symmetry axis $S^0$};
    
    \draw[line width =0.5mm] (25,-2.1) -- (25,2.1) 
    node[right] at (25,1) {inlet $S^\tin$};

    
\end{tikzpicture}
	\caption{
    Sketch of an axially symmetric profile of an ion channel $\Omega$. The pore wall at $\tilde r= \tilde R(\tilde z)$ is denoted with $S^\twall$, the symmetry axis of the domain at $\tilde r=0$ is marked as $S^0$. At $\tilde z=0$ the pore has an outlet and at $\tilde z= \tilde L$ an inlet boundary.
    }
	\label{fig:Domain}
\end{figure}
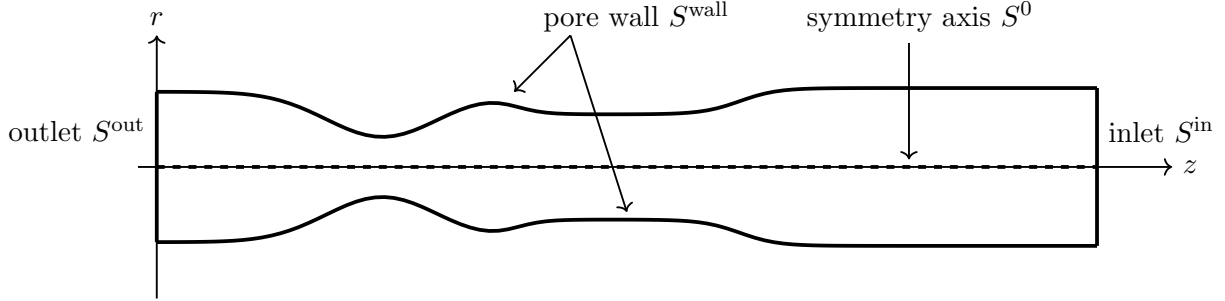

We consider a mixture of anions, cations, solvent molecules and a viscous incompressible fluid. The mixture contains $N$ different ion species with molar density $\tilde n_\alpha(\tilde t,\tilde \vecx)$ for $\alpha = 1, \dots, N$. The ions have molar masses $\tilde m_\alpha$, molar volumes $\tilde v_\alpha$ and carry a charge $z_\alpha e_0$, where $e_0$ is the elementary charge and $z_\alpha$ is the charge number.  
Within this work the solvent density is denoted by $\tilde n_0(\tilde t,\tilde \vecx)$, if present, and is considered to be electrically neutral with a charge number of $z_0=0$.
As an example, consider a mixture of calcium (Ca$^{2+}$), sodium (Na$^+$), chloride (Cl$^-$), and water. The system consists of $N=3$ ionic species with concentrations $n_1$, $n_2$, and $n_3$ corresponding to calcium, sodium, and chloride, and charge numbers $z_1=2$, $z_2=1$, and $z_3=-1$, respectively. The water concentration is denoted by $n_0$.
The ions might in addition be subject to solvation effects, and we denote by $\tilde m_\alpha$ and  $\tilde v_\alpha$ the mass and volume of the solvated ions, respectively \cite{Landstorfer2018}. Hence, the molar mass writes as $\tilde m_\alpha = \hat{m}_\alpha + \kappa_\alpha \tilde  m_0$ since mass is conserved upon solvation, where $ \hat{m}_\alpha$ is the mass of the central ion, $\kappa_\alpha$ the number of solvent molecules bound to the ion, and $\tilde m_0$ the mass of the solvent molecule.
For the partial molar volume of the solvated ions, a similar relation holds, with $\tilde v_\alpha = \hat{v}_\alpha + \kappa_\alpha \tilde  v_0$ with the molar volume $ \hat{v}_\alpha$ of the central ion and $\tilde v_0$ of the solvent. 
We note that the approximation $\tilde m_\alpha / \tilde m_0 = \tilde v_\alpha / \tilde v_0 =: a_\alpha$ is also used.
Furthermore, we define the mass density $\tilde \rho(\tilde t,\tilde \vecx)$, the charge density $\tilde q(\tilde t,\tilde \vecx)$ and the barycentric velocity $\tilde \vecv(\tilde t, \tilde \vecx)$ by
$$\tilde \rho   =\smashoperator{ \sum_{\alpha = 0}^N} \tilde m_\alpha \tilde n_\alpha, \quad \tilde q = F \smashoperator{\sum_{\alpha = 0}^N} z_\alpha \tilde n_\alpha \quad \text{and} \quad \tilde \vecv = \frac{1}{\tilde \rho}\smashoperator{\sum_{\alpha = 0}^N} \tilde m_\alpha \tilde n_\alpha \tilde \vecv_\alpha ~,$$
where $\tilde \vecv_\alpha$ is the velocity field of constituent $\alpha$.
For the domain we consider an axially symmetric geometry $$\Omega = \left\{\tilde \vecx = (\tilde x, \tilde y, \tilde z)^T: 0 \leq \tilde z \leq \tilde L, \, 0 \leq \sqrt{\tilde x^2 + \tilde y^2} \leq \tilde R(\tilde z) \right\},$$ such as given in Figure \ref{fig:Domain}. The radius $\tilde R(\tilde z)$ is a function of $\tilde z$. At $\tilde z=0$ we prescribe an outlet boundary $S^\tout$ and at $\tilde z= \tilde L$ an inlet boundary $S^\tin$. We denote the surface of the ion channel by $S^\twall$ and the symmetry axis at $\tilde r = \sqrt{\tilde x^2 + \tilde y^2}=0$ by $S^0$.
The Nernst-Planck equations describe the evolution of the molar densities $\tilde n_\alpha(\tilde \vecx, \tilde t)$ for $\alpha = 1, \dots, N$ by
\begin{align}\label{NernstPlanck}
	\partial_{\tilde{t}} \tilde n_\alpha = -\tilde \nabla \cdot \tilde{\mathbf{J}}_\alpha = - \tilde \nabla \cdot \left[\tilde{\vecv} \tilde{n}_{\alpha} - \tilde{D}_{\alpha}  \tilde{n}_{\alpha} \left( \tilde \nabla \frac{\tilde \mu_\alpha}{R_G T} + z_\alpha \frac{F}{\kT} \tilde \nabla \tilde{\phi} \right) \right],
\end{align}
with diffusion coefficients $\tilde{D}_{\alpha}$, the effective chemical potential function $\tilde \mu_\alpha$, the molar gas constant $R_G$, the temperature $T$ and the Faraday constant $F$. Furthermore, the definition for the non-convective mass flux $\tilde \vecJ_\alpha = \tilde m_\alpha \tilde n_\alpha (\tilde \vecv_\alpha - \tilde \vecv)$ implies the identity
\begin{align}
    \sum_{\alpha = 0}^N \tilde \vecJ_\alpha = 0.
\end{align}
The chemical potential function $\tilde \mu_\alpha$ is derived from the free energy density, which consists of various contributions such as mixing entropy, polarization of the mixture, and/or mechanical contribution \cite{Borukhov1997,de2013non,Dreyer2013,Liu2013}.
Within this work we will consider effective chemical potential functions of the form
\begin{align}\label{eq:mugeneral}
    \tilde \mu_\alpha = \mu_\alpha^R + \ln{\tilde y_\alpha} - a_\alpha \ln{\tilde y_0},
\end{align} 
with a reference potential $\mu_\alpha^R$ and the mole fraction $\tilde y_\alpha = \tilde n_\alpha / \bar{n}$.  
The total concentration $\bar{n}$ of the mixture is given as $\bar n = \sum_{\alpha = 0}^N \tilde n_\alpha$ which yields $\tilde y_0 = 1 - \sum_{\alpha = 1}^N \tilde y_\alpha$ for the solvent. Since we consider an incompressible mixture, we need to take the incompressibility constraint $\sum_{\alpha=0}^N \tilde v_\alpha \tilde n_\alpha = 1$ into account. From that we can calculate the total concentration in terms of the mole fractions as follows $\bar n = 1/(\tilde v_0 + \sum_{\alpha = 1}^N (\tilde v_\alpha - \tilde v_0 ) \tilde y_\alpha)$.

By different choices of the volume fraction $a_\alpha$ we can recover different extension of the PNPS system from our approach. For example, choosing $a_\alpha = 0$ gives the classical PNP equations where ions are treated as point particles. Note that in this case one has to choose $\bar n = n^R = const.$, with e.g. $n^R = 1/v_0 = 55.4$\,\si{\mole \per \liter} the reference molar density of the pure solvent.
As an extension of the classical model, lattice-based models such as the Bikerman-Freise model \cite{Bikerman1942,Liu2013} are often used in literature, in which it is often assumed that all ions are the same size, such that $a_\alpha = 1$ since $v_\alpha = v_0$. In that case we also find $\bar n = 1/v_0 = n^R = const.$, which would correspond to the total number of lattice sites or the total volume.
Furthermore, considering different ion sizes and solvation effects one would get $a_\alpha > 0$ and $a_\alpha \neq 1$ with varying $\bar n$.
In summary we can recover three different PNP variations
\begin{align}
    \tilde \mu_\alpha = \begin{cases}
        \mu_\alpha^R + \ln{\frac{\tilde n_\alpha}{n^R}} \quad &\text{Nernst-Planck} \\
        \mu_\alpha^R + \ln{\frac{\tilde n_\alpha}{n^R}} - \ln{1 - \sum_{\alpha=0}^N \frac{\tilde n_\alpha}{n^R}} \quad &\text{Bikerman} \\
        \mu_\alpha^R + \ln{\frac{\tilde n_\alpha}{\bar n}} - a_\alpha \ln{1 - \sum_{\alpha=0}^N \frac{\tilde n_\alpha}{\bar n}} \quad &\text{idealized incompressible mixture}.
    \end{cases}
\end{align}
The electrostatic potential $\tilde{\phi}(\tilde \vecx,\tilde t)$ is given by the Poisson equation
\begin{align}\label{Poisson}
	-\tilde \nabla \cdot \left(\eps_0(1+\chi)\nabla \tilde{\phi} \right) = \tilde q(\tilde t, \tilde \vecx) = F \sum_{\alpha = 1}^N z_\alpha \tilde n_\alpha,
\end{align}
with vacuum permittivity $\eps_0$ and a constant dielectric susceptibility $\chi$. Note that $\eps_0(1+\chi) = \eps_r$ is also often referred to as the relative permittivity. The above system is further coupled with the Stokes equations to calculate the velocity $\tilde{\vecv}(\tilde \vecx,\tilde t)$ and the pressure $\tilde{p}(\vecx,t)$. Which are given by the continuity equation
\begin{align}\label{Contunuity}
	\tilde \nabla \cdot \tilde{\vecv} = 0
\end{align}
and the momentum balance
\begin{align} \label{MomentumBalance}
	-\nu \tilde \Delta \tilde{\vecv} + \tilde \nabla \tilde{p} = - \tilde q(\tilde t, \tilde \vecx)  \nabla \tilde{\phi},
\end{align}
with viscosity $\nu$.
We apply Dirichlet boundary conditions for concentrations and potential on the outlet and inlet 
\begin{subequations}
	\begin{align}
		\tilde n_\alpha = \tilde n_\alpha^{\tout}, \; \tilde{\phi} = \tilde \phi^\tout \; \text{on } S^{\tout} \; \text{and }
		\tilde n_\alpha = \tilde n_\alpha^{\tin}, \; \tilde{\phi} = \tilde \phi^\tin \; \text{on } S^{\tin},
	\end{align}
and normal stress boundary conditions for the velocity and the pressure
\begin{align}
	  \nu \tilde \nabla \tilde \vecv \cdot \tilde \vecn - \tilde p \tilde \vecn = \tilde p^\tout \tilde \vecn \; \text{on } S^{\tout}, \; \nu \tilde \nabla \tilde \vecv \cdot \tilde \vecn - \tilde p \tilde \vecn = \tilde p^\tin \tilde \vecn \; \text{on } S^{\tin}.
\end{align}
On the channel wall we apply a constant surface charge
\begin{align}
    -\eps_0(1+\chi) \tilde \nabla \tilde \phi \cdot \tilde \vecn = \tilde{\sigma}(\tilde z) \; \text{on } S^{\text{wall}},
\end{align}
a slip condition $\tilde \vecv_s$ for the velocity and no-flux boundary condition for the ion species 
\begin{align}
	\tilde b \nabla \tilde{\vecv} \cdot \tilde \vecn = \tilde \vecv_s , \;
	\tilde{\vecJ}_\alpha \cdot \tilde \vecn = 0 \; \text{on } S^{\text{wall}}.
\end{align}
Furthermore, we define the following symmetry conditions at $\tilde r=0$
\begin{align}
  	 \tilde \vecJ_\alpha \cdot \tilde \vecn = 0, \;
    \eps_0(1+\chi) \tilde \nabla \tilde \phi \cdot \tilde \vecn = 0, \;
    \nu \tilde \nabla \tilde \vecv \cdot \tilde \vecn - \tilde p \tilde \vecn = \mathbf{0}  \; \text{on } S^0.
\end{align}
\end{subequations}
A crucial property of this system is the ionic current that flows through the pore. The current passing through surface $S$ in $z$-direction can be calculated as follows
\begin{align}
   \tilde I_\alpha  |_{S} = z_\alpha F \int_{S} \tilde \vecJ_\alpha \cdot \tilde \vecn \ud S,
\end{align}
with the total current
\begin{align}
    \tilde I |_{S} = \sum_{\alpha = 1}^N \tilde I_\alpha |_{S}.
\end{align}
\paragraph{Scaled axially symmetric problem in the large-aspect-ratio limit.}
Since we consider a rotational symmetric problem we rewrite the system in cylindrical coordinates. The velocity field is given by $\tilde \vecv = (\tilde w,\tilde u)$ where $\tilde w$ is the radial component and $\tilde u$ the axial component. Similar for the slip condition $\tilde \vecv_s = (\tilde w_s,\tilde u_s)$.
For the large aspect ratio scaling we introduce a typical pore length $L_0 = 10^{-8}-10^{-6}\,\si{\meter}$ and a typical pore radius $R_0=10^{-9}\,\si{\meter}$. In addition, we define a typical time scale $\tau=10^{-9}\,\si{\second}$, a reference concentration $c^R=1\,\si{\mole \per \liter}$ and a reference charge $\sigma^R = 0.16 \,\si{\coulomb} $. We introduce the following dimensionless variables $\tilde{z} = L_0 z, \; \tilde{r} = R_0 r, \; \tilde{t} = \tau t, \; \tilde{D}_\pm = D^R k_\pm, \; \tilde{n}_\alpha = n^R n_\alpha, \; \tilde \mu_\alpha = \mu^R \mu_\alpha, \; \tilde{\phi} = \phi^R \phi, \;\tilde{u} = u^R u, \; \tilde{w} = \delta u^R w, \; \tilde{p} = p^R p \; \tilde{\sigma} = \sigma^R \sigma \; \tilde{b} = b^R b$. 
This gives us the ability to introduce the following dimensionless parameters 
\begin{subequations}
\begin{align}
	D^R = \frac{L_0^2}{\tau}, \;
    \mu^R = \kT, \;
	\phi^R = \frac{\kT}{F},\; 
	p^R = n^R \kT, \\ \nonumber
	u^R = \frac{n^R \kT R_0^2}{\nu L_0 }, \; I^R = \frac{F R_0^2 D^R n^R}{L_0}, \; b^R = L_0
\end{align}
and
\begin{align}
    \Lambda^2 = \frac{\eps_0(1+\chi) \phi^R}{F c^R R_0^2} = \frac{\lambda^2}{R_0^2}, \;
    \gamma = \frac{R_0 \sigma^R}{\eps_0(1+\chi)\phi^R},\;
    P_e = \frac{u^R L_0 }{D^R}, \;
    \delta = \frac{R_0}{L_0 },
	\end{align}
\end{subequations}
where $\mathcal{O}(1/\delta)\gg1$.
The scaled Nernst-Planck equations in cylindrical coordinates are given by
\begin{subequations}\label{eqn:scaled}
	\begin{align}\label{eqn:scaled:a}
		\delta^2 \partial_t n_\alpha = -\delta^2 P_e \left( w \partial_r n_\alpha + u \partial_z n_\alpha \right) + \delta^2 \partial_z f_\pm + \frac{1}{r} \partial_r \left(r g_\pm\right),
	\end{align}
    with the functions
	\begin{align}\label{eqn:scaled:b}
		f_\alpha = k_\alpha n_\alpha \left(\partial_z \mu_\alpha  + z_\alpha \partial_z \phi \right),
	\end{align}
    and
	\begin{align}\label{eqn:scaled:c}
		g_\alpha = k_\alpha n_\alpha \left(\partial_r \mu_\alpha +  z_\alpha \partial_r \phi \right)
	\end{align}
    for ease of notation.
    For the Poisson equation we find
	\begin{align}\label{eqn:scaled:d}
		-\delta^2 \partial_{zz} \phi - \frac{1}{r} \partial_r \left(r \partial_r \phi \right) = \frac{q}{\Lambda^2}.
	\end{align}
    The dimensionless continuity equation in cylindrical coordinates yields
	\begin{align}\label{eqn:scaled:e}
		\frac{1}{r}\partial_r \left(r w\right) + \partial_z u  = 0,
	\end{align}
    and for the momentum balance we have
	\begin{align}\label{eqn:scaled:f}
		\delta^4 \partial_{zz} w + \delta^2 \left[\frac{1}{r} \partial_r \left(r \partial_r w\right) - \frac{1}{r^2} w\right] = \partial_r p + q\partial_r \phi,
	\end{align}
    and
	\begin{align}\label{eqn:scaled:g}
		\delta^2 \partial_{zz} u + \frac{1}{r} \partial_r \left(r \partial_r u \right) =  \partial_z p + q \partial_z \phi.
	\end{align}
\end{subequations}
The scaled boundary conditions on the in- and outlet are
\begin{subequations}
	\begin{align}\label{eqn:scaledbc:a}
		n_\alpha = n_\alpha^{\tout}, \; \phi = \phi^\tout \; \text{on } S^{\tout}, \text{ and }
		n_\alpha = n_\alpha^{\tin}, \; \phi = \phi^\tin  \; \text{on } S^{\tin},
	\end{align}
	\begin{align}\label{eqn:scaledbc:b}
		\delta^2 \partial_z u - p = p^\tout\; \text{on } S^{\tout}, \quad \delta^2 \partial_z u - p = p^\tin \; \text{on } S^{\tin}. 
	\end{align}
    On the channel wall we have
	\begin{align}\label{eqn:scaledbc:c}
		\partial_{r} \phi - \delta^2 R'(z) \partial_{z} \phi = -\sigma(z) \gamma \left(1+ \delta^2 ( R'(z) )^2 \right)^{1/2} \; \text{on } S^{\text{wall}},
	\end{align}
	\begin{align}\label{eqn:scaledbc:d}
		b \partial_r w = \delta w_s, \; b \partial_z u = u_s, \; 
		 \delta^2 R'(z)  f_\alpha - g_\alpha = 0 \; \text{on } S^{\text{wall}},
	\end{align}
        with the normal vector
    \begin{align*}
        \vecn = \left(\mathbf{e}_r - \delta R'(z) \mathbf{e}_z \right)\left(1+ \delta^2 ( R'(z) )^2 \right)^{-1/2},
    \end{align*}
    where $\mathbf{e}_r$ is the normal vector in $r-$direction and $\mathbf{e}_z$ in $z-$direction, 
	\begin{align}\label{eqn:scaledbc:e}
		\partial_r \phi = 0,\;
		\delta^2 \partial_r w - p = 0, \; \partial_{r} u = 0, \text{ on } S^0
	\end{align}
    and 
    \begin{align}\label{eqn:scaledbc:f}
        -P_e \left(\delta^2 w \, n_\alpha + \delta u \, n_\alpha \right) + \delta f_\alpha + g_\alpha = 0 \text{ on } S^0.
    \end{align}
\end{subequations}
The current flowing through $S$ at $z$-direction is given by the integral over the disk
\begin{align}
    I_\alpha = z_\alpha 2\pi \int_0^{R(z=z_S)} (P_e n_\alpha u - f_\alpha) |_{z=z_S} r \ud r.
\end{align}
\section{Large aspect ratio asymptotics of the system}
We use the following ansatz in order to find an asymptotic solution for the system in the limit of 
$\delta \rightarrow 0$:
\begin{align}
    (n_\alpha, \mu_\alpha, \phi, u, w, p ) =& (n^0_\alpha + \delta^2 n_\alpha^1, 
    \mu_\alpha^0 + \delta^2 \mu_\alpha^1,
    \phi^0 + \delta^2 \phi^1, 
    u^0 + \delta^2 u^1, 
    w^0 + \delta^2 w^1, 
    \phi^0 + \delta^2 p^1) \\ \nonumber
    &+ O(\delta^4).
\end{align} 
In this study, the focus is on the case of the Debye length being equivalent or smaller to the radius of the pore, i.e., $\Lambda  \leq \mathcal{O}(1)$ and $\lambda \leq \mathcal{O}(R_0)$. 
Additionally, we assume a substantial surface charge, characterized by the condition $\gamma = \mathcal{O}(1)$ and a large aspect ratio, denoted by $\mathcal{O}(1/\delta)\gg 1$. 

Note that we calculate the asymptotic model for a system with no-slip boundary condition on $S^\twall$ for the velocity, such that $u|_{r=R(z)} = w|_{r=R(z)} = 0$.
\paragraph{Poisson equation.} To leading order $\mathcal{O}(1)$, \eqref{eqn:scaled:d} becomes
\begin{subequations}\label{phizerobvp}
\begin{align} \label{SolutionPhiR}
	-\frac{1}{r} \partial_r \left(r \partial_r \phi^0 \right) = \frac{q^0}{\Lambda^2}.
\end{align}
Since $q^0$ is a function of $\phi^0$ this equation has to be solved later together with the leading order boundary conditions
(from \eqref{eqn:scaledbc:c} and \eqref{eqn:scaledbc:e}, respectively)
\begin{align}\label{phizerobc1}
	\partial_r \phi^0 |_{r=R(z)} = -\gamma \sigma(z),  \quad \text{and} \quad \partial_r\phi^0|_{r=0} = 0.
\end{align}
\end{subequations}
The solution for $\phi^0$ of this boundary value problem is not unique, since the homogeneous
problem for $\phi^0$ is satisfied by an arbitrary $\phi_z^0(z)$. Therefore, 
if there is a particular solution of the inhomogeneous problem,  say $\phi_r^0(r,z)$, then any other solution can be achieved
via 
\begin{align}\label{AnsatzPhi}
    \phi^0(r,z) = \phi_r^0(r,z) + \phi_z^0(z),
\end{align}
with an arbitrary $\phi_z^0(z)$. Moreover, existence requires a solvability condition to
be satisfied, which is obtained by simply integrating \eqref{SolutionPhiR} and
using \eqref{phizerobc1}. The result is the electroneutrality condition
\begin{equation}
\int_0^{R(z)}  q rdr  - \Lambda^2\gamma\sigma R(z) \overset{!}{=}0.
\end{equation}
Notice that this condition includes the surface charges, so the net charge {\em in
the fluid itself} must be negative inside the pipe. 
For the second order $\mathcal{O}(\delta^2)$ we find from equation \eqref{eqn:scaled:d}
\begin{align}
	\frac{1}{r} \partial_r (r \partial_r \phi^1) + 	\partial_z (\partial_z \phi^0 )  = \frac{q^1}{\Lambda^2}.
\end{align}
Integrating the equation with respect to $r$ and using the boundary condition
(from \eqref{eqn:scaledbc:c} and \eqref{eqn:scaledbc:e}, respectively)
\begin{align}\label{eq:CalculationPoisson}
	\partial_r \phi^1 |_{r=R(z)} - R'(z) \partial_z \phi^0 |_{r=R(z)} = 0 \quad \text{and} \quad
	\partial_r\phi^1|_{r=0} = 0.
\end{align}
gives
\begin{align}
    0 &=\int_{0}^{R(z)} \partial_r (r \partial_r \phi^1) \ud r + 	\int_{0}^{R(z)} \partial_z (\partial_z \phi^0 ) r \ud r = \partial_z \int_{0}^{R(z)} \partial_z \phi^0_z r \ud r.
\end{align}   
For a detailed computation see equation \eqref{eq:AppendixPoisson} in the appendix. This can now be used to fix expression for $\phi_z^0$ in \eqref{AnsatzPhi}. Plugging that decomposition
into the preceding identity yields
\begin{align}
	\phi_z^0 = C_1 \int_{0}^{z} \frac{1}{R^2(z)} \ud z + C_2.
\end{align}
Applying the boundary conditions \eqref{eqn:scaled:a}
\begin{align}
	\phi_z^0|_{z=L} = \phi^\tin, \; \phi_r^0|_{z=L} = 0 \; \text{ and } \; \phi_z^0|_{z=0} = \phi^\tout, \; \phi_r^0|_{z=0} = 0,
\end{align}
allows us to determine the constants $C_1$ and $C_2$. We obtain 
\begin{subequations}
    \begin{align}  \label{SolutionPhiZ}
    	\phi_z^0 = (\phi^\tin - \phi^\tout) \frac{I_{\phi_z}(z) - I_{\phi_z}(0)}{I_{\phi_z}(L) - I_{\phi_z}(0)} + \phi^\tout,
    \end{align}
    with 
    \begin{align}
    	I_{\phi_z}(z) =  \int_{0}^{z} \frac{1}{R^2(z)} \ud z.
    \end{align}
\end{subequations}
\paragraph{Nernst-Planck equations} In the leading order $\mathcal{O}(1)$ for the concentrations we have from \eqref{eqn:scaled:a} and \eqref{eqn:scaled:c}
\begin{align}
	\frac{1}{r} \partial_r \left(r g^{0}_\alpha \right) &= 0.
\end{align} 
We integrate with respect to $r$ which gives together with the leading order of
the third boundary condition in \eqref{eqn:scaledbc:d}
\begin{align}
     \mu_\alpha^0 = \ln{y_\alpha^0} - a_\alpha \ln{y_0^0} = - z_\alpha \phi^0.
\end{align}
Using \eqref{AnsatzPhi} for $\phi^0$ we obtain 
\begin{align}\label{SolutionCR}
	y_\alpha^0 = Q_\alpha(z,t) \exp{- z_\alpha \phi_r^0} (y_0^0)^{a_\alpha}
\end{align}
with arbitrary $Q_\alpha(z,t)$ and $y_0^0 = y_0^0(y_\alpha^0)$.
Notice the leading order of \eqref{eqn:scaledbc:d} imposes the condition $g_\alpha^0=0$, but this is automatically the case since $\phi^0$ and hence $\phi_r^0$ satisfies \eqref{phizerobc1}.
For the second order $\mathcal{O}(\delta^2)$ we obtain from equation \eqref{eqn:scaled:a}
\begin{align}
    \partial_t n_\alpha = -P_e \left(w^0 \partial_r n^0_\alpha + u^0 \partial_z n^0_\alpha \right) +  \partial_z f^0_\alpha + \frac{1}{r} \partial_r \left(r g^1_\alpha \right).
\end{align}
In order to solve this equation we will introduce the stream function 
\begin{align}
    s = \int_{0}^{r} u^0 r \ud r,
\end{align}
such that $u^0 = r^{-1} \partial_r s$ and $w = - r^{-1} \partial_z s$.
Integrating above equation gives together with the boundary condition
$R'(z)  f^0_\alpha |_{r=R(z)} - g^1_\alpha|_{r=R(z)} = 0$
the following result
\begin{align}\label{eq:CalculationConcentration}
    \partial_t \int_0^{R(z)} n^0_\alpha r \ud r &= -P_e \int_{0}^{R(z)} \left(-\partial_z s \partial_r n^0_\alpha + \partial_r s \partial_z n^0_\alpha \right) \ud r + \int_{0}^{R(z)} \partial_z f^0_\alpha r \ud r + \int_{0}^{R(z)} \partial_r \left(r g^1_\alpha\right) \ud r \\ \nonumber
    &= P_e \partial_z \int_{0}^{R(z)} n^0_\alpha u^0 r \ud r + k_\alpha \partial_z \int_{0}^{R(z)} n_\alpha^0 \left( \partial_z \mu^0_\alpha + z_\alpha \partial_z \phi^0 \right) r \ud r.
\end{align}
For a detailed calculation see Appendix \ref{sec:AppendixCalculationConcentration}.
Rewriting the molar density in terms of the molar fraction, i.e.$n_\alpha = \bar n y_\alpha$, using solution \eqref{SolutionCR} and the ansatz \eqref{AnsatzPhi} further gives
\begin{align}\label{deqpm}
    \partial_t \left( Q_\alpha \int_0^{R(z)} \bar n^0 (y_0^0)^{a_\alpha} \exp{- z_\alpha \phi^0_r} r \ud r \right) &= \partial_z \left[P_e  Q_\alpha \int_{0}^{R(z)} \bar n^0 (y_0^0)^{a_\alpha} \exp{- z_\alpha \phi_r^0} u^0 r \ud r \right. \\ \nonumber 
    &\left. \qquad \qquad  + k_\alpha z _\alpha Q_\alpha (\partial_z \phi^0_z) \int_{0}^{R(z)} \bar n^0 (y_0^0)^{a_\alpha} \exp{- z_\alpha \phi_r^0} r \ud r \right. \\ \nonumber
    & \left. \qquad \qquad + k_\alpha \partial_z Q_\alpha \int_{0}^{R(z)} \bar n^0 (y_0^0)^{a_\alpha} \exp{- z_\alpha \phi_r^0} r \ud r  \right].
\end{align}
This equation needs to be solved for $Q_\alpha$ together with the boundary conditions \eqref{eqn:scaledbc:a}
\begin{align}\label{Qbc}
    Q_\alpha |_{z=0} = \frac{n^\tout}{\bar n^\tout} \left( \frac{1}{v_0 \bar n^\tout} - \sum_{\beta=1}^N a_\beta \frac{n_\beta^\tout}{\bar n^\tout} \right)^{-a_\alpha} \; \text{ and } \; Q_\alpha |_{z=L} =\frac{n^\tin}{\bar n^\tin} \left( \frac{1}{v_0 \bar n^\tin} - \sum_{\beta=1}^N a_\beta \frac{n_\beta^\tin}{\bar n^\tin} \right)^{-a_\alpha}.
\end{align} 
\paragraph{Momentum balance for no-slip.} We can obtain a solution for the pressure from the leading order $\mathcal{O}(1)$ from equation \eqref{eqn:scaled:f}
\begin{align}
	\partial_r p^0 = -q \partial_r \phi^0.
\end{align}
With a similar argumentation as for the potential we use the following ansatz for the pressure
\begin{align}
    p^0(r,z) = p_z^0(z) + p^0_r(r,z).
\end{align}
Using solution \eqref{SolutionCR}, integrating with respect to $r$ yields 
\begin{align} 
	p_r^0 &= C_1 - \sum_{\alpha=1}^N z_\alpha Q_\alpha \int_0^r \bar n^0 (y_0^0)^{a_\alpha} \exp{-z_\alpha \phi_r^0} \ud r \\ \nonumber
    & = C_1 + \sum_{\alpha=1}^N Q_\alpha \left[\bar n^0 (y_0^0)^{a_\alpha} \exp{-z_\alpha \phi_r^0} - \int_0^r   \exp{-z_\alpha \phi_r^0} \partial_r (\bar n^0 (y_0^0)^{a_\alpha} ) \ud r \right] .
\end{align}
Together with the symmetry condition \eqref{eqn:scaledbc:e} $p_r^0 |_{r=0} = 0$ we get
\begin{subequations}\label{SolutioPR}
    \begin{align}
        p_r^0 = \check Q(r,z) - \check Q(r=0,z),
    \end{align}
    with
    \begin{align}
        \check Q(r,z) = \sum_{\alpha=1}^N Q_\alpha \left[\bar n^0 (y_0^0)^{a_\alpha} \exp{-z_\alpha \phi_r^0} - \int_0^r   \exp{-z_\alpha \phi_r^0} \partial_r (\bar n^0 (y_0^0)^{a_\alpha} ) \ud r \right].
    \end{align}
\end{subequations}
From equation \eqref{eqn:scaled:g} we can derive a solution for the axial component of the velocity $u^0$
\begin{align}
    \frac{1}{r} \partial_r (r \partial_r u^0) &= \partial_z p^0 + q \partial_z \phi^0.
\end{align}
We rewrite the derivatives of the pressure and the potential as follows
\begin{subequations}
    \begin{align}
        \partial_z p^0 = \partial_z p_z^0 + \partial_z p_r^0 = \Pi_z(z) + \Pi_r(r,z) 
    \end{align}
    and
    \begin{align}
        \partial_z \phi^0(r,z) &= \partial_z \phi_z^0 + \partial_z \phi_r^0  =\text{E}_z(z) + \text{E}_r(r,z).
    \end{align}
\end{subequations}
Using equation \eqref{SolutionPhiR} we find
\begin{align}
    \frac{1}{r} \partial_r (r \partial_r u^0) = \Pi_z + \Pi_r - \Lambda^2 (\tE_z + \tE_r)  \frac{1}{r} \partial_r (r \partial_r \phi_r^0) . 
\end{align}
Integrating with respect to $r$ gives
\begin{subequations}
    \begin{align}
        r \partial_r u^0 = \Pi_z \frac{1}{2} r^2 + K_1(r,z) - \text{E}_z\Lambda^2 r \partial_r \phi^0 + h_1,
    \end{align}
    with
    \begin{align}
            K_1(r,z) = \int_0^r (\Pi_r r - \Lambda^2 \text{E}_r \partial_r (r \partial_r \phi_r^0)) \ud r.
    \end{align}
\end{subequations}
We can further rewrite $K_1$ by using equation \eqref{SolutionPhiR} to get
 \begin{align}
        K_1(r,z) &= \int_0^r (\Pi_r r - \Lambda^2 \text{E}_r \partial_r (r \partial_r \phi_r^0)) \ud r \\ \nonumber
        &= \int_0^r \left(\partial_z p_r^0  +  q \partial_z \phi^0_r \right) r \ud r \\ \nonumber
        &= \int_0^r \left(\partial_z \check Q - \partial_z \check Q |_{r=0}  +  \sum_{\alpha=1}^N \bar n^0 z_\alpha Q_\alpha (y_0^0)^{a_\alpha} \exp{-z_\alpha \phi_r^0} \partial_z \phi^0_r \right) r \ud r \\ \nonumber 
        &= \int_0^r \left( \sum_{\alpha=1}^N \partial_z [Q_\alpha \bar n^0 (y_0^0)^{a_\alpha} ] \exp{- z_\alpha \phi_r^0} - \partial_z \check Q |_{r=0}  - \partial_z \sum_{\alpha=1}^N Q_\alpha \int_0^r   \exp{-z_\alpha \phi_r^0} \partial_r [\bar n^0 (y_0^0)^{a_\alpha} ] \ud r \right) r \ud r.
\end{align}
With the symmetry condition \eqref{eqn:scaledbc:e} at $r=0$ we find
$\partial_r u^0  |_{r=0}  = r^{-1} (h_1 +  K_1(r=0,z)) = 0 $ and hence $h_1 \overset{!}{=} - K_1(r=0,z)$. Note that the integrand of $K_1$ is the equilibrium between the axial gradient of the electric double layer (EDL) pressure and the axial electrostatic force. density.
Integrating again with respect to $r$ gives
\begin{subequations}
    \begin{align}
        u^0 = \Pi_z \frac{1}{4} r^2 + K_2(r,z) - \text{E}_z\Lambda^2 \phi^0  + h_2,
    \end{align}
    with
    \begin{align}
            K_2(r,z) = \int_0^r \frac{1}{r} (K_1(r,z) - K_1(r=0,z)) \ud r.
    \end{align}
\end{subequations}
Note that $K_2$ is an EDL-driven axial flow contribution.
Together with the no-slip condition \eqref{eqn:scaled:d} $u^0 |_{r=R(z)} = 0$, we get
\begin{subequations}
    \begin{align}\label{SolutionU}
    	u^0 &=  \Pi_z \frac{1}{4} \left(r^2 - R^2(z) \right) + K_2(r,z) - K_2(r=R(z),z)  - \text{E}_z \Lambda^2 \zeta(r,z),
    \end{align}
    with the $\zeta-$potential
    \begin{align}
        \zeta = \phi_r^0 - \phi_r^0 |_{r=R(z)} ,
    \end{align}
\end{subequations}
\paragraph{Continuity equation for no-slip.}
We need to ensure that the continuity equation \eqref{eqn:scaled:e} is fulfilled. Therefore, we integrate the leading order $\mathcal{O}(1)$ with respect to $r$ from $0$ to $R(z)$ which gives
\begin{align}
	0 &= \int_{0}^{R(z)} \left( \partial_r (r w^0) + (\partial_z u^0) r \right) \ud r \\ \nonumber 
    &= R(z) \underbrace{\left(w^0 |_{r=R(z)} - R'(z) u^0 |_{r=R(z)}\right)}_{=0}  +  \partial_z \int_{0}^{R(z)}  u^0 r \ud r.
\end{align}
Together with no-slip $w^0 |_{r=R(z)} = u^0 |_{r=R(z)} = 0$ and the solution \eqref{SolutionU} for $u^0$ we have
\begin{subequations}
    \begin{align}
    	\partial_z  &\left[\Pi_z \frac{1}{4} \int_{0}^{R(z)} \left(r^2 - R^2(z) \right) r \ud r + K_3(z) - \text{E}_ z \Lambda^2 \int_{0}^{R(z)} \zeta \, r \ud r \right]  = 0,
    \end{align}
    with 
    \begin{align}
        K_3(z) = \int_0^{R(z)} (K_2(r,z) - K_2(r=0,z) ) r \ud r.
    \end{align}
\end{subequations}
Solving this equation gives
\begin{subequations}
    \begin{align}
     	p_z^0 =& C_2 + C_1 I_p(z)  -  \int_{0}^z \left( \frac{16 \Lambda^2\text{E}_z}{R^4(z)}  \int_{0}^{R(z)} \zeta \, r \ud r \right) \ud z + \int_0^z \frac{16 K_3(z)}{R^4(z)}  \ud z ,
    \end{align}
    with 
    \begin{align}
        I_p(z) = \int_{0}^{z} \frac{1}{R^4(z)} \ud z
    \end{align}
\end{subequations}
with the boundary conditions \eqref{eqn:scaledbc:b}
\begin{align}
	p^0_z |_{z=0} = -p^{\tout}, \; p^0_r |_{z=0} = 0 \; \text{ and } \; p^0_z |_{z=L} = -p^\tin, \; p^0_r |_{z=L} =0,
\end{align}
the following solution
\begin{subequations}
    \begin{align}
    	p_z^0 &=  -(p^\tin - p^\tout)  \frac{I_p(z)-I_p(z=0)}{I_p(z=L) - I_p(z=0)} - p^\tout - \Psi(z)  + K_4(z),
    \end{align}
    with 
    \begin{align}
        \Psi(z) = I_{\text{E}_z}(z) - I_{\text{E}_z}(z=0) - \frac{I_{\text{E}_z}(z=L) - I_{\text{E}_z}(z=0)}{I_p(z=L)-I_p(z=0)} \left(I_p(z) - I_p(z=0)\right)
    \end{align}
    \begin{align}
    	I_{\text{E}_z}(z) = \int_{0}^z \left(\frac{16 \Lambda^2\text{E}_z}{R^4(z)}  \int_{0}^{R(z)} \left[\phi_r^0 - \phi_r^0 |_{r=R(z)} \right] r \ud r \right) \ud z ,
    \end{align}
    \begin{align}
        K_4(z) = I_{K_3}(z) - I_{K_3}(z=0) - \frac{I_{K_3}(z=L)- I_{K_3}(z=0)}{I_p(z=L)-I_p(z=0)} \left(I_p(z) - I_p(z=0)\right)
    \end{align}
    and
    \begin{align}
        I_{K_3}(z) = \int_0^z \frac{16}{R^4(z)} K_3(z) \ud z.
    \end{align}
\end{subequations}
\paragraph{Current.} The leading order $\mathcal{O}(1)$ currents are calculated as follows
\begin{align}
    I_\alpha  |_{z=z_S} =& z_\alpha 2 \pi  k_\alpha \left( \partial_z Q_\alpha + z_\alpha Q_\alpha \tE_z \right)  |_{z=z_S} \int_0^{R(z=z_S)} (\bar n^0 (y_0^0)^{a_\alpha}  \exp{- z_\alpha \phi^0_r} ) |_{z=z_S} r \ud r \\ \nonumber
    &- z_\alpha 2 \pi  P_e Q_\alpha  |_{z=z_S} \int_0^{R(z=z_S)} (u^0 \bar n^0 (y_0^0)^{a_\alpha} \exp{-z_\alpha \phi^0_r} ) |_{z=z_S} r \ud r.
\end{align}
\section{The asymptotic model for the generalized PNPS system}
In order to get the solutions for $\phi_r^0(r,z)$ and $Q_\alpha(z)$ we have to solve the following three coupled 1D ordinary differential equations (ODEs) numerically. In the $r$-direction we have to solve
\begin{subequations}
    \begin{align}\label{eq:FinalPhirGeneral}
        - \frac{1}{r} \partial_r ( r \partial_r \phi^0_r) = \frac{\bar n^0}{\Lambda^2} \sum_{\alpha=1}^N z_\alpha Q_\alpha (y_0^0)^{a_\alpha} \exp{-z_\alpha \phi_r^0},
    \end{align}
    with the boundary conditions
    \begin{align}\label{eq:FinalPhirBCGeneral}
        \partial_r \phi^0_r |_{r=R(z)} = -\gamma \sigma(z), \quad \text{and} \quad \partial_r \phi^0_r |_{r=0} = 0
    \end{align}
\end{subequations}
for every point in $z$. In $z-$direction we have to calculate
\begin{subequations}
    \begin{align}\label{eq:FinalQGeneral}
        \partial_t \left( Q_\alpha H_1(y_\alpha^0,\phi^0_r) \right) &= \partial_z \left[ k_\alpha \left( \partial_z Q_\alpha + z_\alpha Q_\alpha \partial_z \phi_z^0 \right) H_1(y_\alpha^0,\phi^0_r)  - P_e Q_\alpha H_2(y_\alpha^0,\phi^0_r,u^0) \right]
    \end{align}
    for every timestep in $t$,
    with
    \begin{align}\label{eq:FinalQIntegralsGeneral}
        H_1 &= \int_0^{R(z)} \bar n^0 (y_0^0(y_\alpha^0) )^{a_\alpha}\exp{-z_\alpha \phi^0_r} r\ud r, \quad
        H_2 = \int_0^{R(z)} u^0 \bar n^0 (y_0^0(y_\alpha^0) )^{a_\alpha}\exp{-z_\alpha \phi^0_r} r\ud r
    \end{align}
    and the boundary conditions
    \begin{align}\label{eq:FinalQBCGeneral}
        Q_\alpha |_{z=0} = \frac{n^\tout}{\bar n^\tout} \left( \frac{1}{v_0 \bar n^\tout} - \sum_{\beta=1}^N a_\beta \frac{n_\beta^\tout}{\bar n^\tout} \right)^{-a_\alpha}, \; Q_\alpha |_{z=L} =\frac{n^\tin}{\bar n^\tin} \left( \frac{1}{v_0 \bar n^\tin} - \sum_{\beta=1}^N a_\beta \frac{n_\beta^\tin}{\bar n^\tin} \right)^{-a_\alpha}.
\end{align} 
\end{subequations}
From that we can calculate the concentrations with
\begin{subequations}
    \begin{align}\label{eq:FinalyalphaGeneral}
        n_\alpha^0 = \bar n^0 \, y_\alpha^0 = Q_\alpha \bar n^0 (y_0^0)^{a_\alpha} \exp{-z_\alpha \phi_r^0} ,
    \end{align}
    \begin{align}\label{eq:Finaly0General}
        y_0^0 = 1 - \sum_{\alpha=1}^N y_\alpha \quad \curvearrowright \quad 1 - y_0^0 - \sum_{\alpha=1}^N Q_\alpha (y_0^0)^{a_\alpha} \exp{-z_\alpha \phi_r^0} \overset{!}{=} 0,
    \end{align}
    \begin{align}\label{eq:FinalnbarGeneral}
        \bar n^0 = \frac{1}{v_0 + \sum_{\alpha=1}^N (v_\alpha - v_0)y_\alpha^0},
    \end{align}
\end{subequations}
and the total electrostatic potential
\begin{subequations}
    \begin{align}\label{eq:FinalPhiGeneral}
        \phi^0(r,z) = \phi^0_r(r,z) + \phi_z^0(z),
    \end{align}
    with
    \begin{align}\label{eq:FinalPhizGeneral}
        \phi^0_z(z) = (\phi^\tin - \phi^\tout) \Delta I_{\varphi_z} + \phi^\tout,
    \end{align}
    and
    \begin{align}
         \Delta I_{\phi_z} = \frac{I_{\phi_z}(z) - I_{\phi_z}(0)}{I_{\phi_z}(z=L) - I_{\phi_z}(z=0)} \quad \text{and} \quad I_{\phi_z} = \int_0^z \frac{1}{R^2(z)} \ud z.
    \end{align}
\end{subequations}
Note that for $R(z) = R=$\,const. we get $\Delta I_{\phi_z} = z/L$. The total pressure is given by
\begin{subequations}
    \begin{align}\label{eq:FinalPGeneral}
        p^0(r,z) = p_r^0(r,z) + p_z^0(z),
    \end{align}
    with
    \begin{align}\label{eq:FinalPRGeneral}
        p_r^0 = \check Q(r,z) - \check Q(r=0,z),
    \end{align}
    \begin{align}
        \check Q(r,z) = \sum_{\alpha=1}^N Q_\alpha \left[\bar n (y_0^0)^{a_\alpha} \exp{-z_\alpha \phi_r^0} - \int_0^r   \exp{-z_\alpha \phi_r^0} \partial_r (\bar n (y_0^0)^{a_\alpha} ) \ud r \right]
    \end{align}
    and
    \begin{align}\label{eq:FinalPzGeneral}
    	p_z^0(z) &=  -(p^\tin - p^\tout) \Delta I_p(z) - p^\tout  - \Psi(z) + K_4(z),
    \end{align}
    with 
    \begin{align}\label{eq:FinalPzIntegralsGeneral}
    	 \Delta I_p = \frac{I_p(z) - I_p(z=0)}{I_p(z=L) - I_p(z=0)} \quad \text{and} \quad I_p(z) = \int_{0}^{z} \frac{1}{R^4(z)} \ud z,
    \end{align}
    such that for $R(z) = R=$\,const. we get $\Delta I_{p} = z/L$,
    \begin{align}\label{eq:FinalPzIntegrals1General}
        \Psi(z) = I_{\text{E}_z}(z) - I_{\text{E}_z}(z=0) - (I_{\text{E}_z}(z=L) - I_{\text{E}_z}(z=0)) \Delta I_p(z),
    \end{align}
    \begin{align}\label{eq:FinalPzIntegrals2General}
        K_4(z) = I_{K_3}(z) - I_{K_3}(z=0) - (I_{K_3}(z=L)- I_{K_3}(z=0)) \Delta I_p(z)
    \end{align}
    and
    \begin{align}\label{eq:FinalPzIntegrals3General}
        I_{\text{E}_z}(z) = \int_{0}^z \left(\frac{16 \Lambda^2 (\partial_z \phi^0_z)}{R^4(z)}  \int_{0}^{R(z)} \zeta \, r \ud r \right) \ud z \quad \text{and} \quad I_{K_3}(z) = \int_0^z \frac{16}{R^4(z)} K_3(z) \ud z,
    \end{align}
    where
    \begin{align}\label{eq:FinalZetaGeneral}
        \zeta = \phi_r^0 - \phi_r^0 |_{r=R(z)} ,
    \end{align}
    \begin{align}\label{eq:FinalK3General}
        K_3(z) = \int_0^{R(z)} (K_2(r,z) - K_2(r=0,z) ) r \ud r,
    \end{align}
    \begin{align}\label{eq:FinalK2General}
        K_2(r,z) = \int_0^r \frac{1}{r} (K_1(r,z) - K_1(r=0,z)) \ud r
    \end{align}
    and
    \begin{align}\label{eq:FinalK1General}
        K_1(r,z) &= \int_0^r \left( \sum_{\alpha=1}^N \partial_z [Q_\alpha \bar n (y_0^0)^{a_\alpha} ] \exp{- z_\alpha \phi_r^0} - \partial_z \check Q (r=0,z) \right. \\ \nonumber
        &\left. \qquad \qquad - \partial_z \sum_{\alpha=1}^N Q_\alpha \int_0^r   \exp{-z_\alpha \phi_r^0} \partial_r [\bar n (y_0^0)^{a_\alpha} ] \ud r \right) r \ud r.
\end{align}
\end{subequations}
The velocity is then calculated as
\begin{subequations}
    \begin{align}\label{eq:FinalVelocityGeneral}
    	u^0(r,z) =  (\partial_z p^0_z) \frac{1}{4} \left(r^2 - R^2(z) \right) - (\partial_z \phi_z^0 )\Lambda^2 \zeta(r,z) + K_2(r,z) - K_2(r=R(z),z) ,
    \end{align}
    with
    \begin{align}\label{eq:FinalPressureGradientGeneral}
        \partial_z p_z^0 = &- (p^\tin - p^\tout) \partial_z (\Delta I_p(z) )- \partial_z \Psi(z) + \partial_z K_4(z)
    \end{align}
    and
    \begin{align}\label{eq:FinalPotentialGradientGeneral}
        \partial_z \phi^0_z = (\phi^\tin - \phi^\tout) \partial_z (\Delta I_{\phi_z}(z)).
    \end{align}
\end{subequations}
For the currents we have
\begin{subequations}
   \begin{align}\label{eq:FinalCurrentGeneral}
        I_\alpha  |_{z=z_S} = 2 \pi  z_\alpha k_\alpha \left( \partial_z Q_\alpha + z_\alpha Q_\alpha \partial_z \phi_z^0 \right)  |_{z=z_S} H_1  |_{z=z_S} - 2 \pi z_\alpha  P_e Q_\alpha  |_{z=z_S} H_2  |_{z=z_S}.
    \end{align}
    with the total current given by
    \begin{align}
        I  |_{z=z_S} = \sum_{\alpha=1}^N I_\alpha  |_{z=z_S}.
    \end{align}
\end{subequations}
\section{Predictions of the asymptotic model}
In this section we illustrate the wide spectrum of applications of the classical PNPS model for typical scenarios in nanotechnological and biological applications. The pseudocode of the numerical sheme used to solve the asymptotic system is given in the appendix by Algorithm \ref{alg:AsymptoticModel}.

Within this work we will consider a mixture with one anion and one cation species with molar densities $n_-$ and $n_+$ and charges numbers $z_- = -1$ and $z_+=+1$, respectively. 

In the following studies we will apply all three different variations of the PNPS system by different choices of the volume fractions $a_\pm=0$, $a_\pm=1$ and $a_\pm > 1$.

\subsection{Electro-osmotic flow transitions arising in cylindrical nanopores}\label{sec:Cylinder}
In order to better understand ion transport through nanopores, whether for improving existing technologies or developing new therapies, it is crucial to identify the driving forces behind EOF and understand how these change when, for example, an additional hydrodynamic pressure gradient is applied. We will first consider a cylindrical pore with constant radius $R(z) = R$ so that no additional effects due to geometry play a role. 
For example, adding a constriction to the pore would lead to an additional pressure drop. This means that in a cylindrical pore, the influence of an applied pressure difference can be considered in isolation, whereas in a pore with variable $R(z)$, it is not clear what additional effects may occur. Moreover, some of the equations in the asymptotic model for a cylindrical pore can be further simplified.
In what follows, we assume for cylindrical pore with a constant radius $R(z) = R = 5$ and a length of $L=25$ with $\delta = 0.1$. For the surface charge we define
\begin{align}\label{eq:homCharge}
    \sigma(z) = \frac{\sigma_0}{2} (\tanh((z-L_1)/\eps) - \tanh((z-L_2)/\eps)),
\end{align}
with $\sigma_0 = 0.15$, $L_1 = 0.2L$, $L_2 = L-L_1$ and $\eps = 5$. All other parameter values are given in Table \ref{tab:ParameterValuesGeneral} and Table \ref{tab:ParameterValuesCylinder}.

\subsubsection{Transitions in the flow regime}
\begin{figure}[H]
	\centering
	\begin{subfigure}{0.44\textwidth}
		\centering
		\includegraphics[width=1\textwidth]{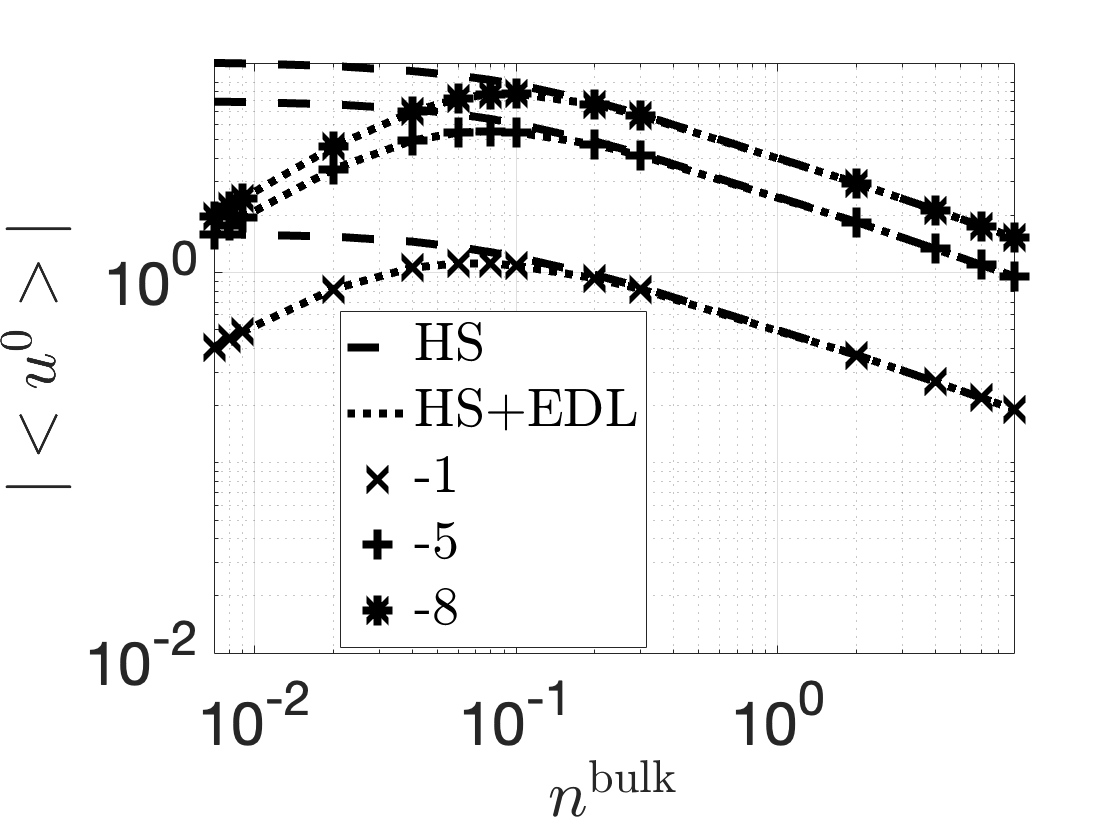}
		 \caption{}
		\label{fig:CylinderVelocityConcentrations}
	\end{subfigure}
	\hfill
    \begin{subfigure}{0.44\textwidth}
		\centering
		\includegraphics[width=1\textwidth]{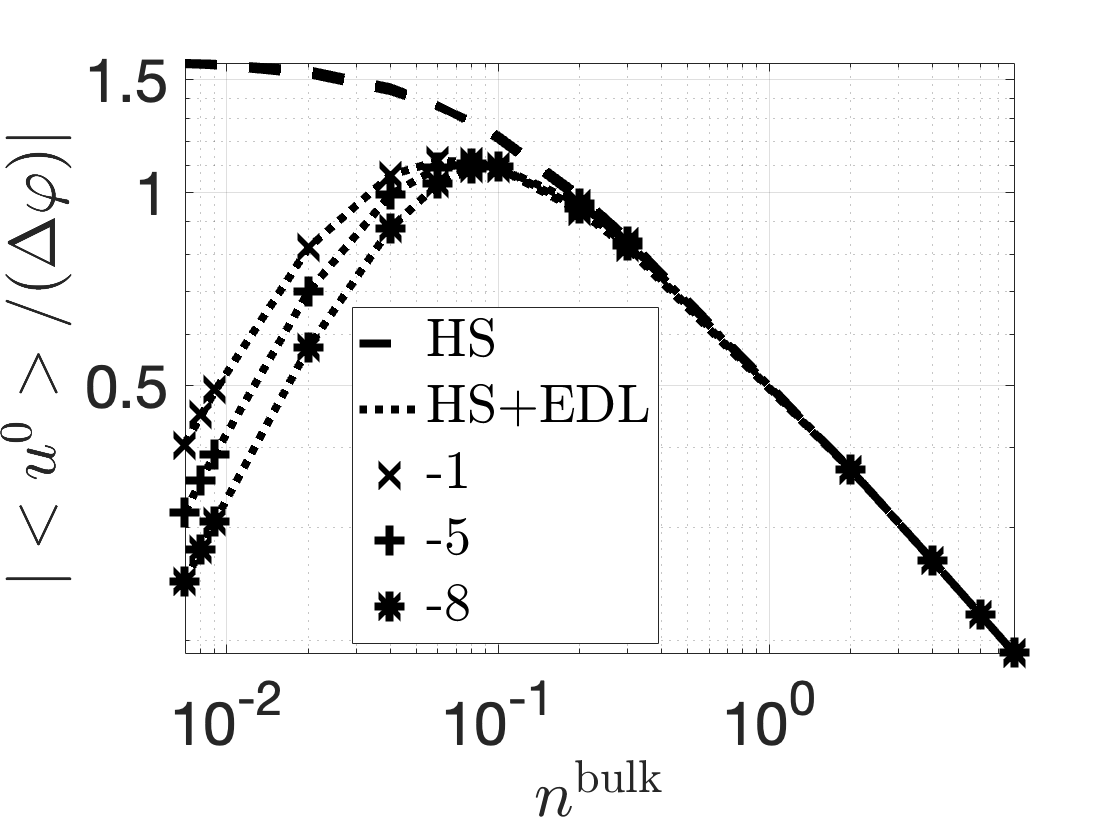}
		 \caption{}
		\label{fig:CylinderVelocityConcentrationsScaled}
	\end{subfigure}
	\caption{(a) Volumetric flow $\langle u^0\rangle$ plotted as a function of the bulk concentration $n^\tbulk$ for $a_\pm = 0$, $\Delta p = 0$ and different potentials $\Delta \varphi = -1$ (crosses), $\Delta \varphi = -5$ (plus signs) and $\Delta \varphi = -8$ (stars). The dashed line corresponds to the $\langle u^{\text{HS}}\rangle$ (Helmholtz-Smoluchowski) flow regime. The dotted line gives $\langle u^{\text{HS}}\rangle + \langle u^{\text{EDL}} \rangle $. (b) Volumetric flow scaled with $\Delta \varphi$. }
	\label{fig:VelocityConcentrations}
\end{figure}

\begin{figure}[H]
	\centering
	\begin{subfigure}{0.44\textwidth}
		\centering
		\includegraphics[width=1\textwidth]{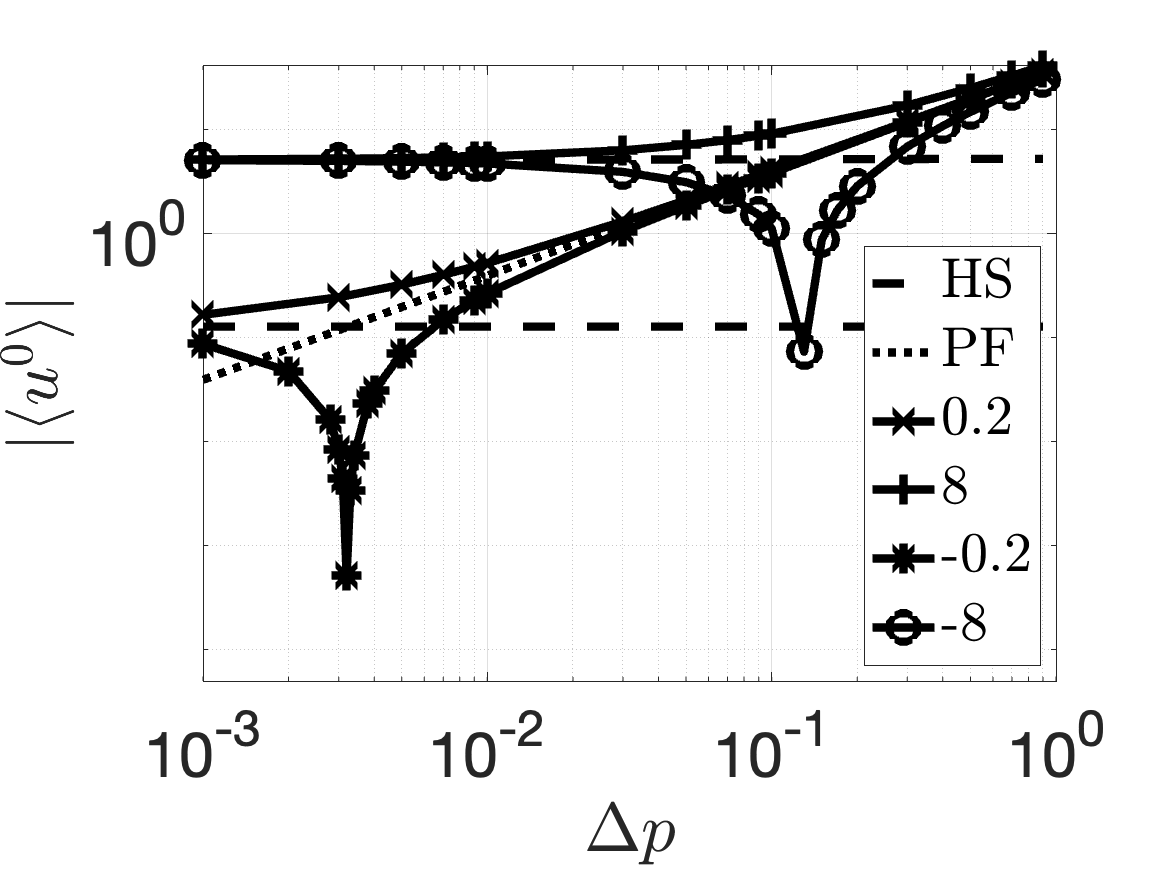}
		 \caption{}
		\label{fig:CylinderVelocity}
	\end{subfigure}
	\hfill
    \begin{subfigure}{0.44\textwidth}
		\centering
		\includegraphics[width=1\textwidth]{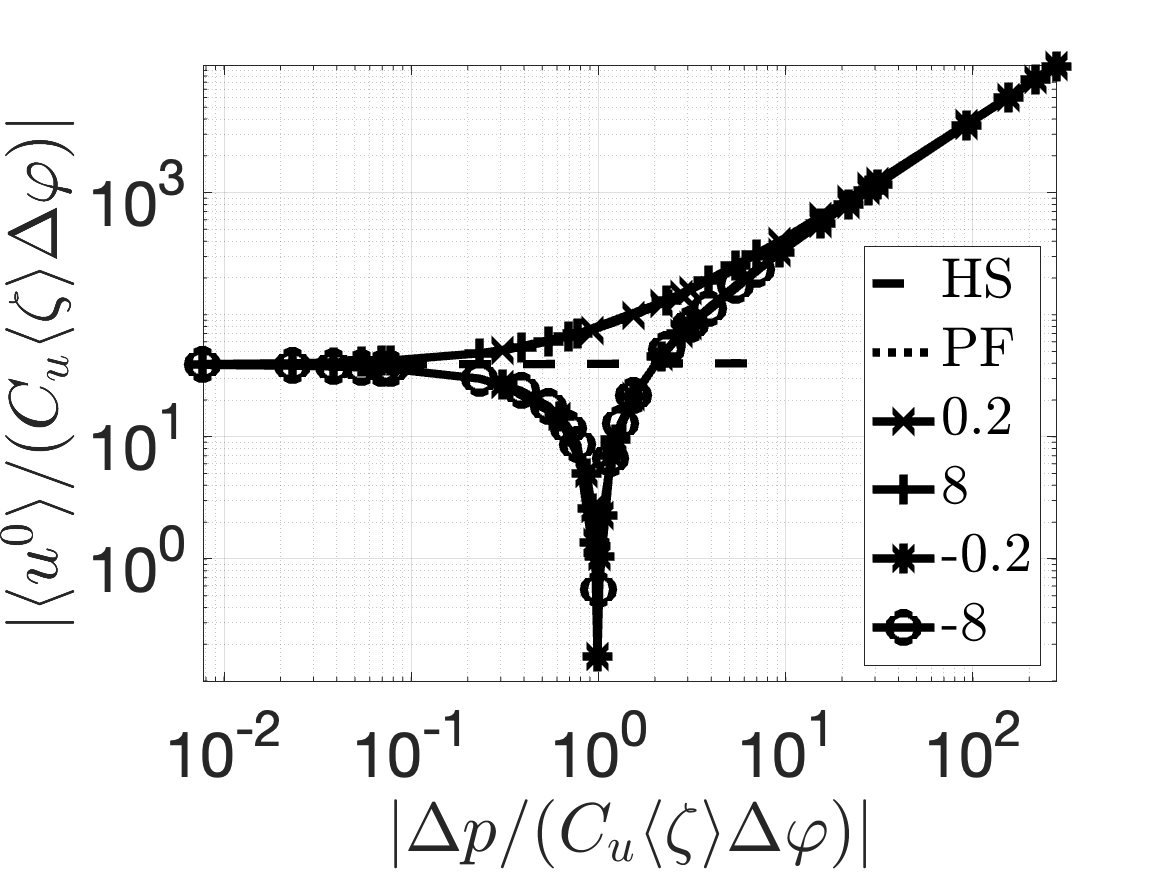}
		 \caption{}
		\label{fig:CylinderVelocityScaled}
	\end{subfigure}
	\caption{(a) Volumetric flow $\langle u^0\rangle$ for $n^\tbulk = 0.6$ ($a_\pm = 0$) plotted as a function of pressure $\Delta p$ for different potentials $\Delta \varphi = 0.2$ (crosses), $\Delta \varphi = 8$ (plus signs), $\Delta \varphi = -0.2$ (stars) and $\Delta \varphi = -8$ (circles). The dotted and dashed lines correspond to the $\langle u^{\text{PF}}\rangle$ (Poiseuille) and $\langle u^{\text{HS}}\rangle$ (Helmholtz-Smoluchowski) flow regimes, respectively. (b) Volumetric flow and pressure scaled with $C_u \Delta \varphi$. }
	\label{fig:Velocity}
\end{figure}

We calculated the volumetric flow $\langle u^0 \rangle = \int_0^L \int_0^R u^0 r \ud r \ud z$ in the pore for different salt concentrations $n_\pm^\tin = n_\pm^\tout = n^\tbulk$ (Figure \ref{fig:VelocityConcentrations}), applied pressure differences $\Delta p = p^\tin - p^\tout$ and potential differences $\Delta \varphi = \varphi^\tin - \varphi^\tout$ between in- and outlet (Figure \ref{fig:Velocity}). 
If we integrate equation \eqref{SolutionU} we get
\begin{align}\label{eq:ApproxVelocity}
    \langle  u^0 \rangle = \langle  u^{\text{PF}}\rangle + \langle  u^{\text{HS}} \rangle + \langle  u^{\text{EDL}} \rangle = \Delta p \frac{1}{16} R^4 -\Lambda^2 \frac{\Delta \varphi}{L} \langle  \zeta \rangle + \langle  u^{\text{EDL}} \rangle.
\end{align}
As illustrated by Figure \ref{fig:CylinderVelocity} the flow can be distinguished between two known regimes: a pressure-controlled regime characterized by a Poiseuille-like profile (PF)
\begin{align}\label{eq:PF}
    \langle  u^{\text{PF}} \rangle = \Delta p \frac{1}{16} R^4
\end{align}
and a potential-driven (EOF) regime. The later is characterized by a Helmholtz-Smoluchowski-like profile (HS) \cite{park2009extension}, also referred to as electroosmotic velocity
\begin{align}
    \langle  u^{\text{HS}}\rangle  = -\Lambda^2 \frac{\Delta \varphi}{L} \langle  \zeta \rangle,
\end{align}
where $\zeta=\varphi_r^0 -\varphi_r^0|_{r=R}$ is the zeta-potential. 
Additionally, a third term $\langle u^{\text{EDL}} \rangle$ with
\begin{align}
    \langle u^{\text{EDL}} \rangle = \int_0^L \int_0^R (K_2(r,z) - K_2(r=R,z))r \ud r \ud z
\end{align}
plays a role as shown in Figure \ref{fig:CylinderVelocityConcentrations}. The function $K_2$ is give by equation \eqref{eq:FinalK2General} and is a flow contribution arising from the EDL axial force density. Thus, the $ \langle u^{\text{EDL}} \rangle$ gives us a correction for the flow due to large Debye length, i.e. for low concentrations with weak screening of the surface charge such that $\Lambda = \mathcal{O}(1)$.
For a detailed computation see Appendix \ref{sec:AppendixMeanVelocity}. 
However, we find that for a salt concentration high enough the EDL--term does not contribute significantly as illustrated in Figure \ref{fig:CylinderVelocityConcentrations}. For this example we did not apply a pressure difference, i.e. $\Delta p = 0$, such that $\langle u^{\text{PF}}\rangle = 0$. We then calculated the volumetric flow for different salt concentrations $n^\tbulk$ with $a_\pm = 0$ (i.e. the classical PNPS system) and applied potential differences. It can be observed that for all $\Delta \phi$, the flow initially increases with increasing concentration until a maximum is reached, and then decreases again. The maximum varies for the different potentials. In addition, we see that after the peak the flow follows the Helmholtz-Smoluchowski-like profile (dashed lines). 
For this particular example and set of parameters our studies showed similar behavior for $a_\pm > 0$. However, this might not always be the case since finite-volume effects might also impact the $\zeta$-potential since it also depends on the concentrations. 
For small concentrations we find that $\langle u^{\text{HS}}\rangle$ is almost constant and that the volumetric flow is driven by the sum of the HS--term and the EDL--term. If we scale flow with $\Delta \varphi $ it can be observed that for all potentials the maxima become equal. However, we still notice a rightward (or leftward) shift of the curves for small concentrations (Figure \ref{fig:CylinderVelocityConcentrationsScaled}).
For the following study we chose the concentrations $n_\pm^\tout=n_\pm^\tin=n^\tbulk = 0.6$ such that $\langle u^{\text{EDL}} \rangle \ll \langle u^{\text{HS}} \rangle$ (see Figure \ref{fig:CylinderVelocityConcentrationsScaled}) and the flow is given by $\langle u^0\rangle = \langle u^{\text{PF}}\rangle + \langle u^{\text{HS}}\rangle$. We then calculated the volumetric flow for different applied pressure and potential differences to study the transition between the pressure-controlled and the potential-controlled regimes. 
We find that if we apply no potential difference ($\Delta \varphi = 0$) that the volumetric flow depends linearly on the pressure difference $\Delta p$ (Figure \ref{fig:CylinderVelocity} dotted line). On the other hand, if we apply zero pressure, we find that the velocity profile is constant which is illustrated in Figure \ref{fig:CylinderVelocity} (dashed lines). 
Furthermore, we find for $\langle u^0 \rangle = 0$ that $\Delta p - C_u \langle \zeta \rangle \Delta \phi = 0$, with $C_u = \Lambda^2 16/(R^4 L)$, such that
\begin{align}
    \langle u^0 \rangle \begin{cases}
        \leq 0 \quad \text{if } \Delta p \leq C_u \langle \zeta \rangle \Delta \varphi \quad &\text{for } \langle \zeta \rangle \Delta \varphi \geq 0\\ 
        \geq 0 \quad \text{if } \Delta p \geq C_u \langle \zeta \rangle \Delta \varphi\quad &\text{for } \langle \zeta \rangle \Delta \varphi \geq 0 \\
        \geq  0 \quad \forall \Delta p \quad &\text{for } \langle \zeta \rangle \Delta \varphi \leq  0.
    \end{cases}
\end{align}
Depending on the sign of the applied potential and pressure difference as well as on the sign of the zeta-potential the volumetric flow changes its sign. For our example we calculated $\langle \zeta \rangle = -88$ and $C_u = 1.8e-4$ for all applied potential differences. As illustrated in Figure \ref{fig:CylinderVelocity} we find that for negative potentials the absolute value of the flow first decreases with increasing pressure until a minimum is reached and then starts to increase again for increasing pressure. The minimum indicates the change in sign as shown in Figure \ref{fig:CylinderVelocityScaled}. For $\Delta \phi = -0.2$ the change in sign is at $\Delta p = 0.0032$ (stars) and for $\Delta \phi = -8$ the change in sign is at $\Delta p = 0.13$ (circles). For $\Delta \phi = 0.2$ (crosses) and $\Delta \phi = 8$ (plus signs), on the other hand, the flow rate increases steadily as the pressure increases.

Figure \ref{fig:CylinderVelocityScaled} shows the volumetric flow scaled by a factor of $C_u \langle \zeta \rangle \Delta \phi$. While the unscaled data show clear quantitative differences between the applied potentials, rescaling the velocities by a factor of $C_u \langle \zeta \rangle \Delta \phi$ results in a clear, uniform behavior. When plotted as normalized flow, all curves for negative potentials converge to a single master curve still reflecting the change in sign in the flow. Similarly, all curves for positive potentials merge into a single curve. This convergence suggests that the influence of the potential difference can be incorporated into this scaling relationship, indicating universal transport behavior once the appropriate dimensionless scale is applied. This rescaling therefore provides a robust way to compare results across different potential landscapes and highlights the underlying physical mechanism that controls the flow. Furthermore, as illustrated in Figure \ref{fig:CylinderVelocityScaled} the flow is dominated by the HS regime for low pressure values until it continuously transitions into the PF regime for increasing pressure.
We were able to show that, even though we use a general derivation for the asymptotic model without further simplifications such as small Debye length or a constant electric field (as in \cite{park2009extension}), we can reconstruct two known solutions for the volumetric flow: Poiseuille and Helmholtz-Smoluchowski. While most studies only consider one of these two regimes, we show that the volumetric flow changes between them depending on the applied voltage and pressure differences. Moreover, our model shows that for small salt concentrations there is an additional contribution to the volumetric flow. 
Note that although the dependence of the flow on $\Delta p$ and $\Delta \varphi$ is evident from equation \eqref{eq:ApproxVelocity}, the zeta-potential still has to be calculated numerically. Since we observe in our numerical simulations that it does not change dramatically for different potentials and pressures, it only needs to be calculated once for a specific salt concentration and surface charge.  
\subsubsection{Transitions in the ion current}
\begin{figure}[H]
	\centering
	\begin{subfigure}{0.44\textwidth}
		\centering
		\includegraphics[width=1\textwidth]{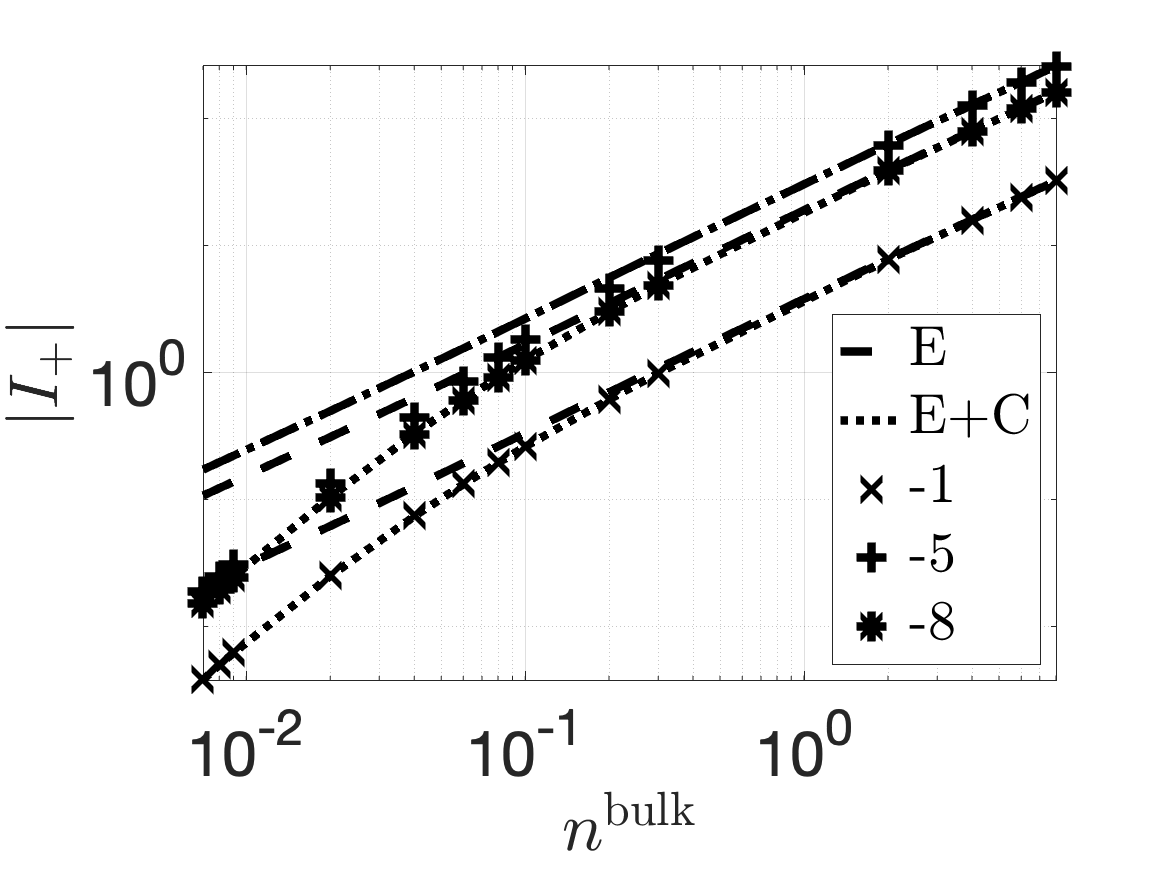}
		 \caption{}
		\label{fig:CylinderCurrentConcentration}
	\end{subfigure}
	\hfill
    \begin{subfigure}{0.44\textwidth}
		\centering
		\includegraphics[width=1\textwidth]{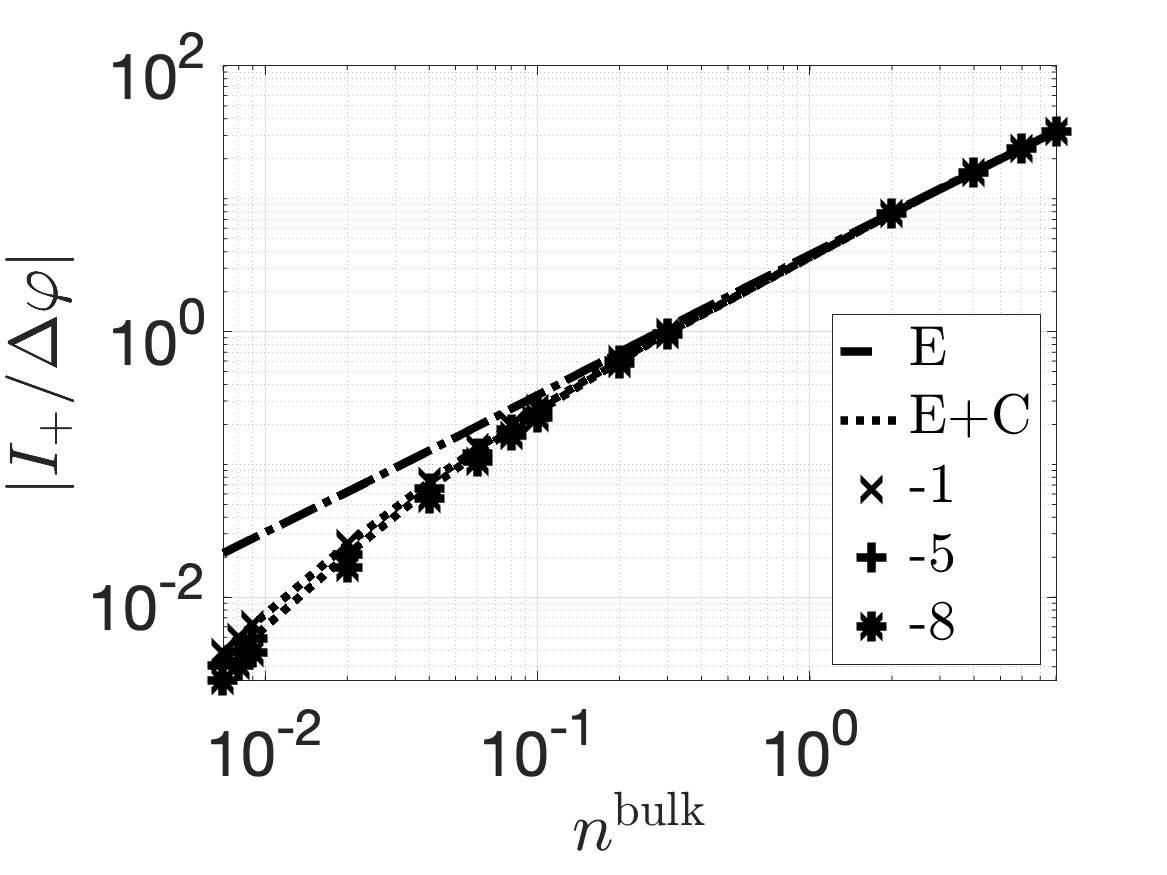}
		 \caption{}
		\label{fig:CylinderCurrentConcentrationScaled}
	\end{subfigure}
	\caption{(a) Cation current $I_+$ plotted as a function of the bulk concentration $n^\tbulk$ for $a_\pm = 0$, $\Delta p = 0$ and different potentials $\Delta \varphi = -1$ (crosses), $\Delta \varphi = -5$ (plus signs) and $\Delta \varphi = -8$ (stars). The dashed line corresponds to the $I^\text{E}$ (potential-driven) regime. The dotted line gives $I^\text{E} + I^\text{C}$. (b) Cation current scaled with $C^\tbulk \Delta \varphi$. }
	\label{fig:CurrentConcentration}
\end{figure}

We also calculated the concentration and pressure dependency of the ionic currents flowing through a circular area at $z=0$. Taking a closer look on equation \eqref{eq:FinalCurrentGeneral} we find that the currents can be approximated as follows
\begin{align}
    I_\pm |_{z=0} &= 2 \pi n^\tbulk_\pm \frac{\Delta \varphi}{L} \left( \frac{k_\pm }{2} R^2  + P_e \Lambda^2 \frac{\langle \zeta\rangle}{L} \right)  \mp  2 \pi P_e n^\tbulk_\pm \frac{\Delta p}{L} \frac{1}{16} R^4 \\ \nonumber
    & \quad \pm 2\pi \bar n^\tbulk \left( k_\pm  (y_0^\tbulk)^{a_\pm} \frac{1}{2} R^2 \partial_z Q_\pm |_{z=0} \mp P_e \frac{\langle u^{\text{EDL}} \rangle}{L} \right) \\ \nonumber
    &= I_\pm^{\text{E}} + I_\pm^{\text{P}} + I_\pm^{\text{C}},
\end{align}
with $n_\pm^\tout = n_\pm^\tin = n_\pm^\tbulk$. A detailed computation can be found in the Appendix \ref{eq:ApproxCurrent}. 
Overall, the current can be written as the sum of a potential-driven (EOF) current $I^{\text{E}}$, a pressure-driven current $I^{\text{P}}$ and a third concentration-dependent term $I^{\text{C}}$ giving us a large Debye layer-correction for the current. 
Similar as for the volumetric flow we now find that $I^{\text{C}}$--term only plays a role for small concentrations when the Debye layer is large (see Figure \ref{fig:CurrentConcentration}). Furthermore, we also found for the current, that choosing $a_\pm = 0 > 0$ gives similar results. We therefore chose $a_\pm = 0$ for this study.
For all three applied voltages $\Delta \varphi = -1$ (crosses), $\Delta \varphi = -5$ (plus signs) and $\Delta \varphi = -8$ (stars), the measured current rise continuously with increasing concentrations  (Figure \ref{fig:CylinderCurrentConcentration}). 

\begin{figure}[H]
	\centering
	\begin{subfigure}{0.44\textwidth}
		\centering
		\includegraphics[width=1\textwidth]{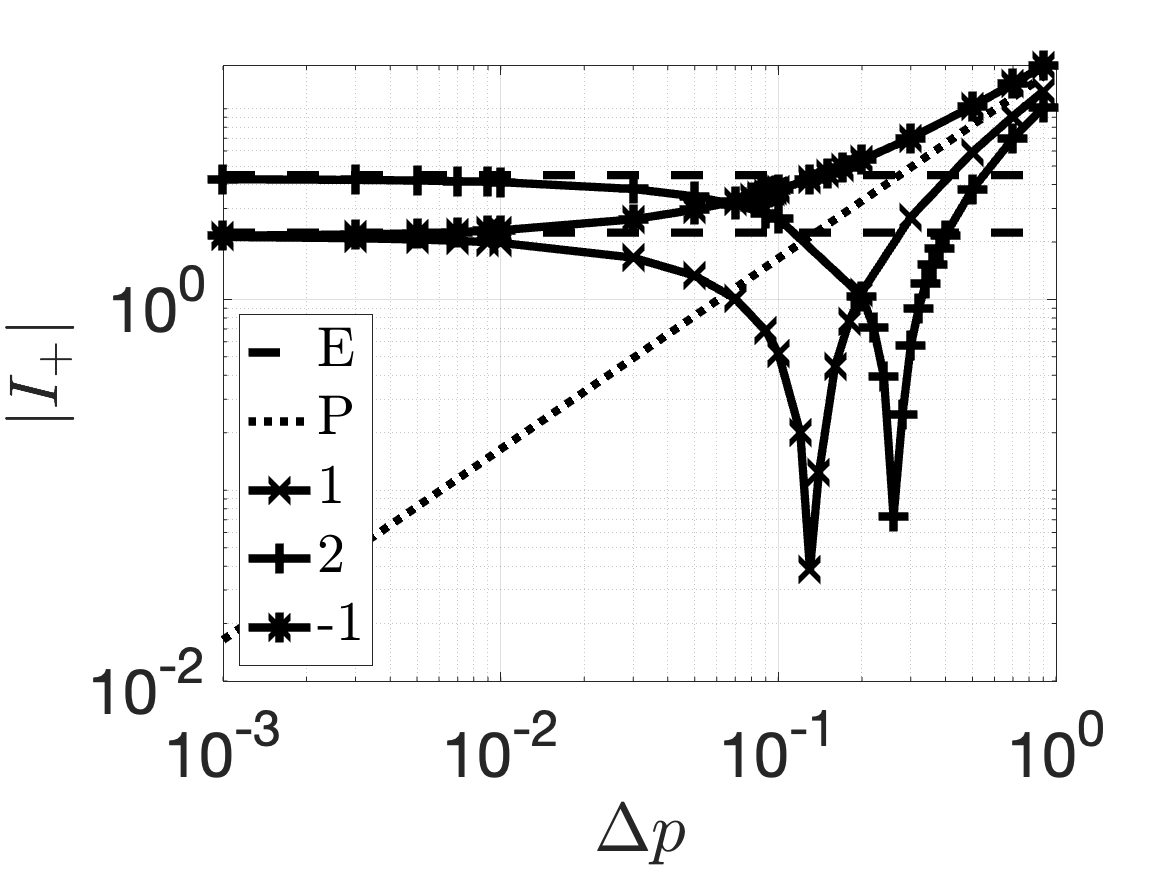}
		 \caption{}
		\label{fig:CylinderCurrent}
	\end{subfigure}
	\hfill
    \begin{subfigure}{0.44\textwidth}
		\centering
		\includegraphics[width=1\textwidth]{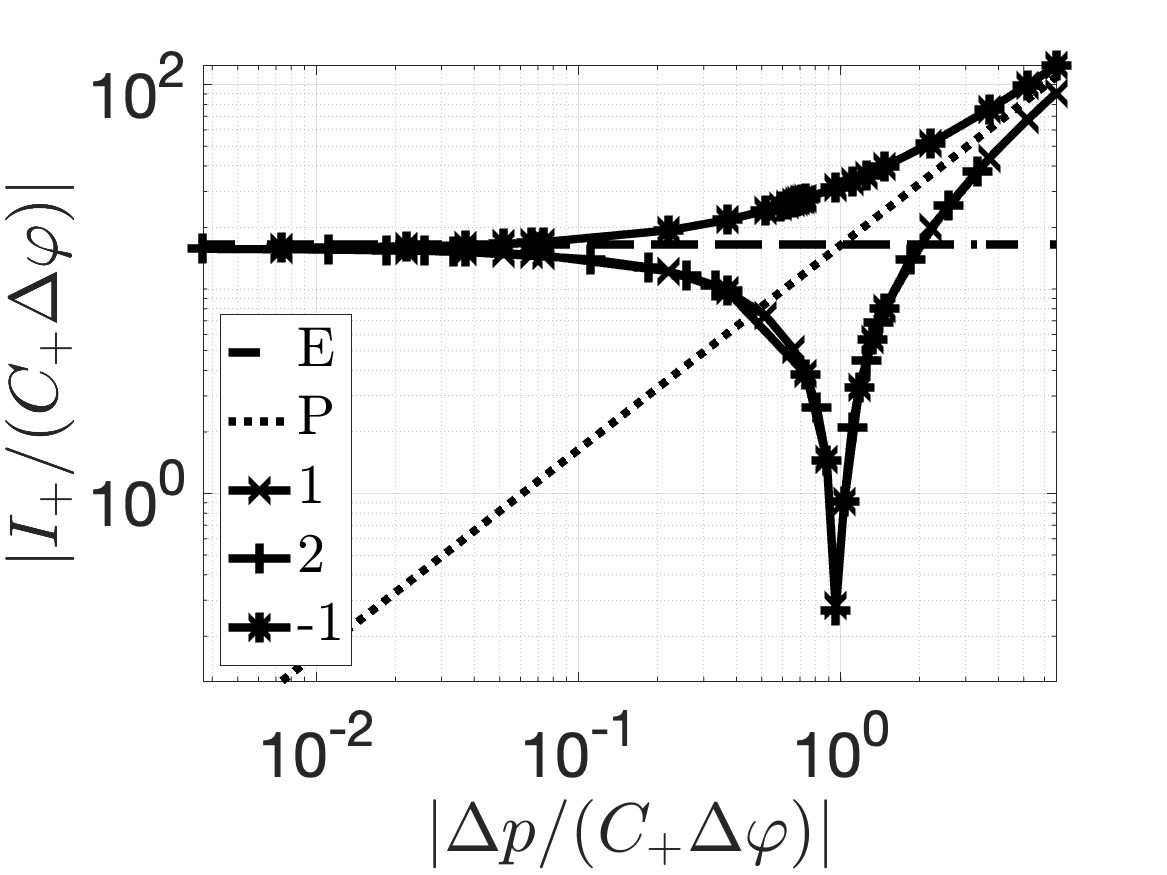}
		 \caption{}
		\label{fig:CylinderCurrentScaled}
	\end{subfigure}
	\caption{(a) Cation current $I_+$ for $n^\tbulk = 0.6$ ($a_\pm = 0$) plotted as a function of pressure $\Delta p$ for different potentials $\Delta \varphi = 1$ (crosses), $\Delta \varphi = 2$ (plus signs), $\Delta \varphi = -1$ (stars). The dotted and dashed lines correspond to the $I^{\text{P}}$ (pressure-driven) and $I^{\text{E}}$ (potential-driven) current regimes, respectively. (b) Current and pressure scaled with $C_+ \Delta \varphi$.}
	\label{fig:Current}
\end{figure}

Scaling the current with $\Delta \phi$ as given in Figure \ref{fig:CylinderCurrentConcentrationScaled} shows that in all three cases the current continuously transitions into the $I^{\text{E}}$-curve (dashed line) indicating that the potential-driven regime dominates there. Note that $I/\Delta \phi$ is often also reffered to as conductance. For small concentrations it can be observed that the $I^{\text{C}}$-term is negligible.
Furthermore, we find that with that scaling all three curves collapse into one curve for high salt concentrations while there is still a slight rightward (or leftward) shift for small concentrations. 
At sufficiently high concentrations, we therefore find that the current consists of a pressure-driven and a potential-driven regime.
For concentrations where $I_\pm^{\text{C}} \ll I_\pm^{\text{E}}$ (see Figure \ref{fig:CylinderCurrentConcentrationScaled}) it can additionally be observed that there is a switch in the sign of the cation current depending on the applied potential and pressure. For $I_\pm = 0$ we have $\Delta p \mp C_\pm \Delta \phi = 0$, with $C_\pm = k_+ 8/(P_e R^2) + \Lambda^2 \langle \zeta\rangle 16/(L R^4) $, such that
\begin{align}\label{eq:ApproxCurrentCation}
    I_+ \begin{cases}
        \leq 0 \quad \text{if } \Delta p \leq C_+ \Delta \varphi \quad &\text{for } C_+ \Delta \varphi \geq 0\\ 
        \geq 0 \quad \text{if } \Delta p \geq C_+ \Delta \varphi\quad &\text{for } C_+ \Delta \varphi \geq 0 \\
        \geq  0 \quad \forall \Delta p \quad &\text{for } C_+ \Delta \varphi \leq  0
    \end{cases}
\end{align}
and 
\begin{align}\label{eq:ApproxCurrentAnion}
    I_- \begin{cases}
        \leq 0 \quad \text{if } \Delta p \leq C_- \Delta \varphi \quad &\text{for } C_- \Delta \varphi \leq 0\\ 
        \geq 0 \quad \text{if } \Delta p \geq C_- \Delta \varphi\quad &\text{for } C_- \Delta \varphi \leq 0 \\
        \geq  0 \quad \forall \Delta p \quad &\text{for } C_- \Delta \varphi \geq  0.
    \end{cases}
\end{align}
We calculated the absolute cation current for $n^\tbulk = 0.6$ as a function of $\Delta p$ for different potential differences $\Delta \varphi = 1$, $\Delta \varphi = 2$ and $\Delta \varphi = -1$ (Figure \ref{fig:CylinderCurrent}). For this example we calculate $\langle \zeta\rangle=-88$ and $n_+= 0.14$.
For positive potentials it can be observed, as illustrated in Figure \ref{fig:CylinderCurrentScaled}, that for increasing pressure, the cation current first decreases, reaches a minimum, and then increases. The respective change in sign here is at approximately $\Delta p = 0.14$ (for $\Delta \varphi = 1$) and $\Delta p = 0.28$ (for $\Delta \varphi = 2$) (Figure \ref{fig:CylinderCurrent}). For $\Delta \varphi = -1$, however, we do not observe any reversal of the current. In this case, the current increases continuously with increasing pressure.
Similar to the flow, we find that when the current is rescaled by a factor of $C_\pm \Delta \phi$, all curves for positive potentials converge to a single curve that continues to capture the sign change for the current. By applying a suitable dimensionless scaling, a universal behavior could thus also be derived for the current.  
For all three cases we find that at low pressure values, the current is almost constant and that the potential-driven current $I^{\text{E}}$ (dashed line) is dominant. As the pressure increases, the current transits into the pressure-driven $I^{\text{P}}$ (dotted line) regime.
Our results show that, in theory, it would be possible to push positively charged ions against their electrostatic gradient. Equations \eqref{eq:ApproxCurrentCation} and \eqref{eq:ApproxCurrentAnion} show that the ion current through the pore can also be differentiated into two regimes. We see that the current can not only be amplified by a hydrodynamic pressure difference, but can also overcome the electrostatic forces.

In \cite{Curk2024}, further regimes and transitions have been discovered, due to the effects of hydrodynamic boundary slip.
In this work we did not take into account transitions triggered by hydrodynamic interfacial slip. Since these are important for a systematic understanding, they will be the subject of our future studies.
\subsection{Impact of geometry, fluid flow and finite-size effects}\label{sec:Trumpet}
\begin{figure}[H]
    \centering
    \begin{tikzpicture}[scale=0.5]

    \draw[->,line width=0.25mm] (-0.5,0) -- (20,0) node[right] {$z$}; 
    \draw[->,line width=0.25mm] (0,-0.5) -- (0,8) node[above] {$r$}; 
    
    \draw[line width =0.5mm,domain=0:17,samples=100] plot (\x,{4*(7-4)*((\x)^2 - \x*17)/17^2 + 7}) 
        node[right] at (5.5,5) {$-\partial_r \varphi_r = \gamma \sigma(z)$}; 

    \draw[line width =0.5mm] (0,7) -- (0,0) node[left] at (-2.5,6){outlet $S^\tout$};
    \node[ ] at (-3.7, 4.5) {$n^R n^\tout = 0.1$\,\si{\mol \per \liter}};
    \node[ ] at (-4.7, 3) {$\phi^R \phi^\tout = 0.0$\,\si{\volt}};
    \node[ ] at (-4.5, 1.5) {$p^R p^\tout = 0.0$\,\si{\pascal}};
    
    \draw[line width =0.5mm,dashed] (0,0) -- (17,0) node[right] at (5.5,-1) {$-\partial_r \phi_r = 0$};
    
    \draw[line width =0.5mm] (17,0) -- (17,7) node[right] at (18.5,6) {inlet $S^\tin$};
    \node[ ] at (21, 4.5) {$n^R n^\tin = 0.1$\,\si{\mol \per \liter}};
    \node[ ] at (20, 3) {$\phi^R \phi^\tin = 0.2$\,\si{\volt}};
    \node[ ] at (20, 1.5) {$p^R p^\tin = 0.0$\,\si{\pascal}};

    \draw[<->,line width=0.5mm] (8.5,0) -- (8.5,4)
    node[right] at (8.7,2) {$R(z)$};
    
\end{tikzpicture}
	\caption{Illustration of a 2D axially symmetric trumpet shaped pore $\Omega$. We apply a surface charge at $r=R(z)$ and a symmetry condition at $r=0$ for the potential. At the inlet at $z=L$ and the outlet at $z=0$ we apply Dirichlet boundary conditions for the concentrations, the potential and the pressure.}
	\label{fig:Trumpet}
\end{figure}
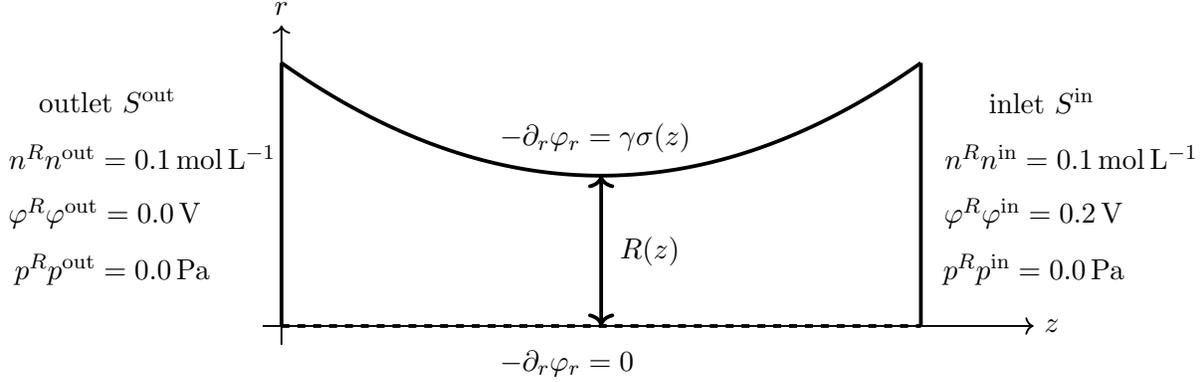

\begin{figure}[H]
	\centering
    \begin{subfigure}{0.44\textwidth}
		\centering
		\includegraphics[width=1\textwidth]{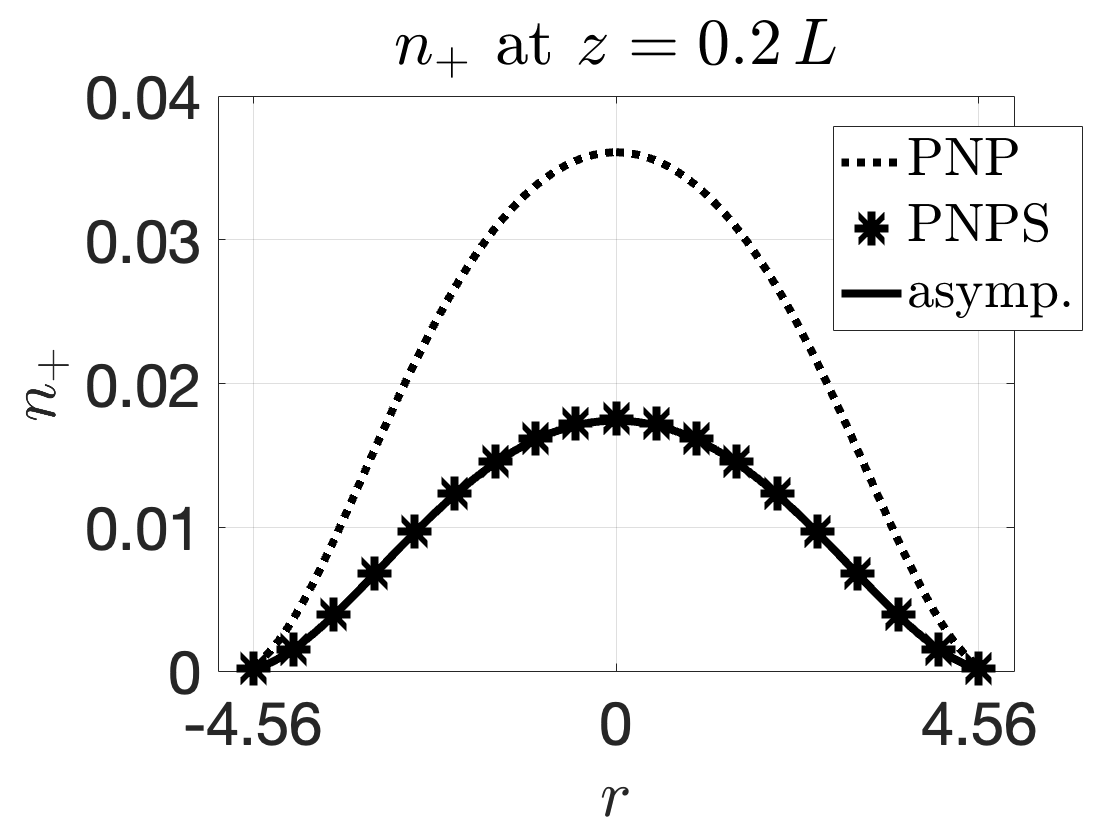}
		 \caption{}
		\label{fig:cpx200}
	\end{subfigure}
	\hfill
     \begin{subfigure}{0.44\textwidth}
		\centering
		\includegraphics[width=1\textwidth]{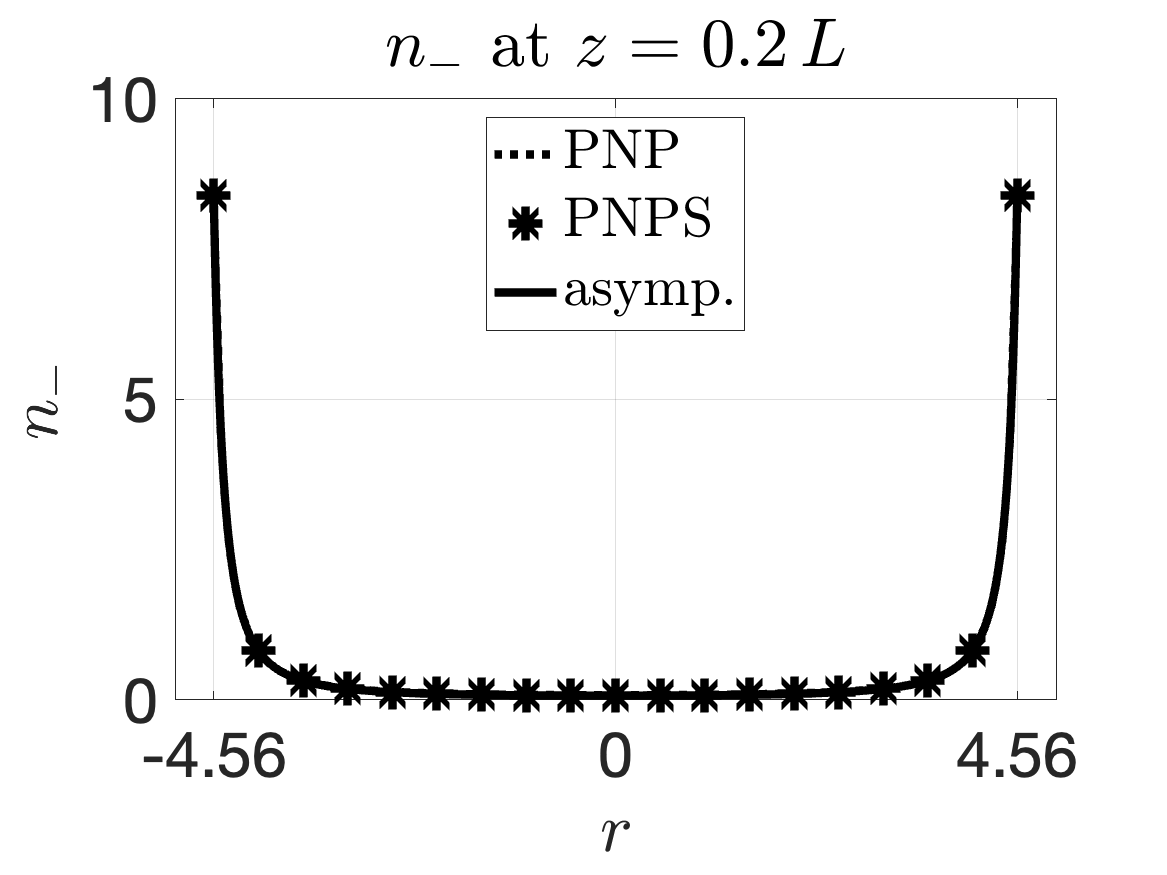}
		 \caption{}
		\label{fig:cnx200}
	\end{subfigure}
    \vfill
	\begin{subfigure}{0.44\textwidth}
		\centering
		\includegraphics[width=1\textwidth]{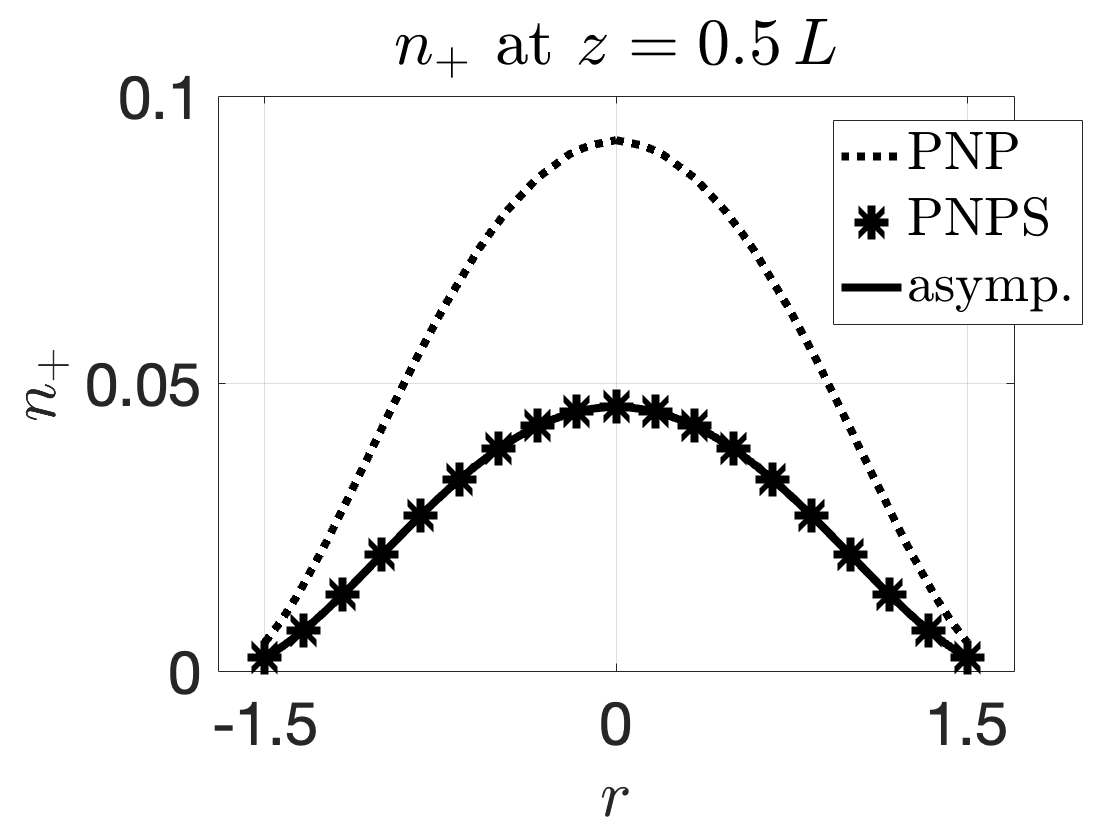}
		\caption{}
		\label{fig:cpx500}
	\end{subfigure}
	\hfill
    \begin{subfigure}{0.44\textwidth}
		\centering
		\includegraphics[width=1\textwidth]{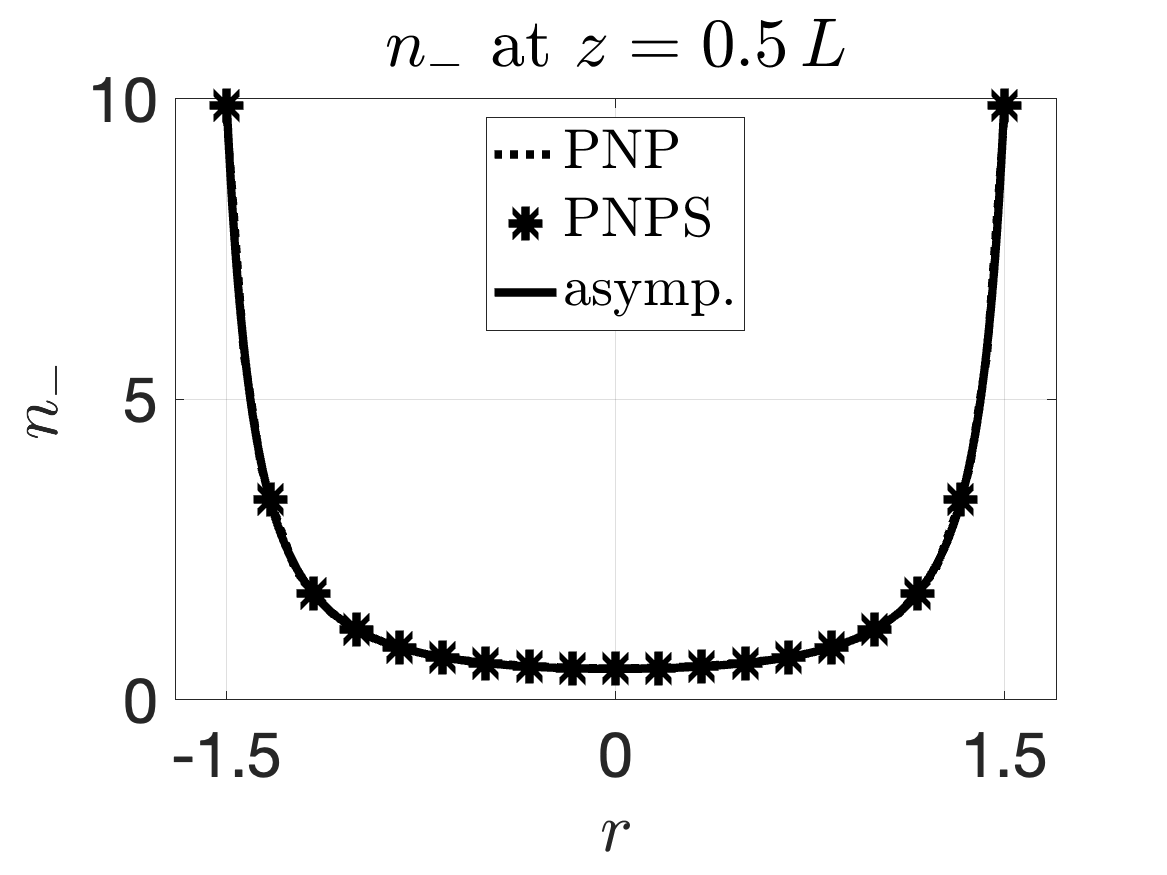}
		\caption{}
		\label{fig:cnx500}
	\end{subfigure}
    \vfill
	\begin{subfigure}{0.44\textwidth}
		\centering
		\includegraphics[width=1\textwidth]{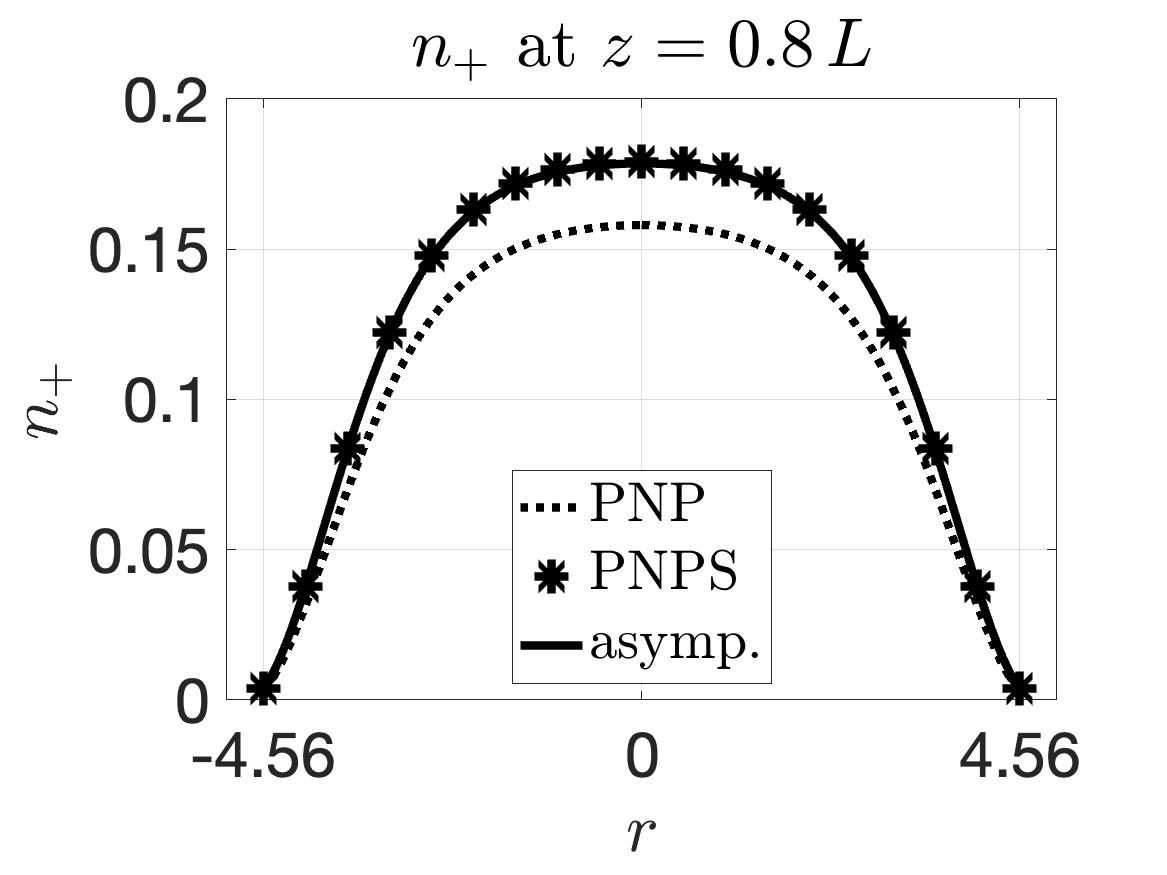}
		\caption{}
		\label{fig:cpx800}
	\end{subfigure}
	\hfill
	\begin{subfigure}{0.44\textwidth}
		\centering
		\includegraphics[width=1\textwidth]{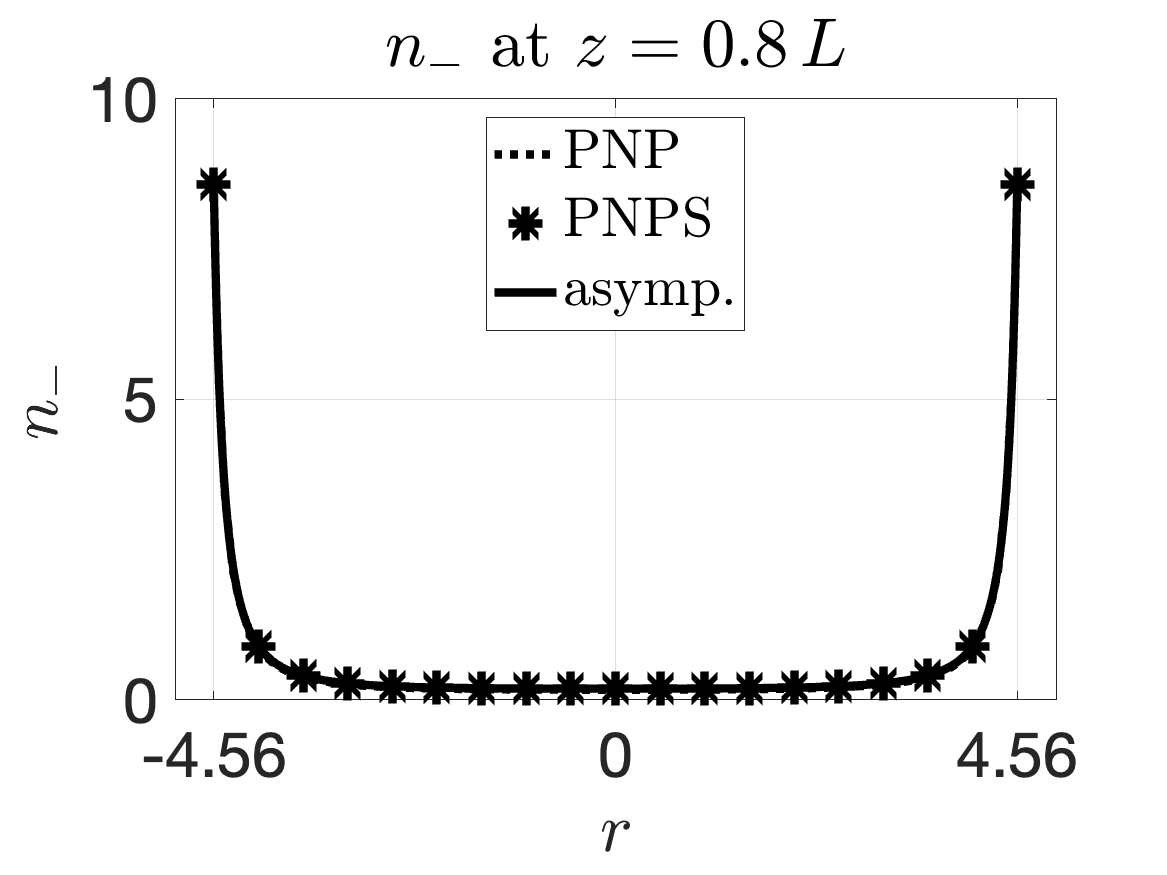}
		\caption{}
		\label{fig:cnx800}
	\end{subfigure}
    \caption{Ion concentrations for a trumpet shaped pore for two different scenarios: classical PNP equations (dotted line) and classical PNPS equations for $a_\pm = 0$ where we compare the 2D FEM solution from FEniCS (stars) to the asymptotic solution (solid line). The left column depicts the cation concentration, the right column the anion concentration.}
	\label{fig:ComparisonMatejczyk}
\end{figure}

\begin{figure}[H]
	\centering
    \begin{subfigure}{0.48\textwidth}
		\centering
		\includegraphics[width=1\textwidth]{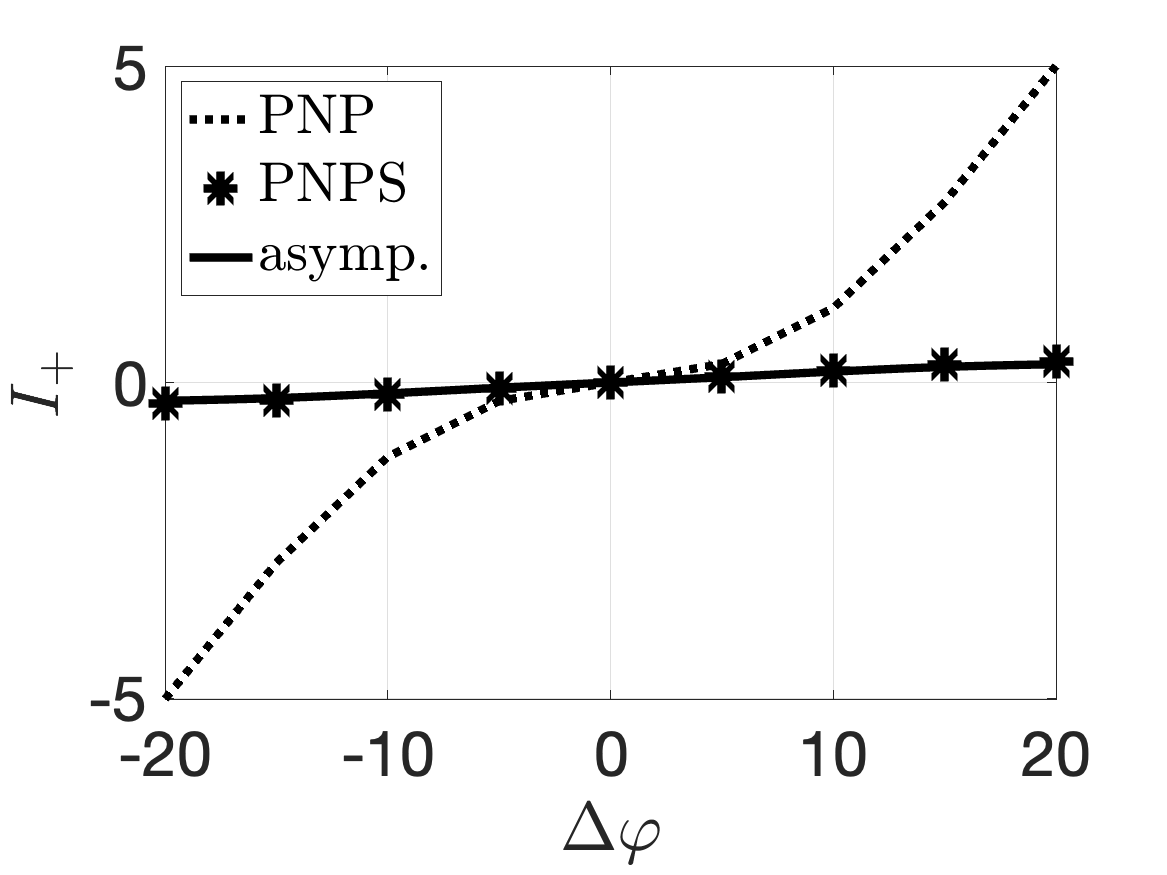}
		 \caption{}
		\label{fig:IV_cp}
	\end{subfigure}
	\hfill
	\begin{subfigure}{0.48\textwidth}
		\centering
		\includegraphics[width=1\textwidth]{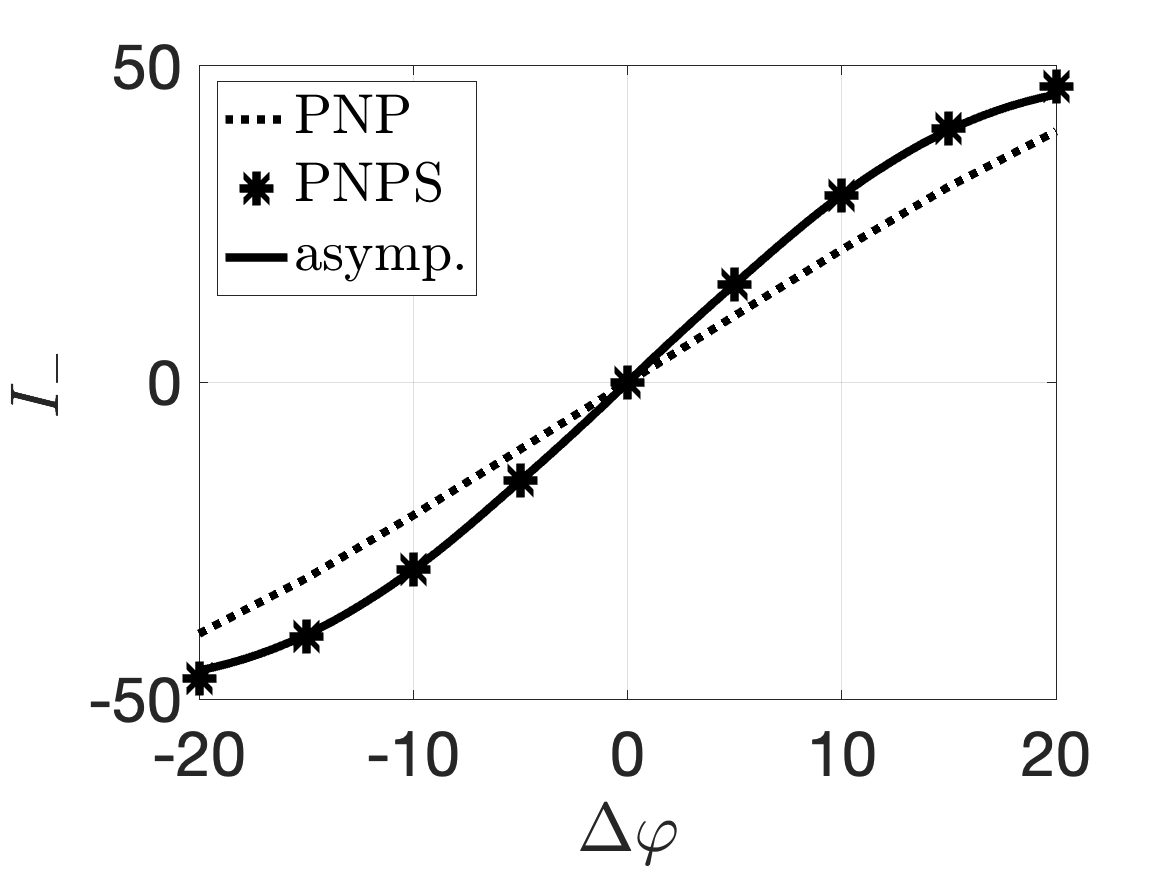}
		\caption{}
		\label{fig:IV_cn}
	\end{subfigure}
    \vfill
    \begin{subfigure}{0.48\textwidth}
		\centering
		\includegraphics[width=1\textwidth]{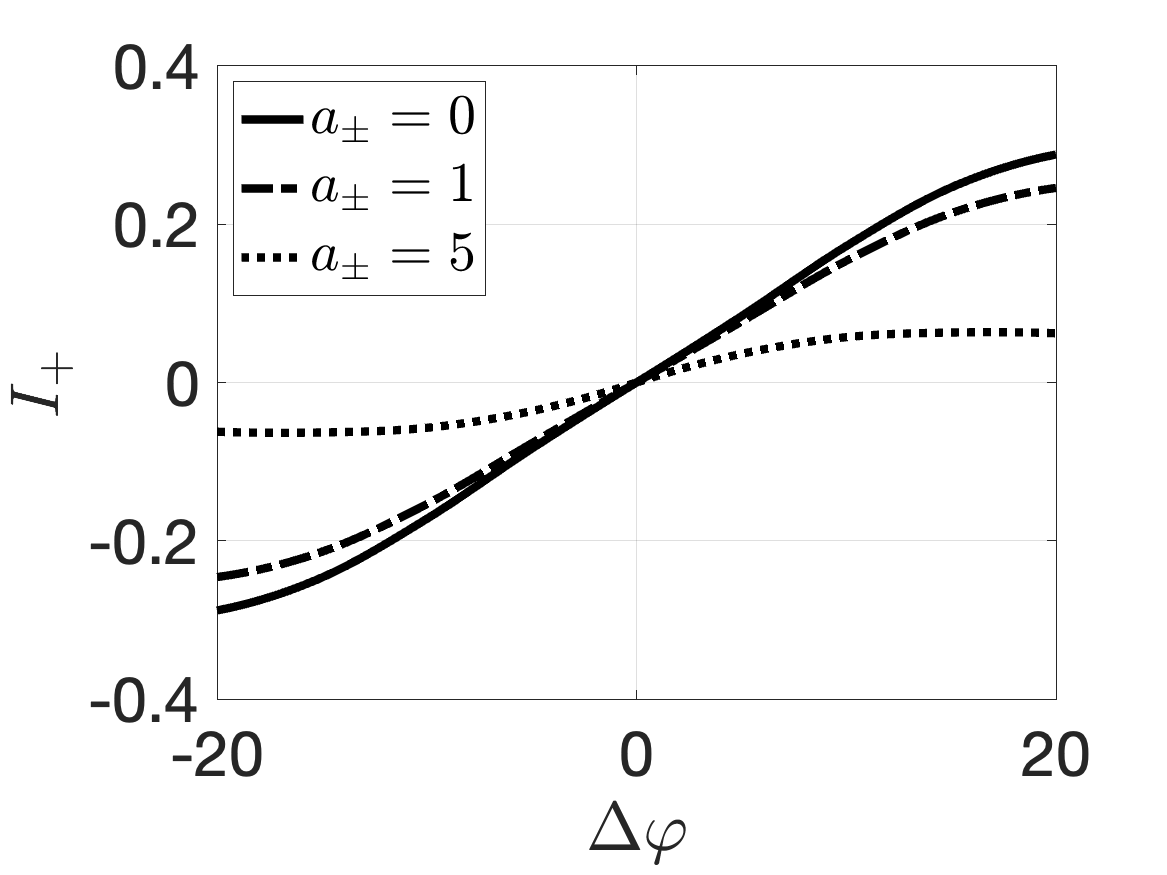}
		 \caption{}
		\label{fig:IV_cp_a}
	\end{subfigure}
	\hfill
	\begin{subfigure}{0.48\textwidth}
		\centering
		\includegraphics[width=1\textwidth]{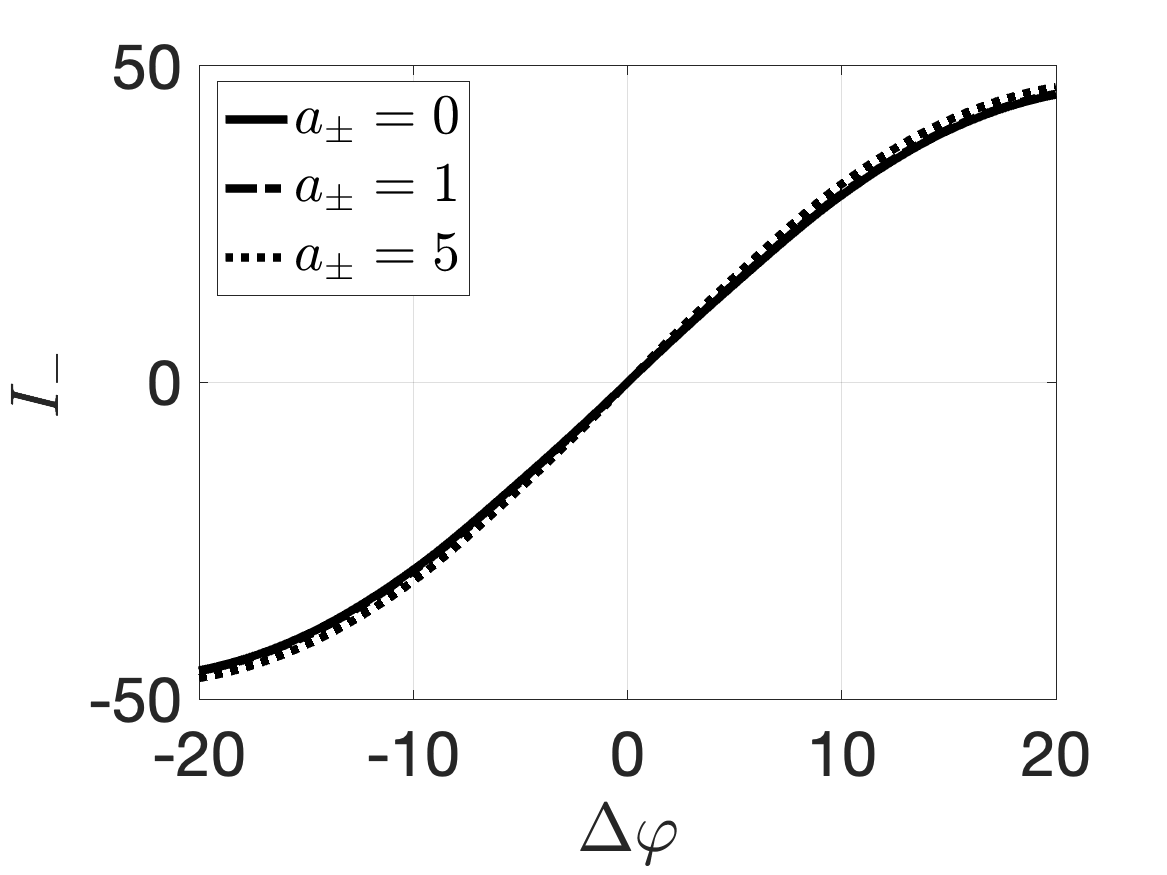}
		\caption{}
		\label{fig:IV_cn_a}
	\end{subfigure}
	\caption{IV-relation for cations (a) and anions (b) for $a_\pm = 0$. The dotted line is the solution of the classical PNP model with $P_e=0$, the stars are the 2D FEM FEniCS solution of the classical PNPS model and the solid line is the solution of the asymptotic PNPS model. (c) Cation current for different values of $a_\pm$ for the generalized PNPS model and (d) anion current for different values of $a_\pm$. The solid line is the asymptotic solution of the classical PNPS model with $a_\pm =0$, the dashed line is for $a_\pm=1$ and the dotted line is for $a_\pm=5$.}
	\label{fig:IV_trumpet}
\end{figure}

In this case study we consider the impact of geometry of the ion channel, where we restrict ourselves to axially symmetric cases. This has been investigated for many decades and the typical test cases are cylindrical, conical and trumpet shapes. To better illustrate the influence of geometry and fluid flow, we first consider the classic PNP or PNPS model, in which we choose $a_\pm=0$.
Within this subsection we thus consider a trumpet shaped pore as illustrated in Figure \ref{fig:Trumpet} with length $L=10$, $\delta = 0.01$ and radius
\begin{align}\label{eq:radiusTrumpet}
    R(z) = 4 \frac{R_1 - R_2}{L^2} (z^2 - z L) + R_1.
\end{align}
We chose $R_1 = 10$ and $R_2 = 1.5$.
For the ion concentration we apply the following Dirichlet boundary conditions $n_\pm^\tout = n_\pm^\tin = 0.1$. We consider a surface charge on the channel wall of the form \eqref{eq:homCharge} 
with $\sigma_0 = 1$, $L_1 = 0.1L$, $L_2 = L-L_1$ and $\eps = 8$.
In addition, we apply a potential difference with $\phi^\tin = 8$, s.t. $\Delta \phi =\phi^\tin- \phi^\tout = 8$. All other parameter values are given in Table \ref{tab:ParameterValuesGeneral} and Table \ref{tab:ParameterValuesTrumpet}.
We calculated the solution for two different scenarios. In the first case, we have solved the classical PNP system given by the following equations \eqref{eqn:scaled:a}-\eqref{eqn:scaled:d} with $P_e = 0$ and $a_\pm = 0$. Note that such a system was also studied by Matejczyk et al. \cite{matejczyk2018} for trumpet shaped pores. In the second case, we have solved the full problem including the coupling with Stokes equations. In both cases we calculated the steady state solution and find that, due to the positive surface charge on the channel wall, the anion concentration becomes high close to the pore wall and decreases towards the center of the pore (figures \ref{fig:cnx200}, \ref{fig:cnx500} and \ref{fig:cnx800}). The cation concentration on the other hand has a maximum on the pore center and decreases towards the wall (figures \ref{fig:cpx200}, \ref{fig:cpx500} and \ref{fig:cpx800}).
Comparing the solutions for the concentrations from the classical PNP (dotted lines) with the classical PNPS system with $P_e \neq 0$ and $a_\pm = 0$ (solid lines and stars) we find that the fluid flow has an impact, especially on the cation concentration (Figure \ref{fig:ComparisonMatejczyk}). For the cation concentration, we observe the greatest changes in the center of the pore at $r=0$. At $z=0.2\,L$ (Figure \ref{fig:cpx200}) we find that the cation concentration decreases by $0.018$. The decrease becomes more significant at the narrowest point where the concentration drops by $0.046$ (Figure \ref{fig:cpx500}). Close to the inlet on the other hand the additional advection term (due to the Stokes coupling) leads to an increase by $0.023$ (Figure \ref{fig:cpx800}).
For the anions, however, the influence of hydrodynamic coupling is minimal. The average deviation between the PNP and PNPS solutions is only $0.18$, which is negligible compared with the maximum concentration ($\approx 10$). Consequently, no substantial change in the anion distribution is observed across the pore cross-sections.
Comparing the asymptotic model (solid lines) with the 2D FEM solution (stars) shows that the asymptotic approximation also captures the overall dynamics of the ion concentrations in a trumpet shaped pore with surface charge and applied potential difference.
A comparison of the IV relations for the calssical PNP and PNPS cases ($a_\pm=0$) also shows that the liquid flow mainly influences the cation dynamics (Figure \ref{fig:IV_cp}). While the cation current for PNP is $I_+ = \pm 5$ at $\Delta \varphi = \pm 20$, we measure a current of only $I_+ = \pm 0.3$ in the PNPS case. For the anions, we measure a slightly higher current for the PNPS case than for the PNP model. For example, at $\Delta \varphi = \pm 20$, the current for PNPS is $I_+ = \pm 46.8$ and for PNP it is $I_+ = \pm 39.7$. Furthermore, Figure \ref{fig:IV_cn} shows that the anion current in the PNPS case appears to saturate for higher potentials, while the PNP current increases linearly.
In addition to the classical PNPS model, we have also calculated the IV relations for different volume fractions to illustrate the impact of finite-size effects. As can be observed in Figure \ref{fig:IV_cp_a} the cation current decreases for increasing $a_\pm$. Furthermore, it seems to saturate for $a_\pm = 5$ (dashed line) at around $\Delta \phi = \pm 20$. For the anion current, however, we do not notice a significant difference (Figure \ref{fig:IV_cn_a}).

\begin{figure}[H]
	\centering
    \begin{subfigure}{0.48\textwidth}
		\centering
		\includegraphics[width=1\textwidth]{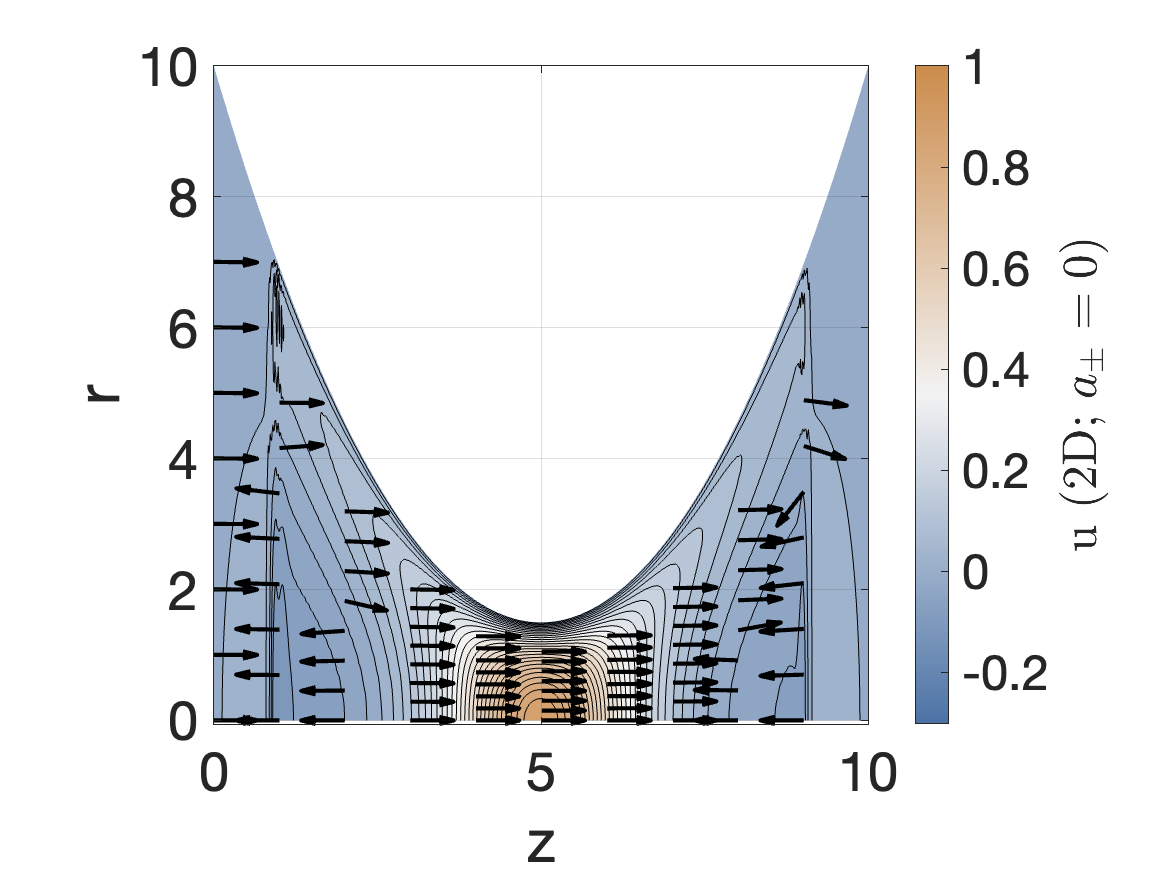}
		 \caption{}
		\label{fig:v_2D}
	\end{subfigure}
	\hfill
	\begin{subfigure}{0.48\textwidth}
		\centering
		\includegraphics[width=1\textwidth]{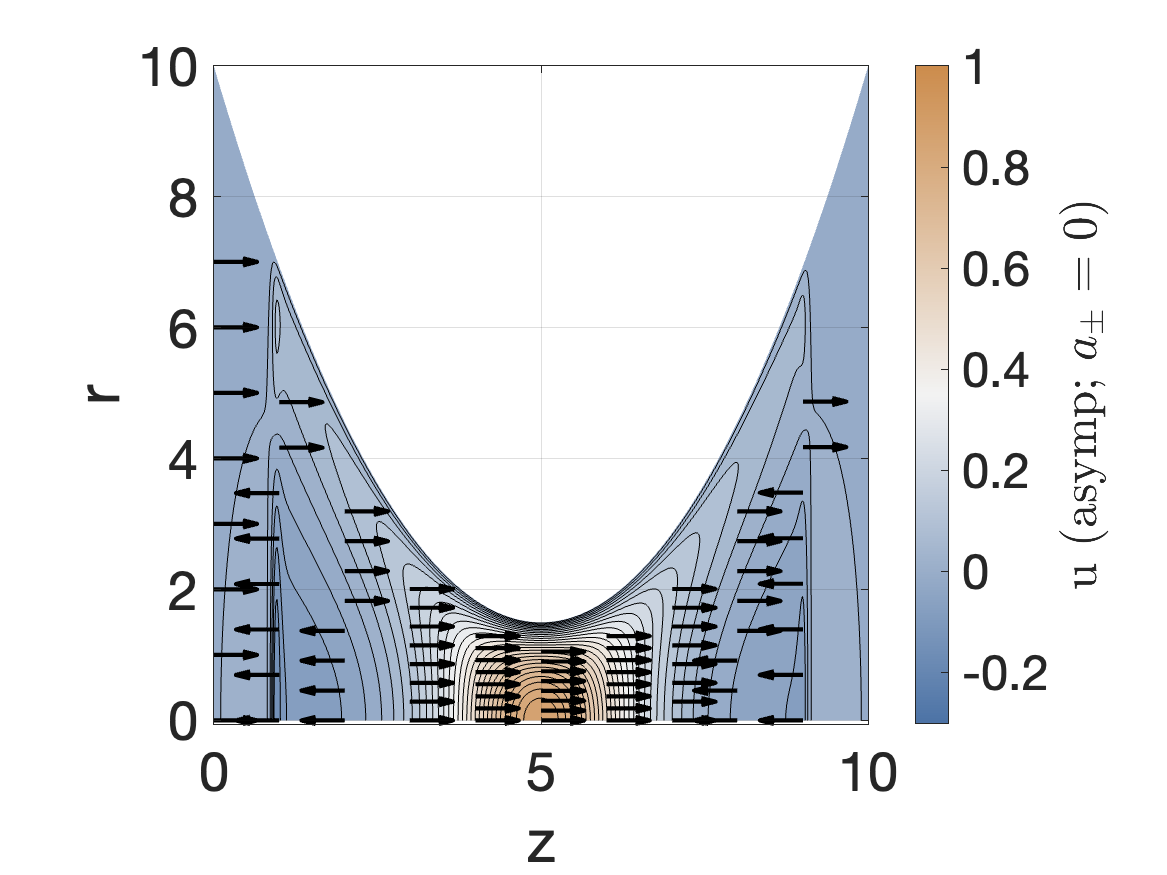}
		\caption{}
		\label{fig:v_asymp}
	\end{subfigure}
    \vfill
    \begin{subfigure}{0.48\textwidth}
		\centering
		\includegraphics[width=1\textwidth]{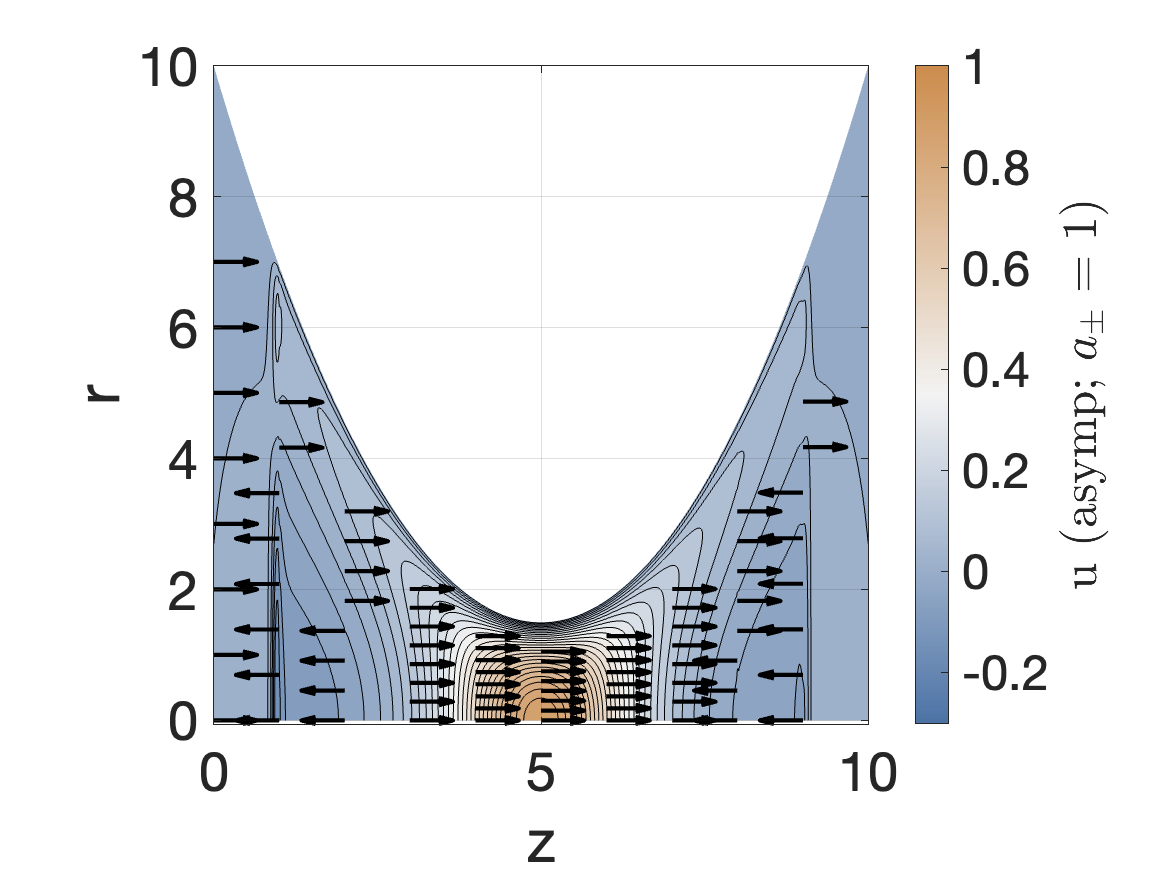}
		 \caption{}
		\label{fig:v_a1}
	\end{subfigure}
	\hfill
	\begin{subfigure}{0.48\textwidth}
		\centering
		\includegraphics[width=1\textwidth]{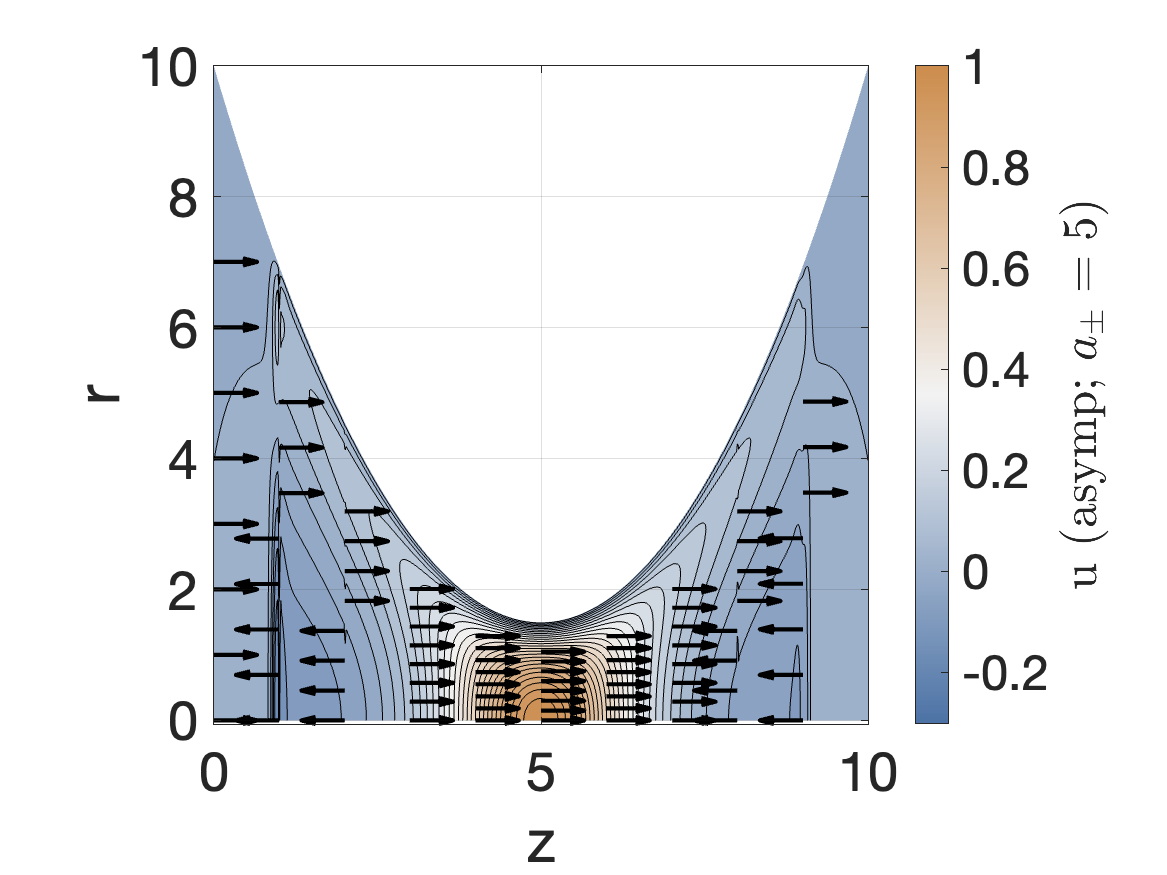}
		\caption{}
		\label{fig:v_a5}
	\end{subfigure}
	\caption{2D FEM FEniCS solution (a) and solution of the asymptotic model (b) of the velocity for a trumpet-shaped pore for $a_\pm = 0$, (c) for $a_\pm=1$ and (d) for $a_\pm=5$. Contourlines of the axial velocity $u$ plotted together with the vector field $\vecv =(w,u)$ as arrows. Note that the arrows are normalized.}
	\label{fig:VelocityColorbar_trumpet}
\end{figure}

In Figure \ref{fig:VelocityColorbar_trumpet} we plotted the contourlines of the axial velocity $u$ together with vector plot of $\vecv = (w,u)$ illustrated as arrows. 
It can be observed that the velocity near the inlet and outlet is approximately zero (blue regions) and increases toward the center of the pore, reaching a maximum of $u=0.88$ (orange region) at $z=0.5\,L$. Near the transition from no surface charge to positive surface charge, it can be seen that the velocity becomes negative (dark blue regions) and, as can be seen from the arrows, also changes direction. While most arrows point in the direction of the inlet, indicating a net transport in that direction, there are also arrows pointing in the opposite direction. 
Calculating the velocity field for different volume fractions $a_\pm$ illustrates that finite-size effects also some impact the velocity field. 
We find that with increasing volume fraction also the maximum velocity in the pore center increases. 
While for $a_\pm=1$ the maxiumum only slightly increases to $u=0.9$ (Figure \ref{fig:v_a1}), we find that for $a_\pm=5$ also the maximum velocity at $z=0.5\,L$ increases to $u=1$ (Figure \ref{fig:v_a5}).

Furthermore, we note that the asymptotic solution (Figure \ref{fig:v_asymp}) accurately reflects the overall dynamics of the velocity in the pores. It clearly illustrates the changes in the transport direction and the influence of the geometry, which leads to a maximum at the narrowest point. In summary, our asymptotic model accurately captures the influence of fluid flow and trumpet shape on ion flow through the pore.
\subsection{Conductivity in protein-based channels}\label{sec:Arbitrary}
\begin{figure}[H]
    \centering
    \begin{tikzpicture}[scale=1.2]

    \draw[->,line width=0.25mm] (-0.5,0) -- (9,0) node[right] {$z$}; 
    \draw[->,line width=0.25mm] (0,-0.5) -- (0,3.0) node[above] {$r$}; 

    \begin{axis}[
        xmin=0,
        ymin=0,
        smooth,
        scale only axis,
        axis lines=none,
        unit vector ratio=1 2 1,
    ]
        \addplot[
            line width=0.5mm
        ]
        table {img/willems_plot_data.txt};
    \end{axis}

    \draw[line width =0.5mm] (0,2.2) -- (0,0);
    \draw[line width =0.5mm] (7.655,0) -- (7.655,2.2);
    \draw[line width =0.5mm,dashed] (0,0) -- (7.65,0);
    
    \draw[line width =0.4mm,dashed] (2.2,0) -- (2.2,0.85) node[left] at (2.5,-0.3){ $2.0$};
    \draw[line width =0.4mm,dashed] (5.5,0) -- (5.5,1.7) node[left] at (5.85,-0.3){ $3.5$};
    \draw[line width =0.4mm,dashed] (0,0) -- (0,0) node[left] at (8.05,-0.3){ $5.5$};

    \draw[line width =0.4mm,dashed] (0,0.68) -- (2.4,0.68) node[left] at (0,0.68){ $1.65$};
     \draw[line width =0.4mm,dashed] (0,0) -- (0.0,0.0) node[left] at (0,2.2){ $5$};



    
\end{tikzpicture}
	\caption{Illustration of a 2D axially symmetric trumpet shaped pore $\Omega$. We apply a homogeneous surface charge at $S^\twall$, i.e. $r=R(z)$ for the potential. At the inlet $S^\tin$ at $\tilde z = 55$\,\si{\nano \meter} and the outlet $S^\tout$ at $\tilde z=0$\,\si{\nano \meter} we apply Dirichlet boundary conditions for the concentrations, the potential and the pressure. The shape of the pore $2 \leq z \leq 3.5$ was taken from the supplementary data from \cite{Willems2020}}
	\label{fig:ArbitraryPore}
\end{figure}
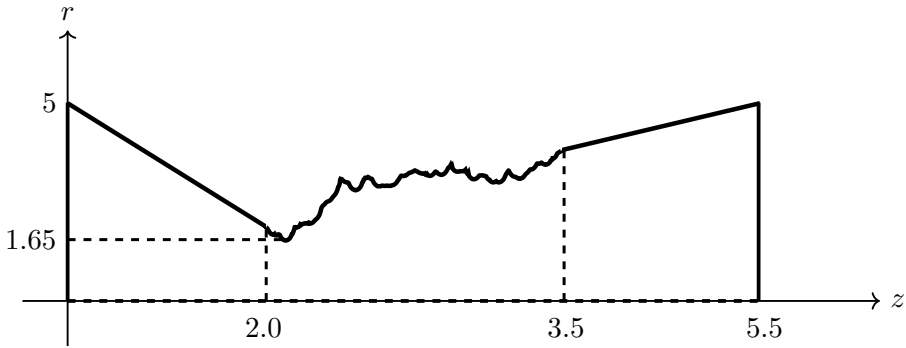

\begin{figure}[H]
	\centering
    \begin{subfigure}{0.44\textwidth}
		\centering
		\includegraphics[width=1\textwidth]{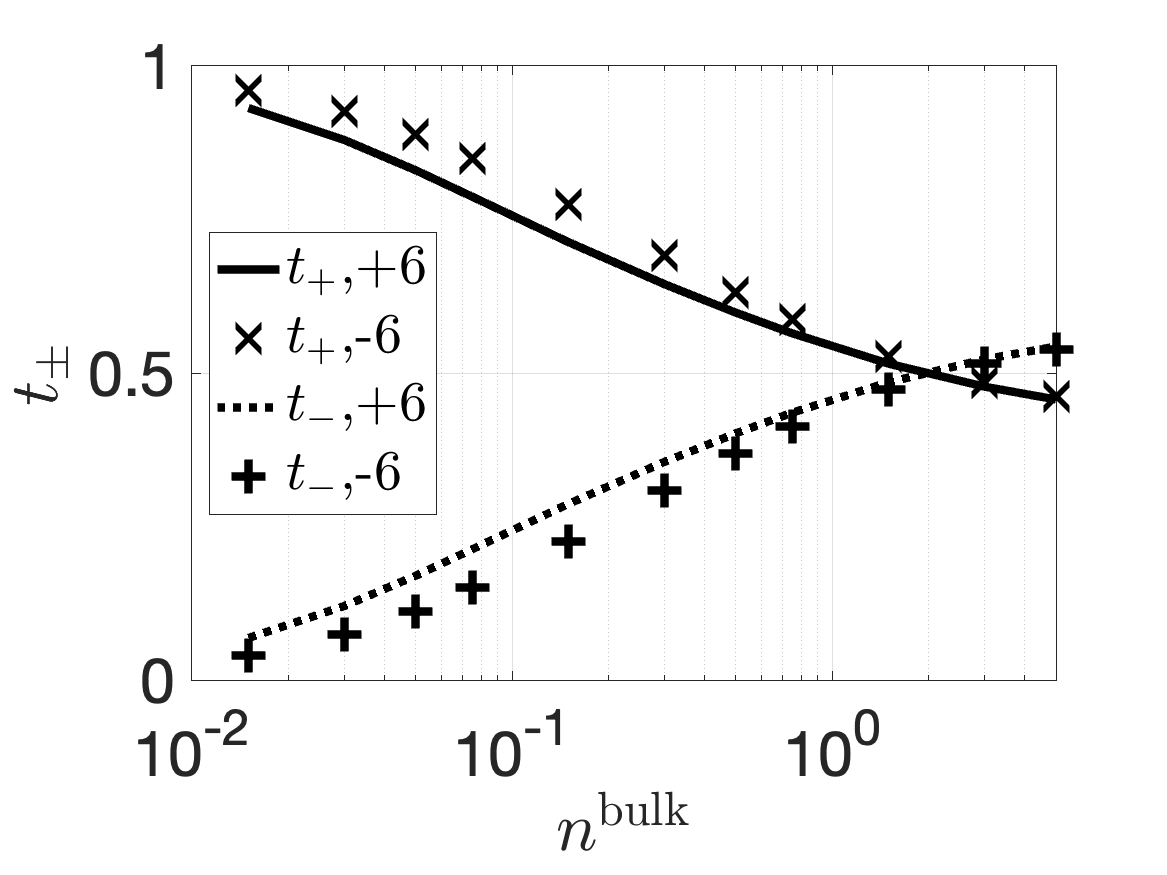}
		 \caption{}
		\label{fig:WillemsTransport}
	\end{subfigure}
    \hfill
    \begin{subfigure}{0.44\textwidth}
		\centering
		\includegraphics[width=1\textwidth]{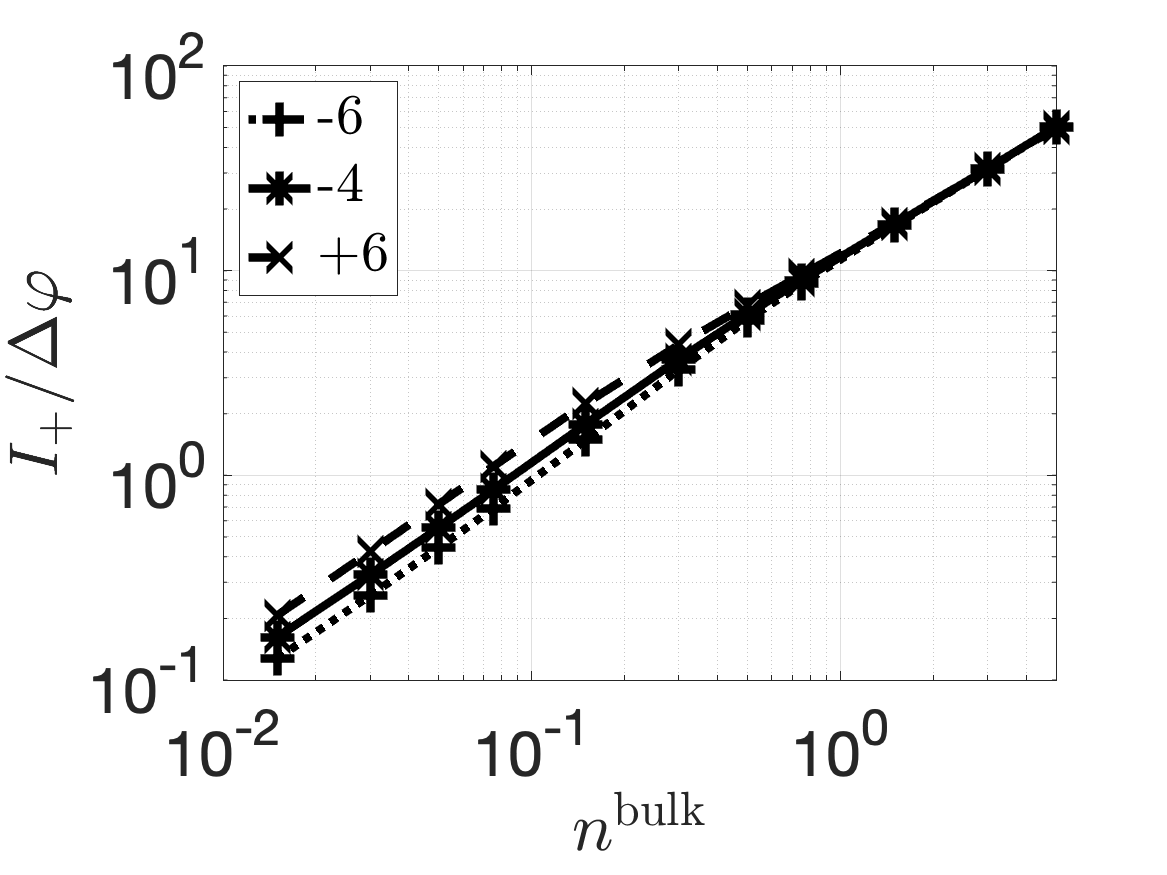}
		 \caption{}
		\label{fig:WillemsConductance}
	\end{subfigure}
	\caption{(a) Transport coefficients for the cations (solid line and crosses) and anions for (dotted line and plus signs) for $\Delta \phi = \pm 6$ for varying salt concentrations. (b) Cation conductance for $\Delta \phi = - 6$ (plus signs), $\Delta \phi = - 4$ (stars) and $\Delta \phi =  6$ (crosses). }
	\label{fig:Willems1}
\end{figure}

Cytolysin A (ClyA) is a pore forming toxin that is released by some strains of \textit{Escherichia coli} bacteria to damage host cells by creating pores in their membrane.
Since ClyA forms stable nanopores, it is also often used in experimental ion channel studies and is therefore an important application example for the study of mathematical models \cite{Willems2020}. 
The hydrophilic pore is around $15$\,\si{\nano \meter} long with a radius of around $1.65$\,\si{\nano \meter} at its narrowest point. More and more computational studies use geometry and charge density obtained by MD simulations to calculate IV-relations using a continuum description of the ion flux through the pore \cite{Liu2020,Xie2022,Willems2020}.
The geometry is often reconstructed using the contour lines of the charge density map of the proteins. However, solving the continuum model for the entire 3D domain is computationally intensive, and meshing the geometry is not trivial. To save computing costs, the MD geometry is therefore often reduced to an axially symmetric 2D domain. Since we have derived the asymptotic model for arbitrarily shaped pores with radius $R(z)$, we can use it to solve the model for more complex geometries derived from MD simulations without the need to define a mesh for it. To demonstrate this, we used the data of \cite{Willems2020} in the supplementary information and reconstructed the axially symmetric 2D shape of the ClyA pore. 
Our simulation domain is given in Figure \ref{fig:ArbitraryPore}. The pore is placed between two reservoirs of length $2$ each with $\delta = 0.1$, i.e. the entry of the pore is at $z=2$ and the outlet at $z=3.5$. The length of the overall domain is $L=5.5$. The two bath have a maximum radius of $R(z=0)=R(z=5.5)=5$ and narrow linearly toward the pore. 
We consider a mixture of K$^+ = n_+$ and Cl$^-= n_-$, with $n_\pm^\tout = n_\pm^\tin = n^\tbulk$ at $z=0$ and $z=5.5$. In addition, we apply a potential difference of $\Delta \phi =\varphi^\tin - \varphi^\tout$ with $\varphi^\tout$ at $z=0$ and $\varphi^\tin$ at $z=5.5$. We choose a solvent volume of $v_0 = 0.018$ ($v_0 = 18e-6$\,\si{\meter \cubed \per \mole}) and a ion volume of $v_\pm = 0.075$ ($v_\pm = 75e-6$\,\si{\meter \cubed \per \mole}), which gives $a_\pm = 4.15$ \cite{Willems2020}.
If not explicitly mentioned otherwise in the text the parameter values were chosen as given in Table \ref{tab:ParameterValuesGeneral} and Table \ref{tab:ParameterValuesArbitrary}. We apply a constant surface charge of $\sigma_0=-0.55$ given by equation \eqref{eq:homCharge} with $L_1=2$, $L_2=3.5$ and $\eps = 0.5$. 
\vspace{-4cm}
We calculated the transport coefficients for varying concentrations through the pore given by $ t_\pm = I_\pm/(I_+ + I_-)$. For our simplified model we find that it is capable to capture the overall behavior of ion flux through the pore. Figure \ref{fig:WillemsTransport} shows the cation transport coefficient $t_+$ as a function of concentration for two different potentials $\Delta \phi = \pm 6$ (solid line and crosses, respectively). 
While for small concentrations the ion conductance of the pore is mainly driven by the cation current, $t_+$ decreases with increasing concentration. Therefore, we find that the anion transport $t_-$ increases with increasing concentration (dotted line and plus signs) until the overall conductance is driven by the anion current. A switch in the selectivity of the pore can be observed at around $n^\tbulk=2$, similar as for the ePNP-NS model in Willems et al. \cite{Willems2020}. However, when comparing the data, it must be taken into account that, unlike Willems et al., we do not solve the full Navier-Stokes equations and also assume a homogeneous surface charge of the channel protein, a constant permittivity, constant diffusion coefficients, and a constant electrophoretic mobility. In addition, we are considering a smaller domain, which could also lead to a slightly different result. Nevertheless, we can capture the overall dynamics.
Since Willems et al. reported a preferred cation selectivity in their work we assume a homogeneous negative surface charge throughout the pore. However, we plan to take the coupling with charge density into account in future work.
As in section \ref{sec:Cylinder} we also calculated the cation conductance which is the current scaled with the applied potential difference $I_\pm / \Delta \phi$ and plotted it for varying salt concentrations (Figure \ref{fig:WillemsConductance}). Again we observe that for all applied potentials the curves fall into one master curve for high enough concentrations. For small concentrations we still observe a rightward (or leftward) shift of the curves.

\begin{figure}[H]
	\centering
    \begin{subfigure}{0.44\textwidth}
		\centering
		\includegraphics[width=1\textwidth]{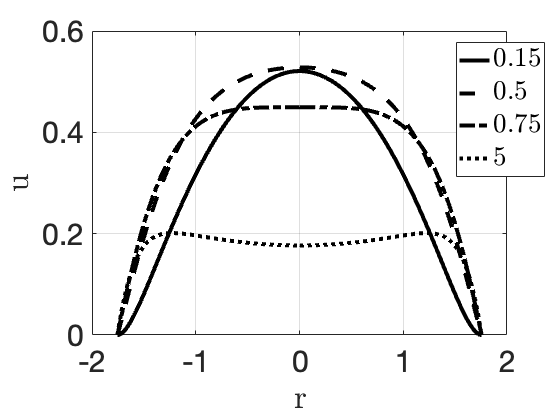}
		 \caption{}
		\label{fig:WillemsVelocity}
	\end{subfigure}
    \hfill
    \begin{subfigure}{0.44\textwidth}
		\centering
		\includegraphics[width=1\textwidth]{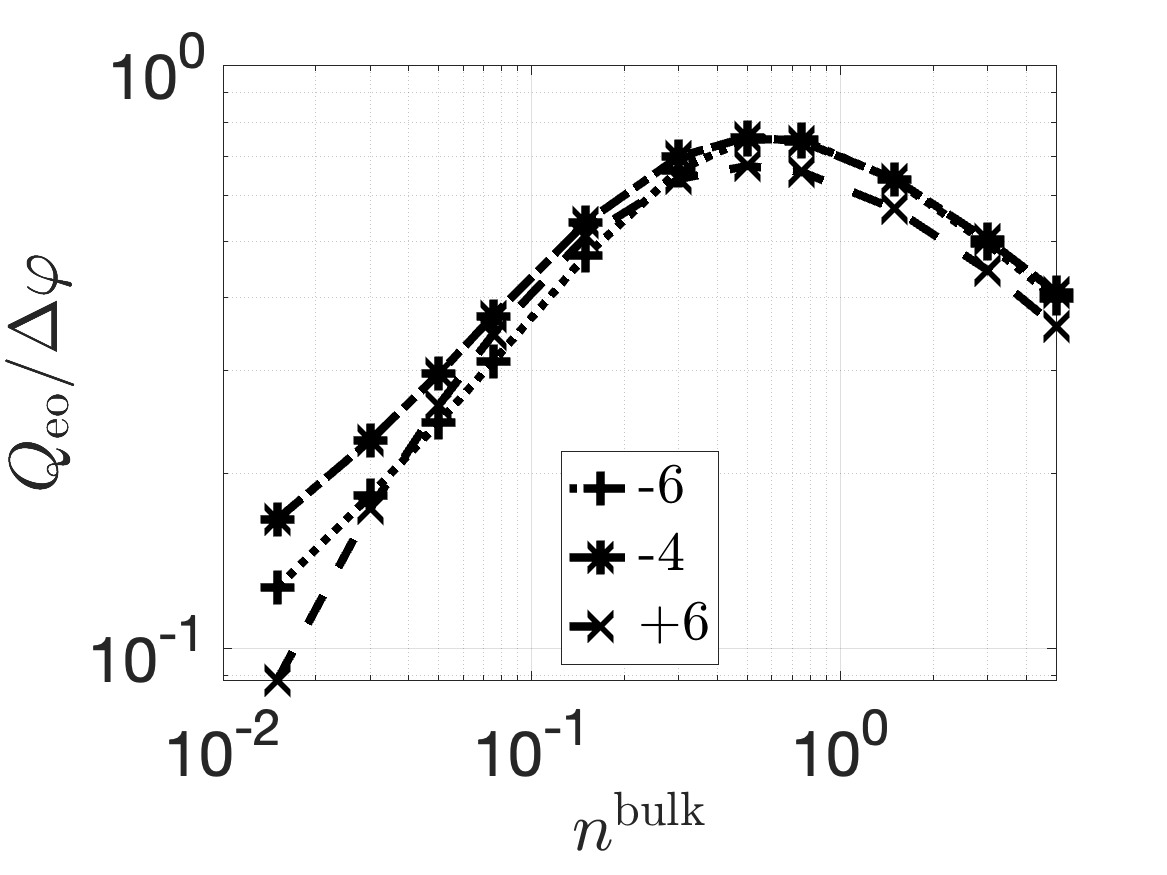}
		 \caption{}
		\label{fig:WillemsFlow}
	\end{subfigure}
	\caption{(a) Velocity calculated for different concentrations $n^\tbulk$ for $\Delta \phi = -4$. (b) Flow rate plotted as a function of $n^\tbulk$ for $\Delta  \phi = 6$ (crosses), $\Delta \phi = -4$ (stars) and $\Delta \phi = -6$ (plus signs).}
	\label{fig:Willems2}
\end{figure}

We also evaluated the velocity at $z=2.08$ for $\Delta \phi =-4$ and four different concentrations (Figure \ref{fig:WillemsVelocity}). Similar as reported by Willems et al. \cite{Willems2020} we find that the velocity in the center of the pore at $r=0$ at first slightly increases by increasing the concentrations from $n^\tbulk=0.15$ (solid) up to $n^\tbulk=0.5$. However, if we further increase the concentration to $n^\tbulk=0.75$ and $n^\tbulk=5$ we find that the velocity has decreased each time. 
Plotting the flow rate $Q_{eo} = \int_S \vecv \cdot \vecn \ud S = \int_0^{R(z)} u^0 r \ud r$ at $z=2.08$ as a function of concentration (Figure \ref{fig:WillemsFlow}) better illustrates this phenomenon. As noted in section \ref{sec:Cylinder}, scaling the flow with $\Delta \phi$ causes the different maxima for different potentials of the same sign. Up to a concentration of $n^\tbulk=0.5$, the flow increases, reaches its maximum of $Q_{eo}/\Delta \varphi = 0.75$ for $\Delta \phi =-4$ (stars) and $\Delta \phi =- 6$ (plus signs), and then decreases again for increasing concentrations. Similar for $\Delta \phi = +6$ where the flow reaches a maximum of $Q_{eo}/\Delta \varphi = 0.67$ (crosses) at $n^\tbulk=0.5$. Values reported by Willems et al. \cite{Willems2020} are again different, they report a maximum of $\tilde Q_{eo} /\Delta \tilde \varphi = 11.5$\,\si{\nano \meter \cubed \per \nano \second \per \volt} at $\tilde n^\tbulk=0.5\,\si{\mol \per \liter}$ while we find a maximum of $\tilde Q_{eo}  /\Delta \tilde \varphi = 7.5$\,\si{\nano \meter \cubed \per \nano \second \per \volt} also at $\tilde n^\tbulk=0.5$\,\si{\mol \per \liter}.

In summary it is possible to use pore-shapes obtained from MD simulations to capture the general dynamics in such channels with our asymptotic model. The observed deviations are due to the 2D model used. Since we only need to know the coordinates for the boundary in order to solve the equations for the asymptotic model, it is not necessary to define a mesh on such an irregular domain, thus saving an often very time-consuming step. 

Furthermore, we could show that scaling laws derived from simpler models, such as those for straight pores, can also be applied to more complex problems.
\section{Conclusion and Outlook}
The asymptotic model developed in this work provides a mathematically systematic and physically transparent framework for describing ion transport in narrow channels. Our approach complements and clarifies earlier one-dimensional reductions, such as those introduced by Fair and Osterle \cite{fair1971reverse} and later extended by Rankin et al. \cite{rankin2016effect}. 
These classical models are based on lubrication theory but are often not derived rigorously from the full PNPS system. Instead, they typically rely on simplifying assumptions such as a small Debye length and/or a uniform electric field, and they further require numerical evaluation of the leading-order axial pressure gradient.
In contrast, the present analysis identifies a distinct asymptotic region in which the Debye length is comparable to the channel width, leading to a closed-form analytical description for pressure and velocity that recovers previous reductions as limiting cases while significantly extending their range of validity.

Our studies show that the volumetric flow in a cylindrical pore with homogeneous surface charge smoothly transitions between an potential-dominated and a pressure-driven regime. For small Debye layers the volumetric flow can be expressed as a combination of two classical limits: in the potential-dominated case, it reduces to the Helmholtz-Smoluchowski approximation, while in the limit of strong hydrodynamic pressure gradient, it approaches the Poiseuille profile.
A similar transition occurs in ion currents. Depending on the applied voltage and the pressure gradient, the total current shifts between an potential-controlled and a pressure-controlled regime, reflecting the coupled electrohydrodynamic response of the system. Furthermore, our calculations show that a correction is necessary for large Debye layers. The EDL axial force density therefore provides an additional flow contribution and is necessary, for example, for low salt concentrations when the surface charge cannot be sufficiently screened.

In addition to electrostatic and hydrodynamic forces, there are other properties of nanopores that influence ion transport through the pore. In particular, geometric constraints and deviations from ideal cylindrical symmetry can have a strong influence on ion flow and IV properties. Fluctuations in pore radius can lead to spatially heterogeneous electric fields and uneven local conductivities. These geometric features can influence ion selectivity and surface charge regulation and promote rectification, thereby affecting overall transport behavior. 
Our results indicate that such effects are particularly pronounced when the characteristic pore dimension is of the order of the Debye length and additional excluded‑volume effects occur. 

By incorporating geometric complexities and finite-size effects into the asymptotic framework, we obtain a more complete picture of how pore shape, surface chemistry, space‑charge competition, and flow collectively control ion transport in realistic nanofluidic systems.
Our generalized approach offers the ability to include various chemical potential functions, such as those that take into account steric effects, i.e., ion-ion interactions, solvation effects, ions of different sizes, etc.
The framework thus covers several models commonly used in the literature, such as the classical PNPS model \cite{Green2022,Li2024}, the Bikerman-Freise model \cite{Bikerman1942,Borukhov1997,Liu2013,Liu2020,Willems2020}, and the incompressible idealized mixture model \cite{Dreyer2013,Keller2025,Landstorfer2016}. Making it easier to compare different approaches.

Furthermore, our asymptotic framework accounts for axial variations in pore radius and, as illustrated by trumpet-shaped and protein-based pores obtained from MD simulations, accurately captures their associated geometric effects.
In our investigation of a trumpet-shaped pore, we observe that ion concentrations vary systematically from the wider areas toward the bottleneck, and this axial variation is further amplified by the coupling between ion transport and local fluid velocity. 
An evaluation of the velocity contourplots also shows that the inclusion of the ion volumes has a significant influence, especially within the confined region. 
The asymptotic model shows that the varying geometry induces an additional pressure gradient that depends on the ion concentrations. 

This becomes particularly clear when examining the fluid velocity and ion conductivity as a function of ion concentration in a ClyA-based channel obtained from MD simulations. Both variables increase with rising salt concentration, but begin to decrease once a certain threshold value is reached. 
Taken together, these observations illustrate how geometry and flow jointly influence the nonlinear transport response of the system.

The flexibility of the asymptotic framework opens up several possibilities for further development. A natural extension concerns the inclusion of hydrodynamic slip at the channel walls. The slip can be included in the reduced Stokes problem as a first approximation, where its effects on the flow field can be profound. Even without electrostatic forces, different slip length scales can fundamentally alter flow resistance, induce transitions between plug-like and shear-dominated profiles, and modify pressure-flow relationships \cite{munch2005lubrication,peschka2019signatures}. Previous studies \cite{datta2013effect,park2009extension} have highlighted the complexity of such slip-induced transitions, and a systematic investigation of their interaction with electroviscous and electroosmotic effects represents an important direction for future work.

By providing analytical insights and computational efficiency across a wide range of parameters, the framework presented here provides a solid foundation for investigating electrohydrodynamic phenomena in biological ion channels, synthetic nanopores, and new nanofluidic technologies.
\section*{Acknowledgements}
C.K. gratefully acknowledges the financial support by the Deutsche Forschungsgemeinschaft (DFG, German Research Foundation) under Germany's Excellence Strategy – The Berlin Mathematics Research Center MATH+ (EXC-2046/1, project ID: 390685689) via the project AA1-14.

The authors would like to thank Dirk Peschka for his support with FEniCS and Python.

AI assistance was used at an early stage of the 2D FEniCS implementation and in parts of the numerical code for the asymptotic model, primarily for technical guidance (e.g., on data interpolation) and initial development. AI tools were also used to suggest color schemes for two contour plots (subsequently adjusted manually) and to help improve phrasing in parts of the manuscript.

\section*{Data availability}
The data and code will be made available upon reasonable request.

\section*{Competing interests}
The authors declare no competing interests.

\bibliographystyle{abbrv}
\bibliography{IonChannel}

@Article{Landstorfer2016,
  author    = {M. Landstorfer and C. Guhlke and W. Dreyer},
  journal   = {Electrochimica Acta},
  title     = {Theory and structure of the metal-electrolyte interface incorporating adsorption and solvation effects},
  year      = {2016},
  month     = {may},
  pages     = {187--219},
  volume    = {201},
  doi       = {10.1016/j.electacta.2016.03.013},
  publisher = {Elsevier {BV}},
}

@Article{Liu2020,
  author    = {Jinn-Liang Liu and Bob Eisenberg},
  journal   = {Entropy},
  title     = {Molecular Mean-Field Theory of Ionic Solutions: A Poisson-Nernst-Planck-Bikerman Model},
  year      = {2020},
  month     = {may},
  number    = {5},
  pages     = {550},
  volume    = {22},
  doi       = {10.3390/e22050550},
  publisher = {{MDPI} {AG}},
}

@Article{Liu2013,
  author    = {Jinn-Liang Liu and Bob Eisenberg},
  journal   = {The Journal of Physical Chemistry B},
  title     = {Correlated Ions in a Calcium Channel Model: A Poisson{\textendash}Fermi Theory},
  year      = {2013},
  month     = {sep},
  number    = {40},
  pages     = {12051--12058},
  volume    = {117},
  doi       = {10.1021/jp408330f},
  publisher = {American Chemical Society ({ACS})},
}

@Article{Dreyer2013,
  author    = {Wolfgang Dreyer and Clemens Guhlke and Rüdiger Müller},
  journal   = {Physical Chemistry Chemical Physics},
  title     = {Overcoming the shortcomings of the Nernst{\textendash}Planck model},
  year      = {2013},
  number    = {19},
  pages     = {7075},
  volume    = {15},
  doi       = {10.1039/c3cp44390f},
  publisher = {Royal Society of Chemistry ({RSC})},
}

@Article{Gillespie2008,
  author          = {Gillespie, Dirk},
  journal         = {Biophysical journal},
  title           = {Energetics of divalent selectivity in a calcium channel: the ryanodine receptor case study.},
  year            = {2008},
  issn            = {1542-0086},
  month           = feb,
  pages           = {1169--1184},
  volume          = {94},
  chemicals       = {Ryanodine Receptor Calcium Release Channel, Calcium},
  citation-subset = {IM},
  completed       = {2008-02-11},
  country         = {United States},
  doi             = {10.1529/biophysj.107.116798},
  issn-linking    = {0006-3495},
  issue           = {4},
  nlm-id          = {0370626},
  pii             = {S0006-3495(08)70635-0},
  pmc             = {PMC2212702},
  pmid            = {17951303},
  pubmodel        = {Print-Electronic},
  pubstate        = {ppublish},
  revised         = {2021-10-20},
}

@Article{Landstorfer2018,
  author    = {M. Landstorfer},
  journal   = {Electrochemistry Communications},
  title     = {On the dissociation degree of ionic solutions considering solvation effects},
  year      = {2018},
  month     = {jul},
  pages     = {56--59},
  volume    = {92},
  doi       = {10.1016/j.elecom.2018.05.011},
  publisher = {Elsevier {BV}},
}

@Book{de2013non,
  author    = {De Groot, S.R. and Mazur, P.},
  publisher = {Dover Publications},
  title     = {Non-Equilibrium Thermodynamics},
  year      = {2013},
  isbn      = {9780486153506},
  series    = {Dover Books on Physics},
  url       = {https://books.google.de/books?id=mfFyG9jfaMYC},
}

@Article{Bikerman1942,
  author    = {J.J. Bikerman},
  journal   = {The London, Edinburgh, and Dublin Philosophical Magazine and Journal of Science},
  title     = {{XXXIX}. Structure and capacity of electrical double layer},
  year      = {1942},
  month     = {may},
  number    = {220},
  pages     = {384--397},
  volume    = {33},
  doi       = {10.1080/14786444208520813},
  publisher = {Informa {UK} Limited},
}

@Article{Willems2020,
  author    = {Willems, Kherim and Ruić, Dino and L. R. Lucas, Florian and Barman, Ujjal and Verellen, Niels and Hofkens, Johan and Maglia, Giovanni and Van Dorpe, Pol},
  journal   = {Nanoscale},
  title     = {Accurate modeling of a biological nanopore with an extended continuum framework},
  year      = {2020},
  issn      = {2040-3372},
  number    = {32},
  pages     = {16775--16795},
  volume    = {12},
  doi       = {10.1039/d0nr03114c},
  publisher = {Royal Society of Chemistry (RSC)},
}

@Article{matejczyk2018,
  author    = {Matejczyk, B. and Pietschmann, J.-F. and Wolfram, M.-T. and Richardson, G.},
  title     = {Asymptotic models for transport in large aspect ratio nanopores},
  journal   = {European Journal of Applied Mathematics},
  year      = {2018},
  volume    = {30},
  number    = {3},
  pages     = {557--584},
  month     = jun,
  issn      = {1469-4425},
  doi       = {10.1017/s0956792518000293},
  publisher = {Cambridge University Press (CUP)},
}

@Article{Green2021,
  author    = {Green, Yoav},
  journal   = {The Journal of Chemical Physics},
  title     = {Ion transport in nanopores with highly overlapping electric double layers},
  year      = {2021},
  issn      = {1089-7690},
  month     = feb,
  number    = {8},
  volume    = {154},
  doi       = {10.1063/5.0037873},
  publisher = {AIP Publishing},
}

@Article{Green2022,
  author    = {Green, Yoav},
  journal   = {Physical Review Fluids},
  title     = {Effects of surface-charge regulation, convection, and slip lengths on the electrical conductance of charged nanopores},
  year      = {2022},
  issn      = {2469-990X},
  month     = jan,
  number    = {1},
  pages     = {013702},
  volume    = {7},
  doi       = {10.1103/physrevfluids.7.013702},
  publisher = {American Physical Society (APS)},
}

@Article{Li2024,
  author    = {Li, Minglun and Muthukumar, Murugappan},
  journal   = {The Journal of Chemical Physics},
  title     = {Electro-osmotic flow in nanoconfinement: Solid-state and protein nanopores},
  year      = {2024},
  issn      = {1089-7690},
  month     = feb,
  number    = {8},
  volume    = {160},
  doi       = {10.1063/5.0185574},
  publisher = {AIP Publishing},
}

@Article{Keller2025,
  author    = {Keller, Christine and Landstorfer, Manuel and Fuhrmann, Jürgen and Wagner, Barbara},
  journal   = {Entropy},
  title     = {A Model Framework for Ion Channels with Selectivity Filters Based on Non-Equilibrium Thermodynamics},
  year      = {2025},
  issn      = {1099-4300},
  month     = sep,
  number    = {9},
  pages     = {981},
  volume    = {27},
  doi       = {10.3390/e27090981},
  publisher = {MDPI AG},
}

@article{siwy2004conical,
  title={Conical-nanotube ion-current rectifiers: the role of surface charge},
  author={Siwy, Zuzanna and Heins, Elizabeth and Harrell, C Chad and Kohli, Punit and Martin, Charles R},
  journal={Journal of the American Chemical Society},
  volume={126},
  number={35},
  pages={10850--10851},
  year={2004},
  publisher={ACS Publications}
}

@article{harrell2004dna,
  title={DNA- nanotube artificial ion channels},
  author={Harrell, C Chad and Kohli, Punit and Siwy, Zuzanna and Martin, Charles R},
  journal={Journal of the American Chemical Society},
  volume={126},
  number={48},
  pages={15646--15647},
  year={2004},
  publisher={ACS Publications}
}

@article{dal2019confinement,
  title={Confinement-controlled rectification in a geometric nanofluidic diode},
  author={Dal Cengio, Sara and Pagonabarraga, Ignacio},
  journal={The Journal of Chemical Physics},
  volume={151},
  number={4},
  year={2019},
  publisher={AIP Publishing}
}

@article{stroock2000patterning,
  title={Patterning electro-osmotic flow with patterned surface charge},
  author={Stroock, Abraham D and Weck, Marcus and Chiu, Daniel T and Huck, Wilhelm TS and Kenis, Paul JA and Ismagilov, Rustem F and Whitesides, George M},
  journal={Physical review letters},
  volume={84},
  number={15},
  pages={3314},
  year={2000},
  publisher={APS}
}

@article{jubin2018dramatic,
  title={Dramatic pressure-sensitive ion conduction in conical nanopores},
  author={Jubin, Laetitia and Poggioli, Anthony and Siria, Alessandro and Bocquet, Lyd{\'e}ric},
  journal={Proceedings of the National Academy of Sciences},
  volume={115},
  number={16},
  pages={4063--4068},
  year={2018},
  publisher={National Academy of Sciences}
}

@article{branton2008potential,
  title={The potential and challenges of nanopore sequencing},
  author={Branton, Daniel and Deamer, David W and Marziali, Andre and Bayley, Hagan and Benner, Steven A and Butler, Thomas and Di Ventra, Massimiliano and Garaj, Slaven and Hibbs, Andrew and Huang, Xiaohua and others},
  journal={Nature biotechnology},
  volume={26},
  number={10},
  pages={1146--1153},
  year={2008},
  publisher={Nature Publishing Group US New York}
}

@article{datta2013effect,
  title={Effect of hydrodynamic slippage on electro-osmotic flow in zeta potential patterned nanochannels},
  author={Datta, S and Choudhary, JN},
  journal={Fluid Dynamics Research},
  volume={45},
  number={5},
  pages={055502},
  year={2013},
  publisher={IOP Publishing}
}

@article{bocquet2010nanofluidics,
  title={Nanofluidics, from bulk to interfaces},
  author={Bocquet, Lyd{\'e}ric and Charlaix, Elisabeth},
  journal={Chemical Society Reviews},
  volume={39},
  number={3},
  pages={1073--1095},
  year={2010},
  publisher={Royal Society of Chemistry}
}

@article{peschka2019signatures,
  title={Signatures of slip in dewetting polymer films},
  author={Peschka, Dirk and Haefner, Sabrina and Marquant, Ludovic and Jacobs, Karin and M{\"u}nch, Andreas and Wagner, Barbara},
  journal={Proceedings of the National Academy of Sciences},
  volume={116},
  number={19},
  pages={9275--9284},
  year={2019},
  publisher={National Academy of Sciences}
}

@article{fair1971reverse,
  title={Reverse electrodialysis in charged capillary membranes},
  author={Fair, JC and Osterle, JF},
  journal={The Journal of Chemical Physics},
  volume={54},
  number={8},
  pages={3307--3316},
  year={1971},
  publisher={American Institute of Physics}
}

@article{rankin2016effect,
  title={The effect of hydrodynamic slip on membrane-based salinity-gradient-driven energy harvesting},
  author={Rankin, Daniel Justin and Huang, David Mark},
  journal={Langmuir},
  volume={32},
  number={14},
  pages={3420--3432},
  year={2016},
  publisher={ACS Publications}
}

@article{park2009extension,
  title={Extension of the Helmholtz-Smoluchowski velocity to the hydrophobic microchannels with velocity slip},
  author={Park, HM and Kim, TW},
  journal={Lab on a Chip},
  volume={9},
  number={2},
  pages={291--296},
  year={2009},
  publisher={Royal Society of Chemistry}
}

@article{schoch2008transport,
  title={Transport phenomena in nanofluidics},
  author={Schoch, Reto B and Han, Jongyoon and Renaud, Philippe},
  journal={Reviews of modern physics},
  volume={80},
  number={3},
  pages={839--883},
  year={2008},
  publisher={APS}
}

@Article{Xie2022,
  author    = {Xie, Dexuan},
  journal   = {Journal of Computational Physics},
  title     = {An efficient finite element iterative method for solving a nonuniform size modified Poisson-Boltzmann ion channel model},
  year      = {2022},
  issn      = {0021-9991},
  month     = dec,
  pages     = {111556},
  volume    = {470},
  doi       = {10.1016/j.jcp.2022.111556},
  publisher = {Elsevier BV},
}

@Article{Curk2024,
  author    = {Curk, Tine and Leyva, Sergi G. and Pagonabarraga, Ignacio},
  journal   = {Physical Review Letters},
  title     = {Discontinuous Transition in Electrolyte Flow through Charge-Patterned Nanochannels},
  year      = {2024},
  issn      = {1079-7114},
  month     = aug,
  number    = {7},
  pages     = {078201},
  volume    = {133},
  doi       = {10.1103/physrevlett.133.078201},
  publisher = {American Physical Society (APS)},
}

@Article{Wang2021,
  author    = {Wang, Ceming and Sensale, Sebastian and Pan, Zehao and Senapati, Satyajyoti and Chang, Hsueh-Chia},
  journal   = {Nature Communications},
  title     = {Slowing down DNA translocation through solid-state nanopores by edge-field leakage},
  year      = {2021},
  issn      = {2041-1723},
  month     = jan,
  number    = {1},
  volume    = {12},
  doi       = {10.1038/s41467-020-20409-4},
  publisher = {Springer Science and Business Media LLC},
}

@Article{Saharia2021,
  author    = {Saharia, Jugal and Bandara, Y. M. Nuwan D. Y. and Karawdeniya, Buddini I. and Hammond, Cassandra and Alexandrakis, George and Kim, Min Jun},
  journal   = {RSC Advances},
  title     = {Modulation of electrophoresis, electroosmosis and diffusion for electrical transport of proteins through a solid-state nanopore},
  year      = {2021},
  issn      = {2046-2069},
  number    = {39},
  pages     = {24398--24409},
  volume    = {11},
  doi       = {10.1039/d1ra03903b},
  publisher = {Royal Society of Chemistry (RSC)},
}

@Article{Koehler2018,
  author    = {Köhler, Mateus Henrique and Bordin, José Rafael and Barbosa, Marcia C.},
  journal   = {The Journal of Chemical Physics},
  title     = {2D nanoporous membrane for cation removal from water: Effects of ionic valence, membrane hydrophobicity, and pore size},
  year      = {2018},
  issn      = {1089-7690},
  month     = feb,
  number    = {22},
  volume    = {148},
  doi       = {10.1063/1.5013926},
  publisher = {AIP Publishing},
}

@Article{Surwade2015,
  author    = {Surwade, Sumedh P. and Smirnov, Sergei N. and Vlassiouk, Ivan V. and Unocic, Raymond R. and Veith, Gabriel M. and Dai, Sheng and Mahurin, Shannon M.},
  journal   = {Nature Nanotechnology},
  title     = {Water desalination using nanoporous single-layer graphene},
  year      = {2015},
  issn      = {1748-3395},
  month     = mar,
  number    = {5},
  pages     = {459--464},
  volume    = {10},
  doi       = {10.1038/nnano.2015.37},
  publisher = {Springer Science and Business Media LLC},
}

@Article{CohenTanugi2012,
  author    = {Cohen-Tanugi, David and Grossman, Jeffrey C.},
  journal   = {Nano Letters},
  title     = {Water Desalination across Nanoporous Graphene},
  year      = {2012},
  issn      = {1530-6992},
  month     = jun,
  number    = {7},
  pages     = {3602--3608},
  volume    = {12},
  doi       = {10.1021/nl3012853},
  publisher = {American Chemical Society (ACS)},
}

@Article{Borukhov1997,
  author    = {Borukhov, Itamar and Andelman, David and Orland, Henri},
  journal   = {Physical Review Letters},
  title     = {Steric Effects in Electrolytes: A Modified Poisson-Boltzmann Equation},
  year      = {1997},
  issn      = {1079-7114},
  month     = jul,
  number    = {3},
  pages     = {435--438},
  volume    = {79},
  doi       = {10.1103/physrevlett.79.435},
  publisher = {American Physical Society (APS)},
}

@article{munch2005lubrication,
  title={Lubrication models with small to large slip lengths},
  author={M{\"u}nch, Andreas and Wagner, BA and Witelski, Thomas P},
  journal={Journal of Engineering Mathematics},
  volume={53},
  number={3},
  pages={359--383},
  year={2005},
  publisher={Springer}
}

@article{craster2009dynamics,
  title={Dynamics and stability of thin liquid films},
  author={Craster, Richard V and Matar, Omar K},
  journal={Reviews of modern physics},
  volume={81},
  number={3},
  pages={1131--1198},
  year={2009},
  publisher={APS}
}

@book{whitham2011linear,
  title={Linear and nonlinear waves},
  author={Whitham, Gerald Beresford},
  year={2011},
  publisher={John Wiley \& Sons}
}

@Article{Nonner1998,
  author    = {Nonner, Wolfgang and Eisenberg, Bob},
  journal   = {Biophysical Journal},
  title     = {Ion Permeation and Glutamate Residues Linked by Poisson-Nernst-Planck Theory in L-Type Calcium Channels},
  year      = {1998},
  issn      = {0006-3495},
  month     = sep,
  number    = {3},
  pages     = {1287--1305},
  volume    = {75},
  doi       = {10.1016/s0006-3495(98)74048-2},
  publisher = {Elsevier BV},
}

@PhdThesis{Gillespie1999,
  author = {Dirk Gillespie},
  school = {Rush University, Chicago},
  title  = {A Singular Perturbation Analysis of thePoisson-Nernst-Planck System: Applications to Ionic Channels},
  year   = {1999},
}

@PhdThesis{Matejczyk2019,
  author = {Bartlomiej Matejczyk},
  school = {The University of Warwick},
  title  = {Mathematical modelling and simulations of the iontransport through confined geometries.},
  year   = {2019},
}
\newpage
\appendix
\section*{Appendix}
\section{Auxiliary calculations for the large aspect ratio asymptotics}
\paragraph{Poisson equation.} In order to solve equation \eqref{eq:CalculationPoisson} 
\begin{align}
    0 &=\int_{0}^{R(z)} \partial_r (r \partial_r \phi^1) \ud r + 	\int_{0}^{R(z)} \partial_z (\partial_z \phi^0 ) r \ud r,
\end{align}  
we went through the following steps.
Switching integration and differentiation in the second integral gives
\begin{align}
    \int_{0}^{R(z)} \partial_z (\partial_z \phi^0 ) r \ud r = - R'(z) \partial_z \phi^0 |_{r=R(z)} + \partial_z \int_{0}^{R(z)} \partial_z \phi^0 r \ud r.
\end{align}
Since the $\mathcal{O}(\delta^2)$-order term of the surface charge boundary condition \eqref{eqn:scaledbc:c} yields
\begin{align}
    \partial_r \phi^1|_{r=R(z)} - R'(z) \partial_z \phi^0 |_{r=R(z)} = 0,
\end{align}
we obtain
\begin{align}\label{eq:AppendixPoisson}
    0 &=\int_{0}^{R(z)} \partial_r (r \partial_r \phi^1) \ud r + 	\int_{0}^{R(z)} \partial_z (\partial_z \phi^0 ) r \ud r  \notag \\ \nonumber
    &= R(z) \underbrace{\left(\partial_r \phi^1|_{r=R(z)} - R'(z) \partial_z \phi^0 |_{r=R(z)} \right)}_{=0}  + \partial_z \int_{0}^{R(z)} \partial_z \phi^0 r \ud r \\ \nonumber
    &= \partial_z \int_{0}^{R(z)} \partial_z \phi^0_z r \ud r.
\end{align} 
\paragraph{Nernst-Planck equation.}\label{sec:AppendixCalculationConcentration}
In order to solve equation \eqref{eq:CalculationConcentration} we introduced the stream function 
\begin{align}
    s = \int_{0}^{r} u^0 r \ud r,
\end{align}
with $u^0 = r^{-1} \partial_r s$ and $w = - r^{-1} \partial_z s$, such that
\begin{align}\label{eq:AppendixConcentrations_1}
     \partial_t n^0_\alpha &= -P_e \int_{0}^{R(z)} \left(w \partial_r n_\alpha^0 + u^0 \partial_z n_\alpha^0 \right) \ud r + \int_{0}^{R(z)} \partial_z f_\alpha^0 r \ud r + \int_{0}^{R(z)} \partial_r \left(r g_\alpha^1\right) \ud r \\ \nonumber
     &= -P_e \int_{0}^{R(z)} \left(-\partial_z s \partial_r n_\alpha^0 + \partial_r s \partial_z n_\alpha^0 \right) \ud r + \int_{0}^{R(z)} \partial_z f_\alpha^0 r \ud r + \int_{0}^{R(z)} \partial_r \left(r g_\alpha^1\right) \ud r.
\end{align}
We rewrite the first integral such that
\begin{align}
    \int_{0}^{R(z)} & \left(-\partial_z s \partial_r n_\alpha^0 + \partial_r s \partial_z n_\alpha^0 \right) \ud r  \\ \nonumber
    &= \int_{0}^{R(z)} \left[-\partial_z(n_\alpha^0 \partial_r s) + n_\alpha^0 \partial_z(\partial_r s) + \partial_r(n_\alpha^0 \partial_z s) - n_\alpha^0 \partial_z(\partial_r s)\right] \ud r \\ \nonumber
    & = -\int_{0}^{R(z)}\partial_z(n_\alpha^0 \partial_r s) \ud r + \int_{0}^{R(z)}\partial_r(n_\alpha^0 \partial_z s) \ud r.
\end{align}
Changing integration and differentiation in the first integral and solving the second integral further yields
\begin{align}
    -\int_{0}^{R(z)}& \partial_z(n_\alpha^0 \partial_r s) \ud r + \int_{0}^{R(z)}\partial_r(n_\alpha^0 \partial_z s) \ud r \\ \nonumber
    &= - \partial_z \int_{0}^{R(z)} n_\alpha^0 \partial_r s \ud r +  n_\alpha^0 \partial_r s |_{r=R(z)} R'(z) + n_\alpha^0 \partial_z s |_{r=R(z)} - n_\alpha^0 \partial_z s |_{r=0} \\ \nonumber
    &= -\partial_z \int_{0}^{R(z)} n_\alpha^0 \partial_r s \ud r.
\end{align}
Due to the no-slip condition $u^0|_{r=R(z)}=w^0|_{r=R(z)}=0$ also $\partial_r s|_{r=R(z)} = 0$ and $\partial_z s|_{r=R(z)} = 0$ and with the symmetry condition $\partial_z s|_{r=0} = 0$. Using this result we have
\begin{align}\label{eq:AppendixConcentrations_2}
     \partial_t n_\alpha^0 &= -P_e \int_{0}^{R(z)} \left(-\partial_z s \partial_r n_\alpha^0 + \partial_r s \partial_z n_\alpha^0 \right) \ud r + \int_{0}^{R(z)} \partial_z f_\alpha^0 r \ud r + \int_{0}^{R(z)} \partial_r \left(r g_\alpha^1\right) \ud r \\ \nonumber
    &= P_e \partial_z \int_{0}^{R(z)} n_\alpha^0 \partial_r s  \ud r + \partial_z \int_{0}^{R(z)} f_\alpha^0 r \ud r - R(z) \underbrace{\left( R'(z) f_\alpha^0 |_{r=R(z)}- g_\alpha^1 |_{r=R(z)} \right)}_{=0} \\ \nonumber
    &= P_e \partial_z \int_{0}^{R(z)} n_\alpha^0 u^0 r \ud r + k_\alpha \partial_z \int_{0}^{R(z)} \left( \partial_z n_\alpha^0 + z_\alpha n_\alpha^0 \partial_z \phi^0 \right) r \ud r,
\end{align}
with the boundary condition \eqref{eqn:scaledbc:d}
\begin{align}
	R'(z)  f_\alpha^0 |_{r=R(z)} - g_\alpha^1 |_{r=R(z)} = 0.
\end{align}
\section{Velocity in a cylindrical pore}\label{sec:AppendixMeanVelocity}
Since we want to reconstruct the Hagen-Poiseuille and Helmholtz-Smoluchowski approximation we will calculate the volumetric flow $\langle u^0 \rangle = \int_0^L \int_0^R u^0 r \ud r \ud z$.
At first we will integrate $u^0(r,z)$ given by equation \ref{eq:FinalVelocityGeneral} in $z$-direction from $0$ to $L$.
Starting with the first term $\partial_z p_z^0$ yields
\begin{align}
    \frac{1}{4} (r^2-R^2) \int_0^L \partial_z p_z^0 \ud z &= - \frac{1}{4} (r^2-R^2) \Delta p.
\end{align}
For the $K_2$--terms we have
\begin{align}
    \langle u^{\text{EDL}} \rangle_r &= \int_{z=0}^L (K_2(r,z) - K_2(r=R,z)) \ud z \\ \nonumber
    &= \int_0^r \frac{1}{r} \int_0^r \int_{z=0}^L ( (\partial_z p_r^0) r - \Lambda^2 (\partial_z \varphi_r^0)  \partial_r(r \partial_r \varphi_r^0)) \ud z (\ud r)^2  \\ \nonumber
    & \qquad - \int_0^r \frac{1}{r} \int_0^r \int_{z=0}^L ( (\partial_z p_r^0) r - \Lambda^2 (\partial_z \varphi_r^0)  \partial_r(r \partial_r \varphi_r^0)) \ud z (\ud r)^2 |_{r=R}  \\ \nonumber
    &=  - \Lambda^2 \int_{z=0}^L \int_0^r \frac{1}{r} \int_0^r (\partial_z \varphi_r^0)  \partial_r(r \partial_r \varphi_r^0) (\ud r)^2 \ud z \\ \nonumber
    & \qquad + \Lambda^2 \int_{z=0}^L \int_0^r \frac{1}{r} \int_0^r (\partial_z \varphi_r^0)  \partial_r(r \partial_r \varphi_r^0) (\ud r)^2 |_{r=R}  \ud z,
\end{align}
since $\int_{z=0}^L \partial_z p_r^0 \ud z = 0$ due to $p_r^0|_{z=0}=p_r^0|_{z=L}=0$.
For the last term in \eqref{SolutionU} we simply have
\begin{align}
    \frac{\Delta \phi }{L} \Lambda^2 \int_{z=0}^L \left[\phi_r^0 - \phi_r^0 |_{r=R} \right] \ud z,
\end{align}
such that the flow in $z$-direction is
\begin{align}
    \langle u^0\rangle_r(r) = - \Delta p \frac{1}{4} (r^2 - R^2) + \frac{\Delta \phi }{L} \Lambda^2 \int_{z=0}^L  \left[\phi_r^0 - \phi_r^0 |_{r=R} \right] \ud z + \langle u^{\text{EDL}} \rangle_r.
\end{align}
Next, we integrate in $r$-direction from $0$ to $R$
\begin{align}
    \langle u^0\rangle = \int_0^R \langle u^0\rangle_r\,r \ud r = \Delta p \frac{1}{16} R^4 - \frac{\Delta \phi }{L} \Lambda^2 \langle \zeta\rangle + \langle u^{\text{EDL}} \rangle, 
\end{align}
with
\begin{align}
    \langle \zeta\rangle = \int_{z=0}^L \int_{r=0}^R \left[\phi_r^0 - \phi_r^0 |_{r=R} \right] r \ud r \ud z
\end{align}
and
\begin{align}
    \langle u^{\text{EDL}} \rangle &=   \int_{z=0}^L \int_{r=0}^R r \int_0^r \frac{1}{r} \int_0^r (\partial_z \varphi_r^0) (n_+^0 - n_-^0) r (\ud r)^3 \ud z \\ \nonumber
    & \qquad - \frac{1}{2} R^2 \int_{z=0}^L \int_0^r  \frac{1}{r} \int_0^r (\partial_z \varphi_r^0)  (n_+^0 - n_-^0) r (\ud r)^2 |_{r=R} \ud z,
\end{align}
with equation \eqref{eq:FinalPhirGeneral}
\begin{align}
    \langle u^0\rangle =  \Delta p \frac{1}{16} R^4 - \frac{\Delta \phi }{L} \Lambda^2 \langle \zeta\rangle + \langle u^{\text{EDL}} \rangle. 
\end{align}
\section{Current for a cylindrical pore}\label{sec:AppendixCurrent}
We find that in a cylindrical pore, i.e. a straight pipe, the current given by equation \eqref{eq:FinalCurrentGeneral} can further be simplified. We therefore calculate the current at the outlet $z=0$, such that
\begin{align}
    I_\pm |_{z=0} &= \pm 2 \pi \int_0^{R} \left[ \bar n^0 (y_0^0)^{a_\pm} \exp{\mp\phi^0_r} \left\{ k_\pm    \partial_z Q_\pm \pm k_\pm    Q_\pm \frac{\Delta \varphi}{L} -  P_e Q_\pm u^0 \right\} \right]_{z=0} r \ud r.
\end{align}
At the outlet we assume that $\phi_r^0$ is small, such that $\exp{\mp\phi^0_r|_{z=0}} \approx 1 \mp \phi^0_r|_{z=0} \approx 1$.
We further have $Q_\pm \exp{\mp\phi^0_r} |_{z=0} = n_\pm^\tout / \bar n^\tout (y_0^\tout)^{-a_\pm} $ and for the velocity we write
\begin{align}
    \int_0^R u^0 |_{z=0} r \ud r = \frac{\langle u^0\rangle}{L}.
\end{align}
All together we can approximate the current flowing through the outlet by
\begin{align}\label{eq:ApproxCurrent}
    I_\pm |_{z=0} &= \pm 2 \pi \int_0^{R} \left[ \bar n^0 (y_0^0)^{a_\pm} \exp{\mp\phi^0_r} \left\{ k_\pm    \partial_z Q_\pm \pm k_\pm    Q_\pm \frac{\Delta \varphi}{L} -  P_e Q_\pm u^0 \right\} \right]_{z=0} r \ud r \\ \nonumber
    &= \pm 2 \pi \int_0^{R} \bar n^\tout (y_0^\tout)^{a_\pm} \left\{ k_\pm    \partial_z Q_\pm \pm k_\pm \frac{n_\pm^\tout}{\bar n^\tout (y_0^\tout)^{a_\pm}} \frac{\Delta \varphi}{L} -  P_e \frac{n_\pm^\tout}{\bar n^\tout (y_0^\tout)^{a_\pm}}  \frac{\langle u^0 \rangle}{L} \right\} r \ud r \\ \nonumber
    & = \pm 2\pi k_\pm \bar n^\tout (y_0^\tout)^{a_\pm} \frac{1}{2} R^2 \partial_z Q_\pm |_{z=0} 
    \\ \nonumber
    & \quad + 2 \pi n^\tout_\pm \left[ \frac{\Delta \varphi}{L} \left( \frac{k_\pm }{2} R^2 + P_e \Lambda^2 \frac{\langle \zeta\rangle}{L} \right)  \mp  P_e \frac{\Delta p}{L} \frac{1}{16} R^4  \mp P_e \frac{\langle u^{\text{EDL}} \rangle}{L} \right].
\end{align}
\section{Numerical scheme for the asymptotic model}
In order to calculate the asymptotic model we need to solve equation \eqref{eq:FinalPhirGeneral} for $\phi_r^0(r,z)$ for every $z \in [0,L]$ with the respective boundary condition \eqref{eq:FinalPhirBCGeneral}. Once we have the solution for the potential we can solve equation \eqref{eq:FinalQGeneral} for $Q_\pm(z)$ with the boundary condition \eqref{eq:FinalQBCGeneral} to obtain $n_\pm^0(r,z)$. Note that in order to solve \eqref{eq:FinalQGeneral} we also have to calculate the integrals \eqref{eq:FinalQIntegralsGeneral} which also depend on $u^0(r,z)$ and $\phi_z^0(z)$. 
Once we have the solution for the concentrations and the potential we solve equation \eqref{eq:FinalVelocityGeneral} for the velocity $u^0(r,z)$ and equation \eqref{eq:FinalPGeneral} for the pressure $p^0(r,z)$.
This procedure is repeated until the maximum number of iteration $max\_iter$ is reached or the solutions have converged, such that $err < \epsilon$.
\begin{algorithm}
\caption{Numerical scheme to calculate $\phi_r^0$ and $Q_\pm$}\label{alg:AsymptoticModel}
\begin{algorithmic}
\State Set $Q^0_\pm = n_\pm^{\text{init}} \exp{\mp \phi_r^{\text{init}}}$ for $z \in [0,L]$;
\State Set $\phi^0 = \phi^{\text{init}}$, $y_0^0 = y_0^{\text{init}}$, $\bar n^0 = \bar n^{\text{init}}$, $p^0 = p^{\text{init}}$ and $u^0 = u^{\text{init}}$ for $r \in [0,R(z)]$ and $z \in [0,L]$;
\State Calculate $\phi_z(z)$ with equation \eqref{eq:FinalPhizGeneral};
\While{$err > \epsilon$ \textbf{and} $i >  max\_iter$}
    \State Calculate $\phi_r^{i+1}(r,z)$ by solving equation \eqref{eq:FinalPhirGeneral} with equation \eqref{eq:FinalPhirBCGeneral} $\forall z \in [0,L]$;
    \State Calculate integrals \eqref{eq:FinalQIntegralsGeneral} with $\phi_z$, $y_0^{i}$, $\bar n^{i}$, $u^{i}$, $p^{i}$ and $\phi_r^{i+1}$;
    \State Calculate $Q_\pm^{i+1}(z)$ by solving equation \eqref{eq:FinalQGeneral} with \eqref{eq:FinalQIntegralsGeneral} and \eqref{eq:FinalQBCGeneral};
    \State Calculate $y_0^{i+1}$ by solving equation \eqref{eq:Finaly0General};
    \State Calculate $y_\pm^{i+1}(r,z) = Q_\pm^{i+1}(z) \exp{\mp \phi_r^{i+1}(r,z)} (y_0^{i+1})^{a_\pm}$;
    \State Calculate $\bar n^{i+1} = 1/(v_0 + (v_+ - v_0) y_+^{i+1} + (v_- - v_0) y_-^{i+1} ) $;
    \State Calculate $\phi^{i+1}(r,z) = \phi_z(z) + \phi_r^{i+1}(r,z)$;
    \State Calculate $u^{i+1}(r,z)$ with equations \eqref{eq:FinalVelocityGeneral} --\eqref{eq:FinalPotentialGradientGeneral};
    \State Calculate $p^{i+1}(r,z)$ with equations \eqref{eq:FinalPGeneral} -- \eqref{eq:FinalK1General};
    \State $err = max(\| Q_{+}^{i+1} - Q_{+}^{i} \|_{max}, \| Q_{-}^{i+1} - Q_{-}^{i} \|_{max})$;
    \State $i = i+1$; 
\EndWhile
\State Calculate $I_\pm$ using \eqref{eq:FinalCurrentGeneral}.
\end{algorithmic}
\end{algorithm}
We solve the three coupled ordinary differential equations (ODEs) in FEniCS using FEM. The postprocessing and calculation of integrals is done using Python.
\section{Parameter values}
The general physical constants, reference quantities and dimensionless parameters used for all simulations in this work are given in Table \ref{tab:ParameterValuesGeneral}. Values for length, radius and boundary values that where changed for the different examples are given in separate tables. Specific parameters for the trumpet shaped pores in Section \ref{sec:Trumpet} are given in Table \ref{tab:ParameterValuesTrumpet}. The parameters for Section \ref{sec:Cylinder} are given in Table \ref{tab:ParameterValuesCylinder} and the ones for Section \ref{sec:Arbitrary} are given in Table \ref{tab:ParameterValuesArbitrary}.
\begin{table*}[h]%
	\centering %
	\caption{Parameter values used for all simulations.\label{tab:ParameterValuesGeneral}}%
	\begin{tabular*}{\textwidth}{@{\extracolsep\fill}llll@{\extracolsep\fill}}
		\hline
		\textbf{Symbol} & \textbf{Meaning}  & \textbf{Value}  & \textbf{Unit} \\
		\hline \vspace{+1ex}
		$T$ & temperature & 298.15 & \si{\kelvin} \\ \vspace{+1ex}
		$e_0$ & elementary charge & 1.602e-19 & \si{\coulomb} \\  \vspace{+1ex}
		$k_B$ & Boltzmann constant & 1.380e-23 & \si{\joule \per \kelvin} \\  \vspace{+1ex}
		$N_A$ & Avogadro constant & 6.022e23 &  \si{\per \mole}\\ \vspace{+1ex}
		$F = e_0 N_A$ & Faraday constant & 9.648e4 & \si{\coulomb \per \mole} \\  \vspace{+1ex}
		$R_G = k_B N_A$ & Gas constant & 8.314 & \si{\joule \per \mole \per \kelvin} \\  \vspace{+1ex}
		$\eps_0$ & vacuum permittivity & 8.854e-12 & \si{\farad \per \meter}\\\vspace{+1ex}
		$\eps_r = 1+\chi$ & relative permittivity & 78.49  & - \\ \vspace{+1ex}
		$\nu$ & viscosity &  0.8904e-3 & \si{\pascal \second} \\ \vspace{+1ex}
		$c^R$ & ref. concentration & 1 & \si{\mole \per \liter} \\ \vspace{+1ex}
		$q^R$ & ref. charge & 0.16 & \si{\coulomb \per \meter \squared} \\ \vspace{+1ex}
		$D^R$ & ref. diffusion coeff. & 1e-9 & \si{\meter \squared \per \second} \\ \vspace{+1ex}
		$R_0$ & ref. radius &   1e-9    &  \si{\meter}  \\ \vspace{+1ex}
		$\phi^R$ & ref. potential & 0.25 & \si{\volt}  \\ \vspace{+1ex}
		$p^R$ & ref. pressure & 24.8 & \si{\bar}  \\ \vspace{+1ex}
		$\gamma$   & -  &  9.23 & - \\ \vspace{+1ex}
		$\Lambda$  & Debye length  &    0.4  & - \\ \vspace{+1ex}
	\end{tabular*}
\end{table*}
\begin{table*}[h]%
	\centering %
	\caption{Parameter values used for the simulations of EOF transitions in a cylindrical pore.\label{tab:ParameterValuesCylinder}}%
	\begin{tabular*}{\textwidth}{@{\extracolsep\fill}llll@{\extracolsep\fill}}
		\hline 
		\textbf{Symbol} & \textbf{Meaning}  & \textbf{Value}  & \textbf{Unit} \\
		\hline \vspace{+1ex}
		$L_0$    & ref. length       & 1e-8      &\si{\meter} \\ \vspace{+1ex}
        $\delta$    & ratio radius and length        & 0.1      & - \\ \vspace{+1ex}
		$v^R$  &   ref. velocity      & 0.278  & \si{\meter \per \second} \\ \vspace{+1ex}
        $P_e   $      &  Péclet number  &  2.78 & - \\ \vspace{+1ex}
		$I^R  $     &  ref. current   & 9.65e-12  & \si{\ampere} \\ \vspace{+1ex}
	    $Q_{eo}^R $ & ref. flow rate   &   2.78e-19  & \si{\meter \cubed \per \second} \\ \vspace{+1ex}
		$L$ & length & 25 & - \\ \vspace{+1ex}
		$L_1 $, $L_2 $ &length for surface charge & $0.2 L$, $L-L_2$ & - \\ \vspace{+1ex}
		$R$ & radius & 5 & - \\ \vspace{+1ex}
		$k_+$, $k_- $ & diffusion coefficients & 1.33, 0.79   &     -      \\ \vspace{+1ex}
		$n^\tout$ , $n^\tin$& bulk concentrations & 0.6, 0.6    &     -       \\ \vspace{+1ex}
		$\sigma_0  $ & surface charge & 0.15    &   -  \\ \vspace{+1ex}
        $v_0$ & molar volume solvent & 0, 0.018 & - \\ \vspace{+1ex}
        $v_\pm$ & molar volume ions & 0, 0.018, 0.09 & - \\ \vspace{+1ex}
        $a_\pm$ & volume fraction & 0 & - \\ \vspace{+1ex}
	\end{tabular*}
\end{table*}
\begin{table*}[h]%
	\centering %
	\caption{Parameter values used for the simulations of a trumpet shaped pore. \label{tab:ParameterValuesTrumpet}}%
	\begin{tabular*}{\textwidth}{@{\extracolsep\fill}llll@{\extracolsep\fill}}
		\hline
		\textbf{Symbol} & \textbf{Meaning}  & \textbf{Value}  & \textbf{Unit} \\
		\hline \vspace{+1ex}
		$L_0$    & ref. length       & 1e-7      & \si{\meter}\\ \vspace{+1ex}
		$\delta$    & ratio radius and length        & 0.01      & - \\ \vspace{+1ex}
		$v^R$  &   ref. velocity      & 0.0278  & \si{\meter \per \second} \\ \vspace{+1ex}
		$P_e   $      &  Péclet number  &  2.78 & - \\ \vspace{+1ex}
		$I^R  $     &  ref. current   & 9.65e-13  & \si{\ampere} \\ \vspace{+1ex}
		$Q_{eo}^R $ & ref. flow rate   &   2.78e-20  & \si{\meter \cubed \per \second} \\ \vspace{+1ex}
		$L$ & length & 10 & - \\ \vspace{+1ex}
		$L_1 $, $L_2 $ &length for surface charge & $0.1 L$, $L-L_2$ & - \\ \vspace{+1ex}
		$R_1$ & maximum radius &10 & - \\ \vspace{+1ex}
		$R_2$ &minimum radius &1.5 &-  \\ \vspace{+1ex}
		$k_+$, $k_- $ & diffusion coefficients & 1.33, 0.79   &     -      \\ \vspace{+1ex}
		$n^\tout$ , $n^\tin$& bulk concentrations & 0.1, 0.1    &     -       \\ \vspace{+1ex}
		$\sigma_0  $ & surface charge & 1.0          &   -  \\ \vspace{+1ex}
		$\phi^\tout  $, $\phi^\tin $ & bulk potential & 8.0, 0.0   &   -  \\ \vspace{+1ex}
		$p^\tout$ , $p^\tin$ & bulk pressure &0.0, 0.0 & - \\ \vspace{+1ex}
		$v_0$ & molar volume solvent & 0 & - \\ \vspace{+1ex}
		$v_\pm$ & molar volume ions & 0 & - \\ \vspace{+1ex}
		$a_\pm$ & volume fraction & 0, 1, 5 & - \\ \vspace{+1ex}
	\end{tabular*}
\end{table*}
\begin{table*}[h]%
	\centering %
	\caption{Parameter values used for the simulations of a ClyA pore.\label{tab:ParameterValuesArbitrary}}%
	\begin{tabular*}{\textwidth}{@{\extracolsep\fill}llll@{\extracolsep\fill}}
		\hline
		\textbf{Symbol} & \textbf{Meaning}  & \textbf{Value}  & \textbf{Unit} \\ 
		\hline \vspace{+1ex}
		$L_0$    & ref. length       & 1e-8      & \si{\meter} \\ \vspace{+1ex}
        $\delta$    & ratio radius and length        & 0.1      & - \\ \vspace{+1ex}$v^R$  &   ref. velocity      & 0.278  & \si{\meter \per \second} \\ \vspace{+1ex}
        $P_e   $      &  Péclet number  &  2.78 & - \\ \vspace{+1ex}
		$I^R  $     &  ref. current   & 9.65e-12  & \si{\ampere} \\ \vspace{+1ex}
	    $Q_{eo}^R $ & ref. flow rate   &   2.78e-19  & \si{\meter \cubed \per \second} \\ \vspace{+1ex}
		$L$ & length & 5.5 & - \\ \vspace{+1ex}
		$L_1 $, $L_2 $ &length for surface charge & 2.0, 3.5 & - \\ \vspace{+1ex}
		$R_1$ & maximum radius & 5 & - \\ \vspace{+1ex}
		$R_2$ &minimum radius & 1.65 &-  \\ \vspace{+1ex}
		$k_+$, $k_- $ & diffusion coefficients & 1.334, 2.032   &     -      \\ \vspace{+1ex}
		$\sigma_0  $ & surface charge & -0.25   &   -  \\ \vspace{+1ex}
		$p^\tout$ , $p^\tin$ & bulk pressure &0.0, 0.0 & - \\ \vspace{+1ex}
        $v_0$ & molar volume solvent & 0.018 & - \\ \vspace{+1ex}
        $v_\pm$ & molar volume ions & 0.075 & - \\ \vspace{+1ex}
        $a_\pm$ & volume fraction & 4.15 & - \\ \vspace{+1ex}
	\end{tabular*}
\end{table*}

\end{document}